\documentclass[preprint,3p,12pt]{elsarticle}
\usepackage{mathrsfs}
\usepackage{amsmath}
\usepackage{stmaryrd}
\usepackage{bbding}
\usepackage{dcolumn}
\usepackage{graphicx}
\usepackage{amsfonts}
\usepackage{amssymb}
\usepackage{psfrag}
\usepackage{wrapfig}
\usepackage{subfigure}
\usepackage{makeidx}
\usepackage{bm}
\usepackage{epsf}
\usepackage{epsfig}
\usepackage{setspace}
\usepackage{graphicx}
\usepackage{epstopdf}
\usepackage{psfrag}
\usepackage{subfigure}
\usepackage{color}

\begin{document}
\title{A generalized ENO reconstruction in compact GKS for compressible flow simulations}

\author[HKUST1]{Fengxiang Zhao}
\ead{fzhaoac@connect.ust.hk}

\author[HKUST1,HKUST2]{Kun Xu\corref{cor}}
\ead{makxu@ust.hk}

\address[HKUST1]{Department of Mathematics, Hong Kong University of Science and Technology, Clear Water Bay, Kowloon, Hong Kong}
\address[HKUST2]{Shenzhen Research Institute, Hong Kong University of Science and Technology, Shenzhen, China}
\cortext[cor]{Corresponding author}

\begin{abstract}

This paper presents a generalized ENO (GENO)-type nonlinear reconstruction scheme for compressible flow simulations. The proposed reconstruction preserves the accuracy of the linear scheme while maintaining essentially non-oscillatory behavior at discontinuities. By generalizing the adaptive philosophy of ENO schemes, the method employs a smooth path function that directly connects high-order linear reconstruction with a reliable lower-order alternative. This direct adaptive approach significantly simplifies the construction of nonlinear schemes, particularly for very high-order methods on unstructured meshes.
A comparative analysis with various WENO methods demonstrates the reliability and accuracy of the proposed reconstruction, which provides an optimal transition between linear and nonlinear reconstructions across all limiting cases based on stencil smoothness. The consistency and performance of the GENO reconstruction are validated through implementation in both high-order compact gas-kinetic schemes (GKS) and non-compact Riemann-solver-based methods. Benchmark tests confirm the robustness and shock-capturing capabilities of GENO, with particularly superior performance when integrated with compact schemes.
This work advances the construction methodology of nonlinear schemes and establishes ENO-type reconstruction as a mature and practical approach for engineering applications.

\end{abstract}


\maketitle

\section{Introduction}

Compressible flows can simultaneously contain broadband spatial flow structures and strong discontinuities such as shock waves. While high-order schemes offer advantages in accuracy and efficiency compared to second-order schemes, creating a demand for their application in complicated flow simulations \cite{wang2007highorder,wang2021benchmark,wang2022highorder}, high-accuracy discretizations can produce numerical oscillations and computational instability near discontinuities, posing a challenge for high-order numerical schemes. Within the study of numerical schemes for compressible flow simulations, "high-order" typically refers to schemes with a convergence order of two or greater.
Research into high-order numerical schemes for compressible flow simulations has received sustained attention for nearly four decades, resulting in the development of a diverse array of such schemes, such as essentially non-oscillatory (ENO) and weighted ENO (WENO) schemes \cite{harten2,liu-WENO,jiang-WENO}, discontinuous Galerkin (DG) methods \cite{reed,cockburn1,cockburn2}, Flux Reconstruction (FR) and Correction Procedure via Reconstruction (CPR) methods \cite{huynh2007fr,CPR_wang}, and compact gas-kinetic scheme (GKS) \cite{zhao2019compact,zhao_compact-tri,zhao2023direct,zhao2023AIA}.

Nonlinear reconstruction is one of the most crucial foundational components in the construction of high-order schemes, influencing their robustness, accuracy, and dispersive and dissipative properties, particularly for flow simulations involving broadband spatial wavenumber structures.
The objective of nonlinear reconstruction is to implement a special treatment that adaptively switches from a high-order or high-resolution reconstruction to a more reliable low-order reconstruction in the vicinity of discontinuities. Furthermore, for high-order temporal discretization methods developed for time-accurate solvers (such as the GKS solver), it is also necessary to apply an adaptively nonlinear limiting to the high-order linear temporal discretization at strong shock solutions \cite{zhao2023direct}.

ENO schemes are a class of methods used in the computation of compressible flows to achieve high-order accuracy while avoiding spurious oscillations at discontinuities \cite{harten2}. The core principle of ENO is to adaptively select the smoothest stencil from a set of candidates to perform the reconstruction.
Building upon the ENO framework, WENO schemes were subsequently developed \cite{liu-WENO,jiang-WENO,shu-review-1,shu-review-2}.
Through a nonlinear combination of several lower-order candidate polynomials, WENO reconstruction adaptively achieves high-order linear reconstruction in smooth regions and degenerates to a one-sided, lower-order reconstruction that avoids crossing discontinuities in their vicinity.
Initially proposed for structured meshes, WENO reconstruction was later extended to unstructured meshes \cite{hu1999-WENO-tri}, followed by numerous optimization efforts. Optimization research on WENO reconstruction methods primarily focuses on designing nonlinear weights and optimizing the selections of sub-stencils. The representative modified versions of WENO schemes includes WENO-M, WENO-Z, targeted ENO (TENO), central WENO and simplified WENO, etc. \cite{WENO-M,WENO-Z,fu2016teno,fu2023review,levy1999central,capdeville2008central,zhu2016new}.
As WENO schemes and their optimized versions recover linear reconstruction via the nonlinear combination of candidate polynomials from sub-stencils, they exhibit sensitivity to the calculation of nonlinear weights and the selection of sub-stencils. Achieving superior accuracy and resolution often necessitates intricate strategies of designing nonlinear weights. Furthermore, the direct extension of the original WENO algorithm to unstructured meshes poses challenges \cite{zhao2018}, and simple versions that obviate the need for optimal weight calculation offer a more readily applicable form for unstructured meshes \cite{zhu2016new}.
Building upon ENO and WENO reconstruction methods, another notable class of nonlinear reconstruction methods comprises hybrid methods, such as hybrid ENO and hybrid WENO methods \cite{bauer1997hybrid,hill,taylor,ren}. Such hybrid methods typically intend to employ ENO/WENO reconstruction in discontinuous flow regions and linear reconstruction in smooth regions, facilitated by specially designed switching criteria.

This study proposes a generalized ENO (GENO)method for nonlinear high-order reconstruction by adaptively selecting smoothness-based stencils.
Specifically, the GENO reconstruction achieves high accuracy and adaptivity by connecting a high-order linear reconstruction with a lower-order reconstructions with ENO property via a delicately designed continuous path function.
This path function is governed by a linearity-preserving principle, which ensures that the underlying linear reconstruction is recovered for smooth waves, even on the coarse mesh for the under-resolved flow structure.
This property is crucial for turbulent flow simulations characterized by multi-scale eddies.
In addition, compared to WENO reconstruction and its modified versions, the current GENO reconstruction is directly formulated as a hybrid of high-order linear and reliable lower-order reconstructions. This formulation reduces the scheme's dependence on the quality of sub-stencils and their candidate polynomials. Consequently, GENO is particularly well-suited for unstructured meshes, where it can employ simple linear polynomials for sub-stencil reconstructions without compromising accuracy.

In this study, the GENO reconstruction is primarily applied to compact GKS, while also being employed in high-order schemes that utilize approximate Riemann solvers for flux computation.
The compact GKS updates cell-averaged conservative flow variables and their gradients obtained from a time-dependent evolution solution, which are  utilized for compact high-order linear spatial reconstruction. It provides a new framework for constructing compact schemes, emphasizing consistency between the numerical scheme and physical evolution \cite{zhao2023direct}.
The compact GKS achieves genuine spatial-temporal coupling.
Furthermore, the high-order linear reconstruction in compact GKS only involves the immediate neighbors of the reconstruction cell, maintaining the consistency between the sizes of the numerical and the physical domains of dependence. High-order compact GKS has been developed for both structured and unstructured meshes. This study will employ an established linear compact GKS to implement and validate GENO reconstruction.

This paper is organized as follows.
section 2 presents the finite volume method and compact GKS for compressible flow simulations.
The GENO reconstruction will be introduced in Section 3.
Section 4 will present the implementation of high-order schemes on structured and unstructured meshes with the GENO reconstruction.
Section 5 presents an analysis of the GENO method and a comparative study against other reconstruction methods.
In Section 6, the compact GKS will be tested in a wide range of cases. The last section is the conclusion.

\section{Finite volume method and compact GKS}
In the finite volume method for compressible flow simulations, the governing equations for the conservative variables of mass, momentum, and energy, based on conservation laws applied to control volumes defined on the mesh cells, are given as:
\begin{equation}\label{semifvs}
\frac{\text{d}\textbf{W}_{j}}{\text{d}t}=-\frac{1}{\big|\Omega_j\big|}\int_{\partial \Omega_j} \textbf{F}\cdot \textbf{n} \mathrm{d} \Gamma,
\end{equation}
where $\textbf{W}_{j}$ is the cell-averaged conservative variables, $\textbf{F}$ is the numerical fluxes at cell interfaces, $\big|\Omega_j\big|$ is the size of the mesh cell, and $\textbf{n}$ is the unit outer normal vector of the interface $\partial \Omega_j$. Depending on the specific flow being modeled, $\textbf{F}$ can represent either the purely convective fluxes from the Euler equations or the total fluxes from the Navier-Stokes (N-S) equations, which includes both convective and viscous terms.
The cell-averaged conservative variable $\textbf{W}_{j}$ is defined as
\begin{align*}
\textbf{W}_{j}= \frac{1}{\big| \Omega_j \big|} \int_{\Omega_j} \textbf{W}(\mathbf{x}) \mathrm{d}\Omega.
\end{align*}
The integral on the cell interfaces on the right-hand side of Eq.(\ref{semifvs}), which is a line integral in two dimensions and a surface integral in three dimensions, is discretized using Gaussian quadrature as
\begin{align*}
\int_{\partial \Omega_j} \textbf{F}\cdot \textbf{n} \mathrm{d} \Gamma=\sum_{l=1}^{l_0}\big( \big|\Gamma_{l} \big| \sum _{k=1}^{k_0} \omega_k \textbf{F}(\mathbf{x}_k)\cdot \textbf{n}_l \big).
\end{align*}
$l_0$ is the interface number of cell $\Omega_j$, and $\big|\Gamma_{l} \big|$ represents the length or area of the cell interface. $k_0$ and $\omega_k$ are the number of Gaussian quadrature points and weight of the quadrature rule.

To update the conservative variables in Eq.(\ref{semifvs}), the procedure within the finite volume method involves \cite{harten2}: (1) spatial reconstruction to obtain the values (and derivatives) of the conservative variables at the cell interfaces; (2) using the reconstructed values as initial conditions for flow evolution and numerical flux calculation; and (3) temporal discretization.
In this study, the numerical flux is computed using the HLLC approximate Riemann solver \cite{toro2013riemann} and a flux solver based on the gas distribution function of the GKS \cite{xu2,xu1,zhao_compact-tri}, respectively. Correspondingly, the high-order time discretization is performed using a four-stage fourth-order Runge-Kutta (RK) method \cite{liu-WENO} and a nonlinear two-stage fourth-order method (S2O4) \cite{zhao2023direct}, respectively.
The focus of this study is spatial reconstruction. The flux calculation and temporal discretization employed in these methods, which have been published in the literature, will not be presented in this paper.

In the compact GKS employed for this study, the cell-averaged gradients are updated in addition to the cell-averaged conservative variables.
The update is performed using the Gauss theorem:
\begin{equation}\label{slope}
 (\nabla {\bf W})_j^{n+1} = \frac{1}{|\Omega_j|}\int _{\Omega_j} \nabla {\bf W} ({\bf x}, t^{n+1} ) \mathrm{d}\Omega =
 \frac{1}{|\Omega_j|} \int_{\partial \Omega_j} {\bf W}({\bf x},t^{n+1}) {\bf n} \mathrm{d} \Gamma.
\end{equation}
Here, ${\bf W} ({\bf x}, t^{n+1} )$ represents the time-dependent solution within cell $\Omega_j$ at the cell interface, which is directly modeled by using the time-evolving gas distribution function of the GKS, the details of which can be found in \cite{zhao2023direct}.

\section{GENO reconstruction}

In this section, a new high-order nonlinear reconstruction method, designated GENO method, for compressible flow simulations is presented.
Addressing the challenges posed by the coexistence of smooth waves and discontinuities, the GENO method adaptively provides a desirable reconstruction that transitions from a high-order linear polynomial in smooth regions to a lower-order, essentially non-oscillatory reconstruction near discontinuities.
The core idea of the GENO method is to directly connect the high-order linear reconstruction and the low-order essentially non-oscillatory reconstruction via a delicately designed smooth path function.
This function ensures that the GENO method preserves the high-order linear one over a broad range of wavenumbers, particularly when the mesh resolution is insufficient.

The GENO reconstruction is given by
\begin{equation}\label{HLP}
R(\mathbf{x})=\chi q^H(\mathbf{x}) +(1-\chi) q^L(\mathbf{x}),
\end{equation}
where $\chi$ $(0\leq \chi \leq1)$ is the path function, and $\chi$ approaches $0$ at discontinuities and approaches $1$ in smooth regions.
$q^H(x)$ is a high-order polynomial by linear reconstruction, and $q^L(x)$ is a lower-order one with ENO property.
The path function is determined by the smoothness of the high-order linear reconstruction and the specific transition design between the high-order and low-order reconstructions.

\subsection{Path function}

The path function is constructed here.
The transition between a high-order linear reconstruction (optimal for smooth waves) and a lower-order one (for discontinuities) is guided by a principle we term ``linearity preserving". The principle makes the nonlinear reconstruction preserve the high-order linear reconstruction over a broad range of wavenumbers, ensuring accuracy in smooth regions while adaptively switching to a robust, non-oscillatory reconstruction near discontinuities.
The path function $\chi$ is given as
\begin{equation}\label{HLP-weight}
\chi=\mathrm{\mathbf{Tanh}}(C \alpha)/ \mathrm{\mathbf{Tanh}}(C),
\end{equation}
where $\alpha$ is the ultimate smoothness indicator, introduced to quantify the smoothness of the solution within the local reconstruction stencil.
Its detailed formulation is presented later in this section.
$\alpha$ is bounded in the interval $[0,1]$, where the lower and upper bounds correspond to the limiting cases of a sharp discontinuity and a perfectly smooth solution, respectively.
The free parameter $C$ controls the steepness of the hyperbolic tangent function, thereby governing the sharpness of the transition between the high- and lower-order reconstructions.
In this study, $C$ is set to $20$. Numerical tests verified that the results exhibit negligible differences for $C$ values between $10$ and $20$ in the numerical examples considered.

\begin{figure}[!htb]
\centering
\includegraphics[width=0.60\textwidth]{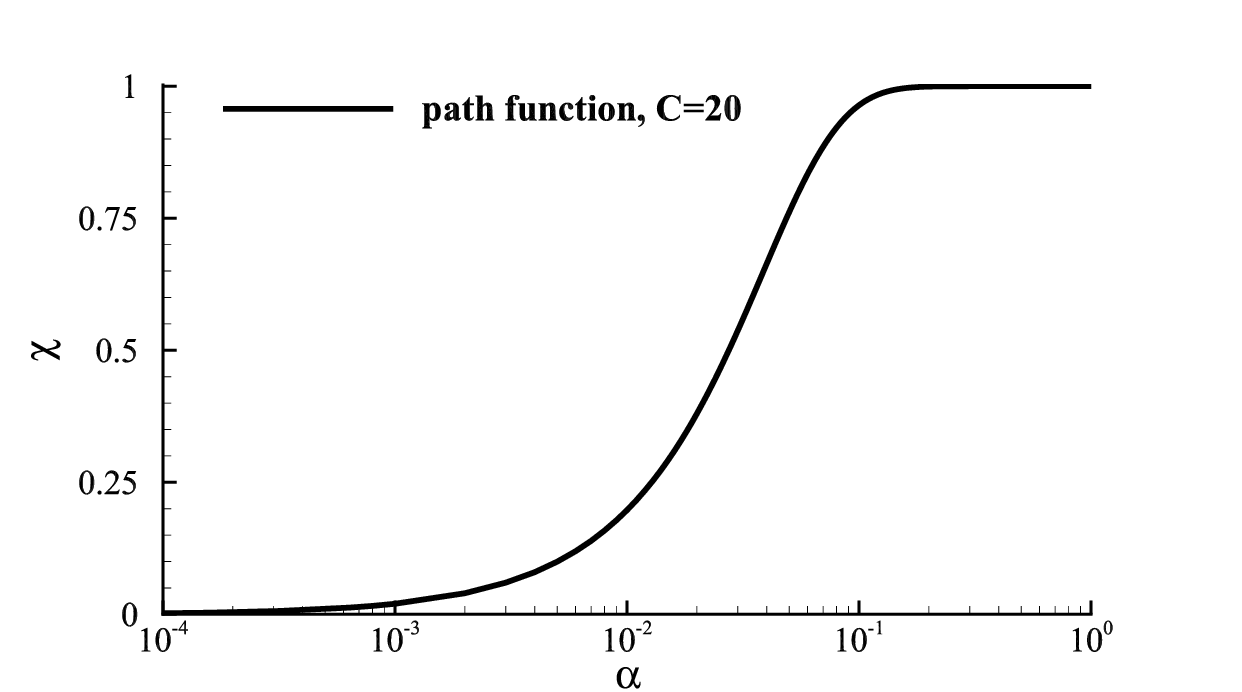}
\caption{\label{HLP-path-function} The path function $\chi$ as a function of the ultimate smoothness indicator $\alpha$. The plot illustrates the transition from the discontinuous limit ($\alpha=0$) to the smooth limit ($\alpha=1$).}
\end{figure}
The path function is plotted in Fig. \ref{HLP-path-function}, with the horizontal axis on a logarithmic scale. The figure illustrates that the path function preserves the high-order linear reconstruction over a wide range of the ultimate smoothness indicator. This linearity-preserving property ensures that the nonlinear scheme preserves to its linear counterpart in simulations with under-resolved flow structures, thereby improving computational accuracy.
Furthermore, the GENO method simplifies the construction of high-order nonlinear schemes, as it only requires a high-order linear reconstruction and a robust low-order one (e.g., a second-order ENO scheme). This feature is particularly beneficial for extending high-order methods in unstructured meshes.

\subsection{Ultimate smoothness indicator}
The ultimate smoothness indicator of the high-order linear reconstruction on the cell is defined as follows:
\begin{equation}\label{HLP-sp}
\begin{split}
&\alpha=\frac{2\alpha^H}{\alpha^H+\alpha^L},\\
&\alpha^H=1+(\frac{IS^{\tau}}{IS^H+\epsilon})^r, \alpha^L=1+(\frac{IS^{\tau}}{IS^L+\epsilon})^r,
\end{split}
\end{equation}
where $\epsilon$ is a small positive constant to avoid division by zero, set to $10^{-15}$.
The value of the exponent $r$ depends on the definitions of $IS^{\tau}$. In this study, $r$ is uniformly set to $2$ and $3$ on structured and unstructured meshes, respectively. These values provide a good balance between accuracy and the ability to handle discontinuities.
$IS^{\tau}$ represents a higher-order smoothness indicator associated with the large stencil of the high-order linear reconstruction.
$IS^H$ and $IS^L$ represent the smoothness indicators of the least smooth and a smooth candidate reconstruction polynomial, respectively, such as a high-order linear reconstruction on the large stencil and a smooth reconstruction on a sub-stencil.
By comparing the smoothness indicators $IS^H$ and $IS^L$, the ultimate smoothness indicator $\alpha$ is determined and takes values between $0$ and $1$.
This formulation was initially employed in \cite{zhao2023direct} to assess the smoothness of flow variables at cell interfaces for nonlinearly limiting high-order flux functions.

Drawing upon the definition of conventional smoothness indicators within WENO schemes \cite{jiang-WENO}, the smoothness indicator of the high-order linear reconstruction $q^H(\mathbf{x})$ in the GENO method is given as:
\begin{equation}\label{IS-highorder-1}
IS^H = \sum_{|l|=1}^{3} h^{2|l|-m} \int_{\Omega} \big( D^{|l|}q^H(\mathbf{x}) \big)^2\mathrm{d}\Omega,
\end{equation}
where $h$ is the mesh cell size, $l$ is a multi-index and $m$ is the dimension, and $D^{|l|}q(\mathbf{x})$ represents the partial derivatives.
For example, in three dimensions ($m=3$), a multi-index of $l=(1,2,0)$ corresponds to $|l|=3$ and $D^{|l|}q(\mathbf{x})=\partial^{3}q(\mathbf{x})/\partial x \partial y^2$.
In the present study, high-order reconstruction is exclusively considered for cases of fourth order and above.
The smoothness indicator in Eq.(\ref{IS-highorder-1}) includes only terms up to the third-order derivative, thereby excluding terms associated with higher-order derivatives.
The rationale for this formulation of Eq.(\ref{IS-highorder-1}) is to ensure that the smoothness indicators $IS^H$ and $IS^L$ remain of comparable magnitude, allowing the resulting ultimate smoothness indicator $\alpha$ in Eq.(\ref{HLP-sp}) to provide a balanced and reliable smoothness measurement. This design choice has been empirically validated through numerical tests, demonstrating enhanced accuracy and robustness.

On structured meshes, when all polynomials on the sub-stencils are of quadratic or higher degree, the determination of $IS^H$ can be further simplified as follows:
\begin{equation}\label{IS-highorder-2}
IS^H = \mathrm{\mathbf{Max}}\{IS_k\}, ~k=1,2,\cdots,
\end{equation}
where $IS_k$ are the smooth indicators of the polynomials on sub-stencils used for the lower-order ENO-type reconstruction.
Numerical examples validate that this simplification maintains both accuracy and robustness.

\subsection{Lower-order reconstructions with the ENO property}

In the GENO reconstruction, the lower-order reconstruction with ENO property adaptively dominates in regions of solution discontinuities.
In this study, the following approach is adopted:
\begin{equation}\label{HLP-q-loworder}
\begin{split}
&q^L(\mathbf{x})=\sum_{k=1}^{k_0} \omega_k p_k(\mathbf{x}), \\
&\omega_k=\frac{\widetilde{\omega}_k}{\sum_k\widetilde{\omega}_k}, \widetilde{\omega}_k=\frac{d_k}{(IS_k+\epsilon)^{r^L}},
\end{split}
\end{equation}
where exponent $r^L$ is set to $2$.
$p_k(\mathbf{x})$ represent the reconstruction polynomials on the sub-stencils and $k_0$ is the total number of sub-stencils.
$d_k$ are linear weights for optimizing the low-order reconstruction with the ENO property.
For simplicity, the weights $d_k$ can be uniformly set to $1$ for all $p_k$.
However, an improved strategy is used when a central sub-stencil is available.
Assuming the sub-stencil corresponding to $p_1$ is central, its weight is boosted by setting $d_1 = C_0$, while the other weights remain $d_k=1$ for $k>1$.
A value of $C_0< 10$ is recommended; in this study, we use $C_0 =8$.
At discontinuities of the solution, $q^L(\mathbf{x})$ plays a dominant role in the reconstruction given by Eq.(\ref{HLP}).
The nonlinear combination in Eq.(\ref{HLP-q-loworder}) ensures that the smoothest one or several of the $p_k(\mathbf{x})$ will dominate the reconstruction of $q^L(\mathbf{x})$, thereby guaranteeing the ENO property of the GENO method.


$IS^L$ in Eq.(\ref{HLP-sp}) is obtained from the low-order, ENO-type reconstruction $q^L$.
On a structured mesh, a choice for $IS^L$ that allows for the effective identification of discontinuities is given by:
\begin{equation}\label{HLP-IS-loworder-2}
IS^L=\mathrm{\mathbf{Min}}\{IS_k\}.
\end{equation}
Due to the simplicity of the reconstruction on structured meshes, a polynomial of at least second order can be used to determine $q^L$, which will be detailed in the next section. Consequently, the $IS^L$ formulation above achieves a good balance between accuracy and discontinuity-capturing capability.

However, for unstructured meshes, given the higher complexity of multi-dimensional high-order reconstruction, we recommend using a first-order polynomial to define $q^L$.
Since the smoothness indicator corresponding to a first-order polynomial is of low accuracy, an $IS^L$ with improved accuracy is given by:
\begin{equation}\label{HLP-IS-loworder-3}
IS^{L}=\big(\sum_{k=1}^{k_0}IS_k -\mathrm{\mathbf{Min}}\{IS_k\} -\mathrm{\mathbf{Max}}\{IS_k\} \big)/(k_0-2).
\end{equation}

\section{High-order schemes based on GENO reconstruction}

In this section, the specific implementation of GENO reconstruction within several existing high-order schemes will be presented.
The GENO reconstruction will be used in these high-order schemes, including high-order compact GKS and non-compact finite volume schemes employing approximate Riemann solvers in one dimension, as well as high-order compact GKS on two- and three-dimensional unstructured meshes.

\subsection{High-order non-compact scheme in one dimension}
A 6th-order non-compact WENO scheme was developed for structured meshes in \cite{zhao2017}. The scheme utilizes a six-cell stencil, symmetric about the cell interface, to achieve 6th-order reconstruction.
This subsection details the GENO reconstruction for the left-side interface value in the 6th-order reconstruction. The symmetric right-side reconstruction is readily obtained and thus omitted here.
The large stencil $S_0$ and sub-stencils $S_k$ $(k=1,2,3,4)$ used for the reconstruction are as follows:
\begin{equation*}
\begin{split}
&S_0=\{I_{j-2},I_{j-1},I_{j},I_{j+1},I_{j+2},I_{j+3}\},\\
&S_1=\{I_{j-2},I_{j-1},I_{j}\},\\
&S_2=\{I_{j-1},I_{j},I_{j+1}\},\\
&S_3=\{I_{j-1},I_{j},I_{j+1},I_{j+2}\},\\
&S_4=\{I_{j},I_{j+1},I_{j+2},I_{j+3}\}.
\end{split}
\end{equation*}
With this sub-stencils selection, the 6th-order WENO reconstruction was presented in \cite{zhao2017}.

The GENO reconstruction uses the stencil $S_0$ for a high-order linear reconstruction $p^5$, and a nonlinear combination of polynomials $p_k$ on $S_k$ $(k=1,2,3,4)$ for lower-order reconstructions with the ENO property. Therefore, the 6th-order GENO reconstruction is
\begin{equation}\label{6th-noncompact-HLP}
\begin{split}
&R(x)= \chi p^5(x) +(1-\chi)\sum_{k=1}^4\omega_kp_k(x),\\
\end{split}
\end{equation}
The specific expressions for $p^5(x)$ and $p_k(x)$ are detailed in \cite{zhao2017}.
In addition, $IS^{\tau}$ is given as
\begin{equation}\label{6th-noncompact-path}
IS^{\tau}=|(IS_1+IS_3)/2-IS_2|.
\end{equation}

Furthermore, the 5th-order scheme with GENO reconstruction was also tested in numerical examples of the current study.
The high-order linear reconstruction and the lower-order reconstruction within the GENO method were obtained from the large stencil and sub-stencils, respectively, of the 5th-order WENO scheme \cite{jiang-WENO}. Details are omitted here for brevity.
The ultimate smoothness indicator for the 5th-order GENO scheme utilizes the same definition as Eq.(\ref{6th-noncompact-path}).
While $IS_1$, $IS_2$, and $IS_3$ represent the smoothness indicators for the polynomials on the three sub-stencils of the 5th-order WENO reconstruction, respectively.

\subsection{High-order compact GKS in one dimension}
The high-order compact GKS for structured meshes was presented in \cite{zhao2023direct}, which employs a simplified WENO (sWENO) reconstruction avoiding the computation of optimal linear weights.
The 8th-order compact GKS using GENO reconstruction is presented in this study.
The high-order linear reconstruction and the lower-order part of the GENO reconstruction are given by the following stencils:
\begin{equation*}
\begin{split}
&S_0=\{I_{j-1},I_{j},I_{j+1},I_{j+2}\},\\
&S_1=\{I_{j-1},I_{j}\},\\
&S_2=\{I_{j-1},I_{j},I_{j+1}\},\\
&S_3=\{I_{j},I_{j+1},I_{j+2}\}.
\end{split}
\end{equation*}
It should be emphasized that only three linear polynomials are used here to obtain the lower-order reconstruction in GENO, without requiring the reconstruction from all sub-stencils as in \cite{zhao2023direct}, thus simplifying the algorithm.
In the numerical examples presented in this study, the sub-stencils for sWENO and TENO were kept exactly the same as in \cite{zhao2023direct}.
The cell average and its gradient on each cell within the $S_0$ stencil are used for the 8th-order reconstruction.
The cell average and its gradient on cell $I_{j-1}$ within the $S_1$ stencil are used for the reconstruction.
Three quadratic polynomials are obtained from $S_k$ $(k=1,2,3)$.
Thus the high-order GENO reconstruction is obtained as
\begin{equation}\label{8th-compact-HLP}
\begin{split}
&R(x)= \chi p^7(x) +(1-\chi)\sum_{k=1}^3\omega_kp_k(x).\\
\end{split}
\end{equation}
And $IS^{\tau}$ is given as follows:
\begin{equation}\label{6th-compact-path}
IS^{\tau}=|(IS_2+IS_3)/2-IS_1|.
\end{equation}

\subsection{High-order compact GKS on unstructured meshes}

The compact GKS for 2-D and 3-D unstructured meshes have been developed in \cite{zhao2023direct,zhao2023AIA}. The nonlinear reconstructions of these schemes are based on the sWENO approach.
Since high-order reconstructions on unstructured meshes are more computationally and memory intensive than on structured meshes, these compact GKS employ a nonlinear combination of linear polynomials and a higher-order polynomial for the nonlinear reconstruction.
However, the limited accuracy of smoothness indicators for linear polynomials, particularly in the vicinity of extrema where they tend towards zero, can result in the non-linear reconstruction being unduly influenced by the linear components, leading to a reduction in overall accuracy.
The GENO method directly combines the high-order linear reconstruction with lower-order reconstructions, avoiding explicit smoothness comparisons and polynomial assembly. This approach fundamentally mitigates the accuracy degradation caused by the use of very low-order polynomials on the sub-stencils.

\begin{figure}[!htb]
	\centering
    \includegraphics[width=0.30\textwidth]{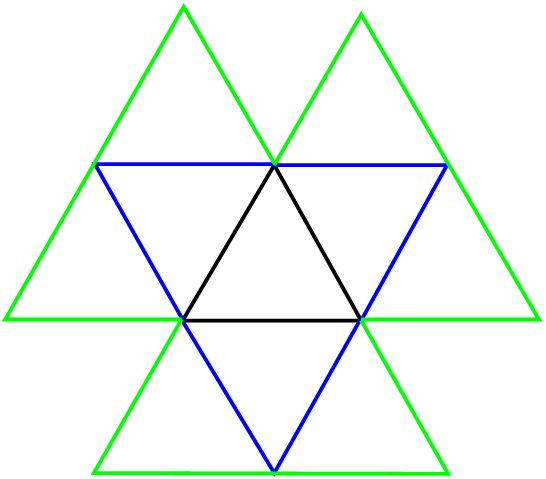}
    \hspace{2.0cm} 
	\includegraphics[width=0.30\textwidth]{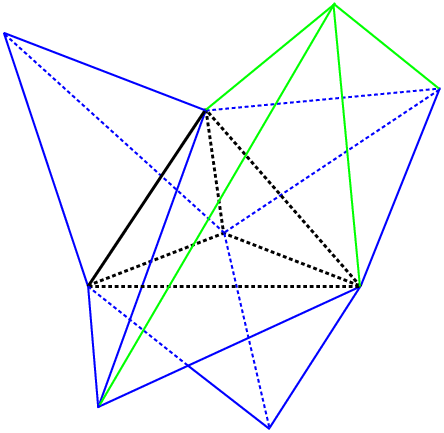}
	\caption{\label{schematic-compact-stencil} A schematic of reconstruction stencils of 4th-order compact GKS. For the sake of simplicity and clarity of the illustration in 3-D case, only two cells of the second-level neighboring cells (green) are shown. }
\end{figure}

This study considers two-dimensional triangular and three-dimensional tetrahedral unstructured meshes. The reconstruction cell is denoted as $\Omega_j$, the face-adjacent neighboring cells are denoted as $\Omega_{j_k}$, and the compact reconstruction stencil is denoted as $S_0$ which gives the high-order linear reconstructions, as illustrated in Fig. \ref{schematic-compact-stencil}.
While the GENO applies to both 2D and 3D, this subsection will focus on the 3D case, as the 2D implementation is formally analogous.
The low-order part of the GENO reconstruction is obtained by the following four sub-stencils:
\begin{equation*}
\begin{split}
&S_1=\{\Omega_{j},\Omega_{j_1}\},\\
&S_2=\{\Omega_{j},\Omega_{j_2}\},\\
&S_3=\{\Omega_{j},\Omega_{j_3}\},\\
&S_4=\{\Omega_{j},\Omega_{j_4}\}.
\end{split}
\end{equation*}
The cell averages and the gradient on cell $\Omega_{j_k}$ within the $S_k$ stencil are used for the reconstruction.
Four linear polynomials are obtained from $S_k$ $(k=1,2,3,4)$.
The high-order GENO reconstruction is given by:
\begin{equation}\label{6th-compact-HLP}
\begin{split}
&R(x)= \chi p^3(x) +(1-\chi)\sum_{k=1}^4\omega_kp_k(x).\\
\end{split}
\end{equation}

Furthermore, an auxiliary smoothness indicator, denoted as $\widetilde{IS}^H$, is introduced for the purpose of determining $IS^{\tau}$. This indicator, $\widetilde{IS}^H$, is defined by
\begin{equation}\label{IS-highorder-2}
\widetilde{IS}^H = \sum_{|l|=1}^{2} h^{2|l|-m} \int_{\Omega} \big( D^{|l|}q^H(\mathbf{x}) \big)^2\mathrm{d}\Omega,
\end{equation}
where the calculation of $\widetilde{IS}^H$ exclusively includes terms corresponding to the first and second-order derivatives.
Thus, a higher-order smoothness indicator $IS^{\tau}$ is constructed and given as:
\begin{equation}\label{IS-tau}
IS^{\tau}= |IS^H-\widetilde{IS}^H|.
\end{equation}
The formulation of $IS^{\tau}$ as defined above is advantageous for unstructured meshes because $IS^H$ and $\widetilde{IS}^H$ are directly derived from the high-order linear reconstruction, without requiring contributions from lower-order sub-stencils. At smooth extrema, $IS^{\tau}$ is a higher-order small quantity relative to $IS^H$ and $IS^L$, which enhances computational accuracy in flows involving wave propagation and vortex evolution. A similar construction methodology for $IS^{\tau}$ was initially presented in \cite{zhao2023AIA}.

\section{Comparison of GENO method with other nonlinear methods}

\subsection{Overview for WENO and TENO reconstructions}
The successful treatment of discontinuities by WENO methods and other high-order non-linear schemes hinges on the adaptive procedure. This procedure nonlinearly combines several lower-order polynomials to effectively recover a high-order linear reconstruction or achieve a reliable low-order upwind reconstruction, where the design of the combination weights is the critical component.
For smooth solutions, preserving the prescribed high-order linear reconstruction on the large stencil is the optimal objective for non-linear schemes such as WENO, enabling the attainment of optimal accuracy, dispersion, and dissipation properties.

This subsection reviews the WENO scheme and its variants in the one dimension. The conventional high-order WENO scheme is given as
\begin{equation}\label{5th-nonlinear}
R(x_{j+1/2})=\sum_{k=1}^{k_0}w_k q_k(x_{j+1/2}),
\end{equation}
where $x_{j+1/2}$ is the cell interface, $q_k(x)$ are the quadratic polynomials on the candidate sub-stencils.
The WENO-Z weight is designed as
\begin{equation}\label{wgt-WENO-Z}
\begin{split}
w_k&=\frac{\widetilde{w}_k}{\sum_{k=1}^{k_0}\widetilde{w}_k}, \\
\widetilde{w}_k&=d_k\big(1+ (\frac{\tau_Z}{IS_k+\epsilon}) \big),
\end{split}
\end{equation}
where $d_k$ are the optimal linear weights.

To optimize WENO schemes such that they maintain accuracy close to, or even comparable to, linear schemes when the mesh resolution is insufficient to fully resolve smooth flow structures, significant research has been devoted to optimize the nonlinear weights.
A notable example of this research is the TENO scheme \cite{fu2023review}. The nonlinear weights in TENO are obtained by selecting either $0$ or the optimal linear weight $d_k$. For the nonlinear reconstruction in Eq.(\ref{5th-nonlinear}), the TENO weights are defined as follows.
First, the scale-separation indicator for the TENO scheme is defined as
\begin{equation}\label{wgt-TENO-1}
\begin{split}
\gamma_k&=\frac{\widetilde{\gamma}_k}{\sum_{k=1}^{k_0}\widetilde{\gamma}_k}, \\
\widetilde{\gamma}_k&=(\frac{1}{IS_k+\epsilon})^7,
\end{split}
\end{equation}
Based on the reference \cite{fu2023review} and numerical examples in this study, the exponent in Eq.(\ref{wgt-TENO-1}) is set to $7$.
The nonlinear weights, characterized by the cutting process, are determined by the scale-separation indicator as follows:
\begin{equation}\label{wgt-TENO-2}
\begin{split}
&w_k=\frac{\widetilde{w}_k d_k}{\sum_{k=1}^{k_0}\widetilde{w}_k d_k}, \\
&\widetilde{w}_k=\left\{\begin{array}{ll}
0,  &  \gamma_k \leq C_T,\\
1,  &  \mathrm{else}.
\end{array} \right.
\end{split}
\end{equation}
The parameter $C_T$ is the threshold used for cutting. Based on the reference \cite{fu2023review} and numerical experiments in this study, the value of $C_T$ is set to $10^{-6}$.
A notable exception is the Titarev-Toro test case in Section 6.3, where the use of TENO weights with the eighth-order compact scheme requires $C_T$ to be set to $10^{-3}$.
Smaller values of $C_T$ lead to significant numerical oscillations.

\subsection{Comparison study of the nonlinear reconstructions}

To conduct a comparative study of WENO, TENO, and GENO, the proportion of linear high-order reconstruction within the nonlinear reconstruction of each method is presented, with this proportion varying according to the solution's smoothness on the reconstruction stencil.
For the GENO method, the path function $\chi$ serves as a parameter to represent the proportion of the linear reconstruction. At $\chi=1$, the nonlinear reconstruction is entirely determined by the linear one.
For Eq.(\ref{5th-nonlinear}), corresponding to the WENO and TENO reconstructions, we define a parameter to represent the extent to which the linear reconstruction is preserved
\begin{equation}\label{linearity-paprameter}
\chi=\mathrm{\mathbf{Min}}(w_k/d_k),~k=1,\cdots,k_0.
\end{equation}
In the TENO reconstruction, the values of $w_k$ are either $d_k$ or $0$. Therefore, as long as the reconstruction on at least one sub-stencil is non-smooth, the TENO reconstruction is entirely dominated by the smooth sub-stencils, resulting in the maximum deviation from the linear high-order reconstruction.
In the WENO reconstruction, the values of $w_k$ continuously transition from $d_k$ to $0$. Only when $w_k=0$ does the reconstruction deviate maximally from the linear reconstruction.

The following two scenarios with different smoothness of the solution on the reconstruction stencil will be considered.
\begin{enumerate}
    \item $IS^{max}=\phi \cdot IS^{min}$, $\tau_Z\sim IS^{max}$,
    \item $IS^{max}=\phi \cdot IS^{min}$, $\tau_Z\sim IS^{min}$.
\end{enumerate}
The parameter $\phi$ is the ratio of smoothness indicators for the least and most smooth candidate polynomials used for nonlinear reconstructions.
The two scenarios represent extreme cases that pose challenges to nonlinear reconstructions.
Near discontinuities, the values of smoothness indicators on candidate stencils spanning the discontinuity become relatively large.
Near extrema, the value of the smoothness indicator on the smoothest sub-stencil approaches zero. The difference of theses two scenarios is reflected in the parameter $\tau_Z$.

\begin{figure}[!htb]
	\centering
    \includegraphics[width=0.495\textwidth]{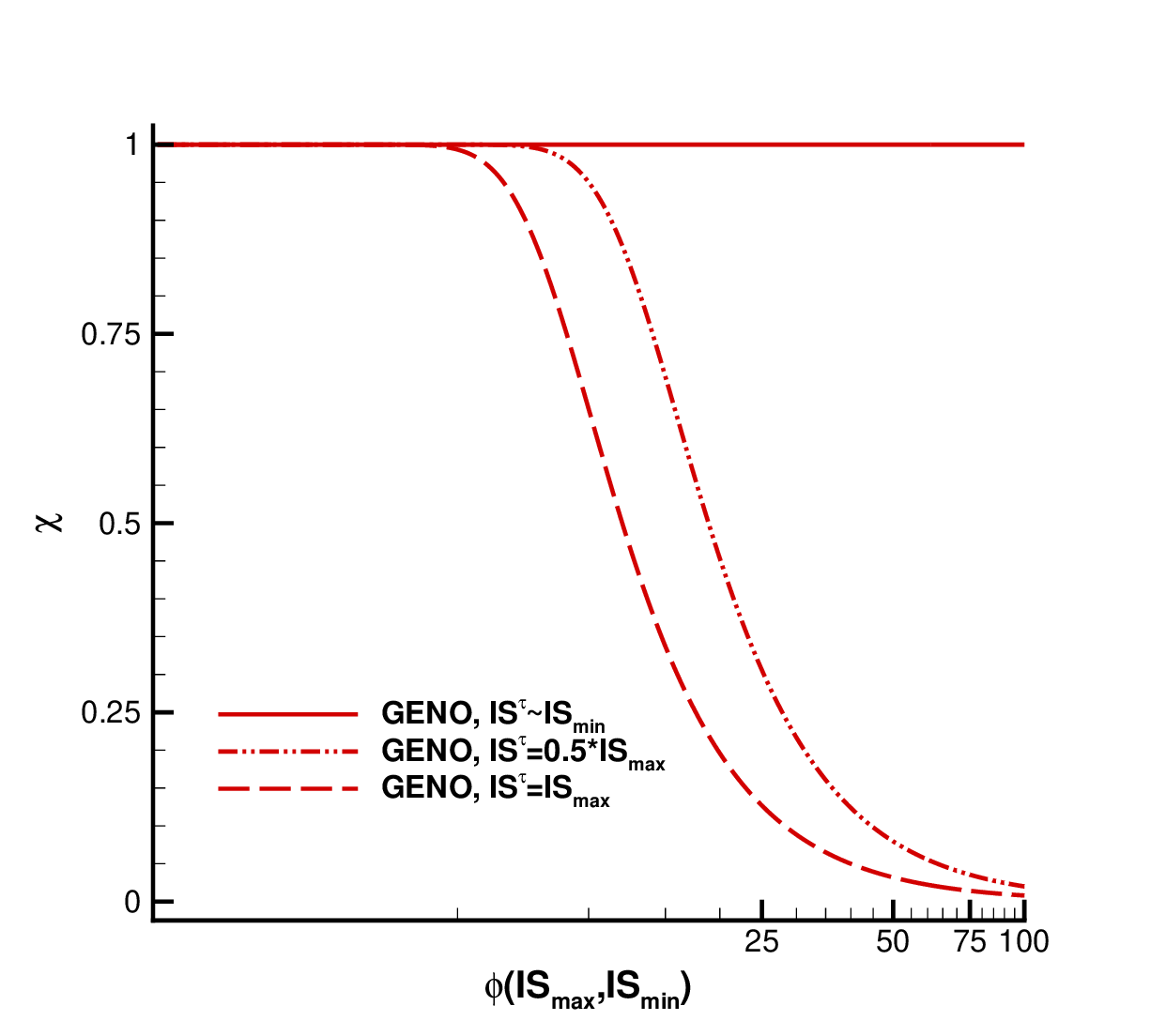}
	\caption{\label{Comparison-nonlinear-1} The proportion of the linear reconstruction part in the GENO reconstruction method for solutions with different levels of smoothness. $\phi$ is the ratio of smoothness indicators for the least and most smooth candidate polynomials. }
\end{figure}

\begin{figure}[!htb]
	\centering
    \includegraphics[width=0.495\textwidth]{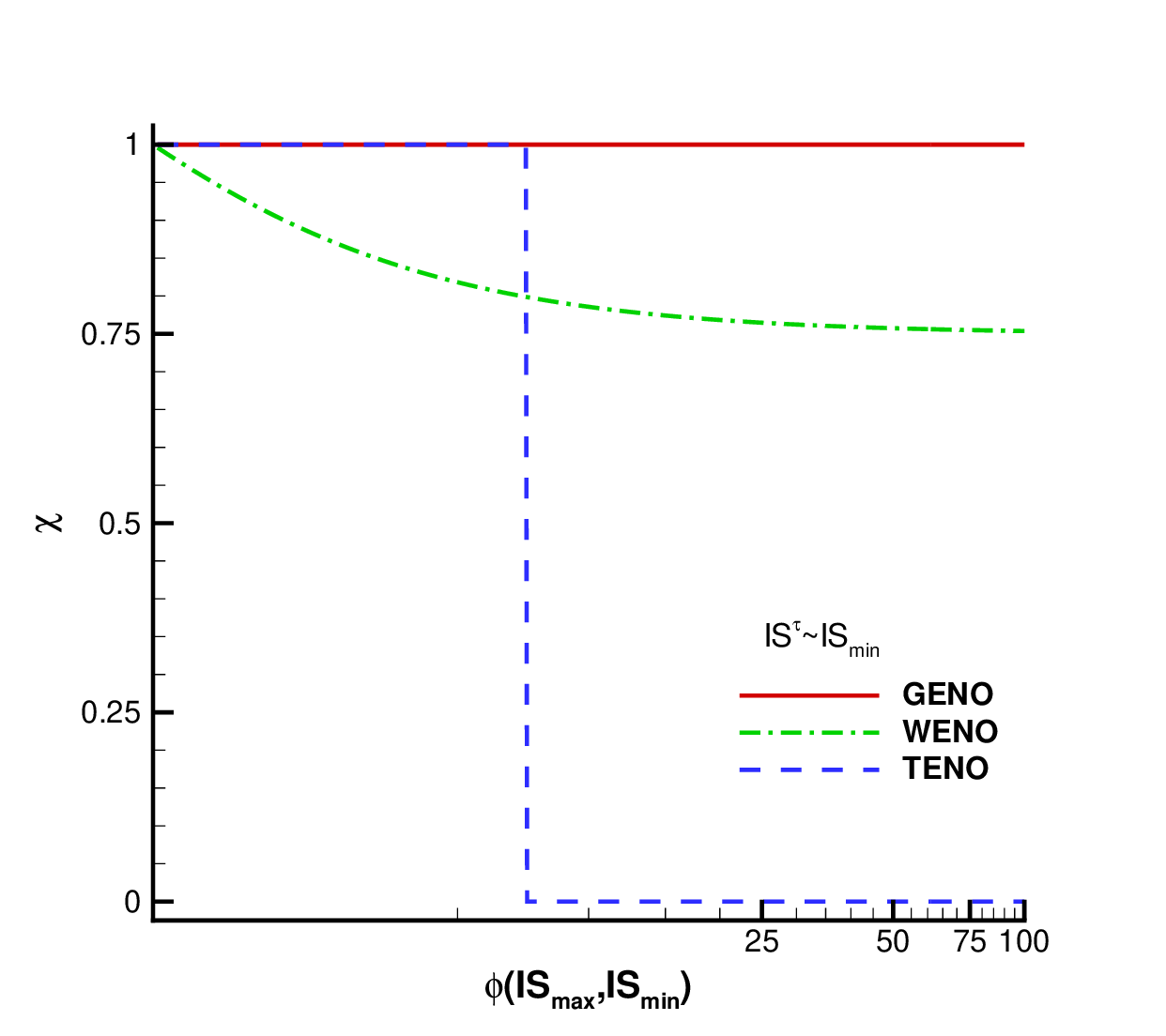}
    \includegraphics[width=0.495\textwidth]{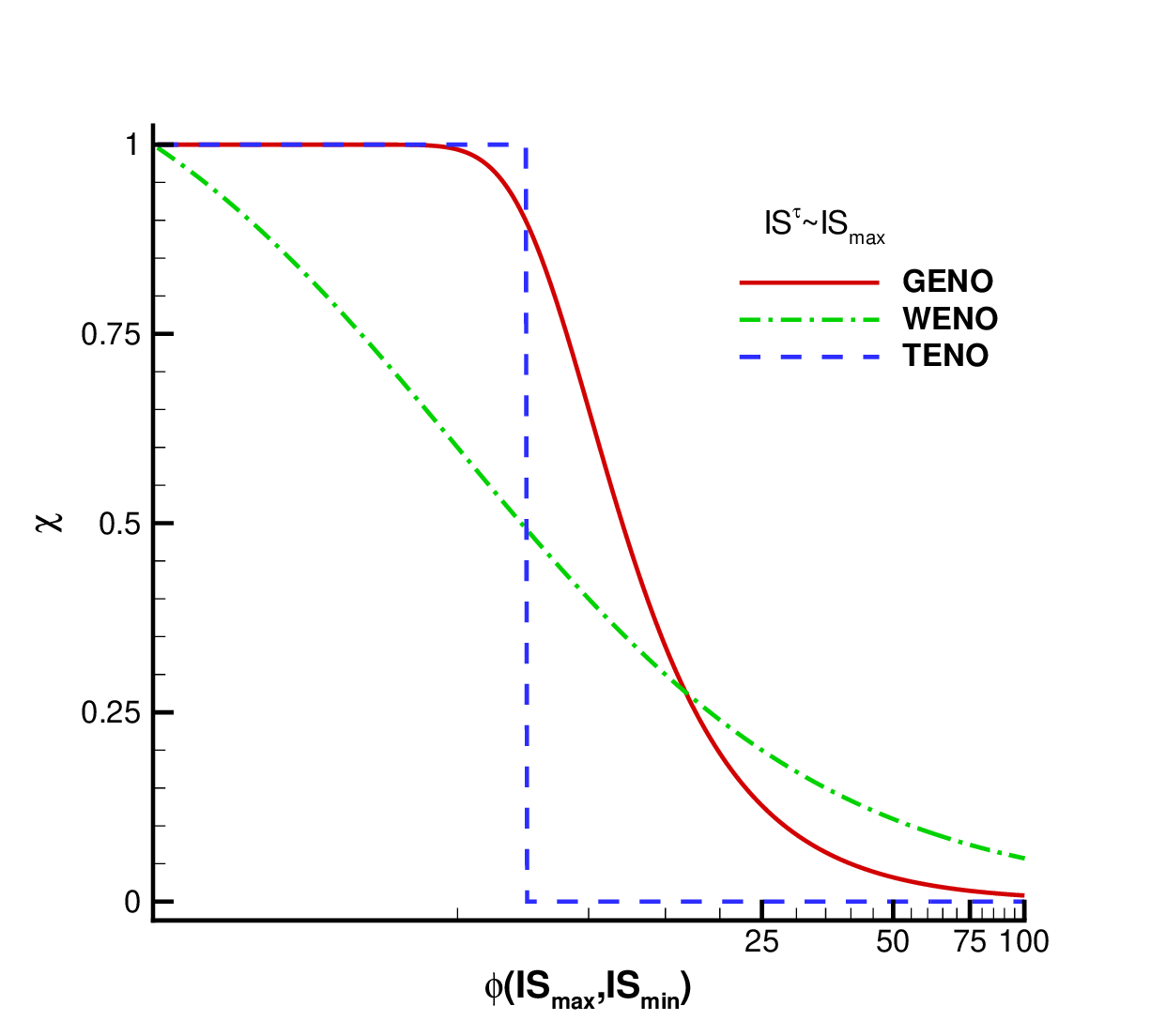}
	\caption{\label{Comparison-nonlinear-2} The proportion of the linear reconstruction part in the conventional WENO, TENO and the GENO reconstruction methods for solutions with different levels of smoothness. Left: $IS^{\tau} \sim IS_{min}$ ;Right: $IS^{\tau} \sim IS_{max}$.}
\end{figure}

Fig. \ref{Comparison-nonlinear-1} illustrates the proportion of the linear reconstruction in the final reconstruction for the GENO method under different $IS^{\tau}$.
The solid and dashed lines in Fig. \ref{Comparison-nonlinear-1} represent the proportion of the high-order linear reconstruction for smooth and discontinuous solutions, respectively, and curves for any other complex scenarios lie between these two.
A key feature of the GENO reconstruction is clearly demonstrated by the solid line in Fig. \ref{Comparison-nonlinear-1}, which perfectly handles the extreme case of smooth extrema by preserving the high-order linear reconstruction. 
For WENO and TENO reconstructions, the deviation from linear reconstruction in the two above extreme cases is shown in Fig. \ref{Comparison-nonlinear-2}.
For quantitative analysis, it assumes a simplified yet illustrative scenario employing three sub-stencils, with the linear weight assumed to be $d_k = 1/3$.
The performance of GENO is also presented for comparison.
In the case of smooth extrema, as shown in the left of Fig. \ref{Comparison-nonlinear-2}, the TENO method switches the reconstruction from a high-order linear reconstruction to a low-order reconstruction with ENO properties when $\phi$ reaches a certain value.
The WENO method deviates slightly from the high-order linear reconstruction even when $\phi$ is small. GENO reconstruction, however, maintains the optimal linear reconstruction throughout.
In the case of discontinuities, corresponding to the right of Fig. \ref{Comparison-nonlinear-2}, all three reconstructions eventually deviate completely from linear reconstruction, although the transition paths differ. For higher mesh resolution corresponding to smaller $\phi$ values, both GENO and TENO maintain the target linear reconstruction.

\section{Numerical examples}
This section presents benchmark computations using compact and non-compact high-order schemes with GENO reconstruction.
For comparison, results from corresponding schemes equipped with WENO (or sWENO) and TENO reconstructions are also presented.
The high-order schemes using the WENO (or sWENO) reconstruction have been detailed in the published literatures.
In contrast, the TENO-based reconstructions are implemented herein by replacing the calculation of the nonlinear weights with the formulation of TENO weights.
Specifically, the TENO weights for the fifth-order scheme are taken from \cite{fu2016teno}, whereas those for the sixth- and eighth-order schemes employed in this study are given by Eq.(\ref{wgt-TENO-2}).

First, one-dimensional accuracy test cases involving smooth, high-wavenumber wave propagation are presented. In addition, shock-entropy and shock-shock interactions are simulated using a high-order compact GKS and non-compact high-order finite volume schemes (with HLLC flux).
Subsequently, GENO reconstruction is applied within a high-order compact GKS on two- and three-dimensional unstructured meshes.
Two-dimensional examples, including hypersonic flow past a cylinder, flow over a forward-facing step, and a viscous shock tube problem, are then presented.
Following this, an initial three-dimensional accuracy test verifies GENO's linearity-preserving property on unstructured meshes. Finally, other three-dimensional examples, encompassing Taylor-Green vortex flow and explosion and implosion scenarios, are computed to validate the accuracy and robustness of GENO reconstruction.

The time step is determined by the CFL condition with $CFL\geq 0.5$ in all test cases.
For viscous flow, the time step is limited by the viscous term as $\Delta t=CFL \cdot h^2/(2(m-1) \nu)$ as well, where $m$ denotes the dimension, $\nu$ is the kinematic viscosity coefficient, and $h$ is the cell size.
For one-dimensional uniform mesh, $h=\Delta x$; and for the two- and three-dimensional meshes, $h_j=m|\Omega_j |/|\Gamma_j|$ which is the radius of the inscribed circle or sphere in the control volume $\Omega_j$.

\subsection{1-D accuracy test}
The one-dimensional advection of density perturbation is tested first. The initial condition is given as follows
\begin{align*}
\rho(x)=1+0.2\sin(\pi x),\ \  U(x)=1,\ \ \  p(x)=1, x\in[0,2].
\end{align*}
The periodic boundary condition is adopted.

\begin{figure}[!htb]
\centering
\includegraphics[width=0.495\textwidth]{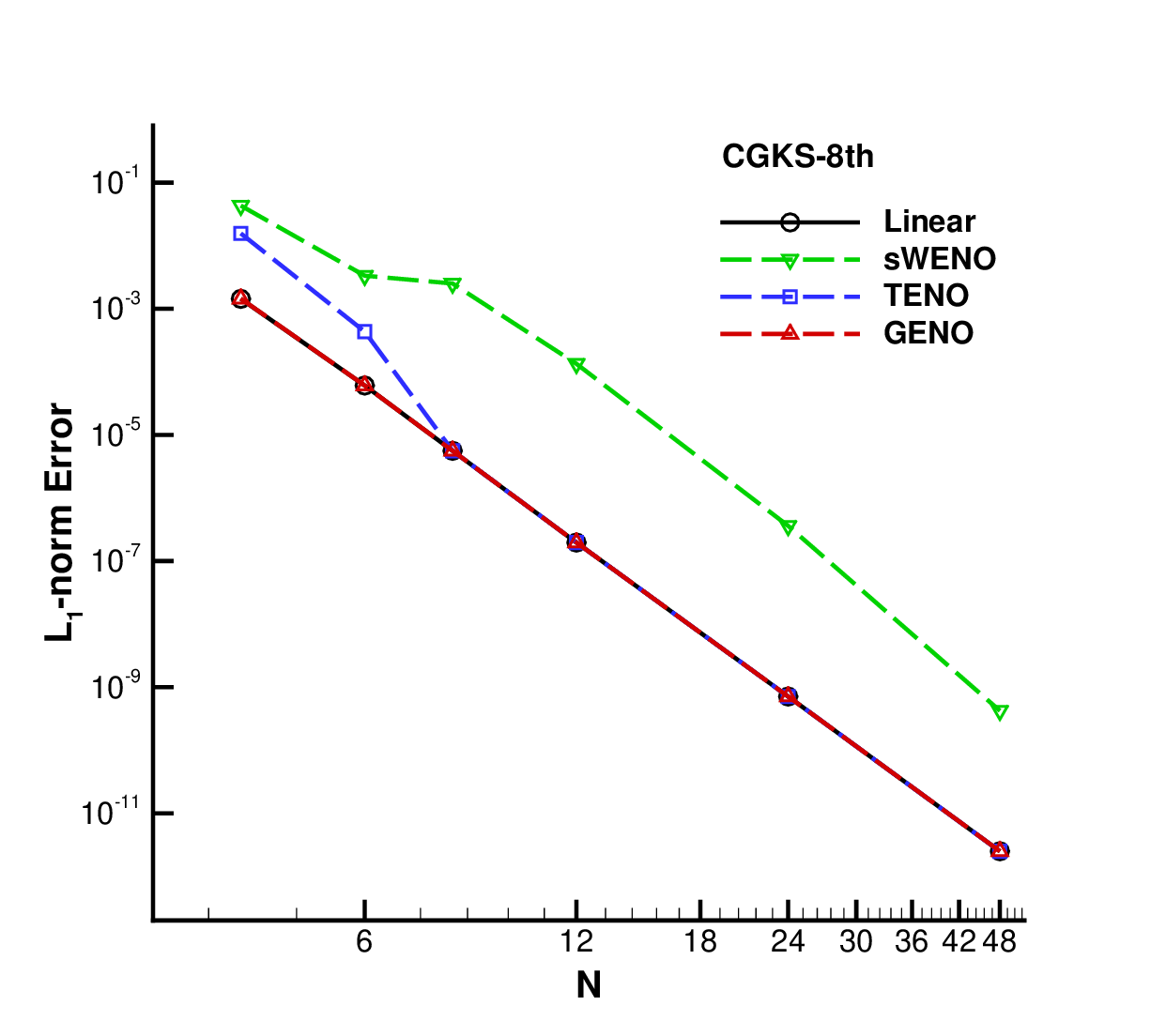}
\includegraphics[width=0.495\textwidth]{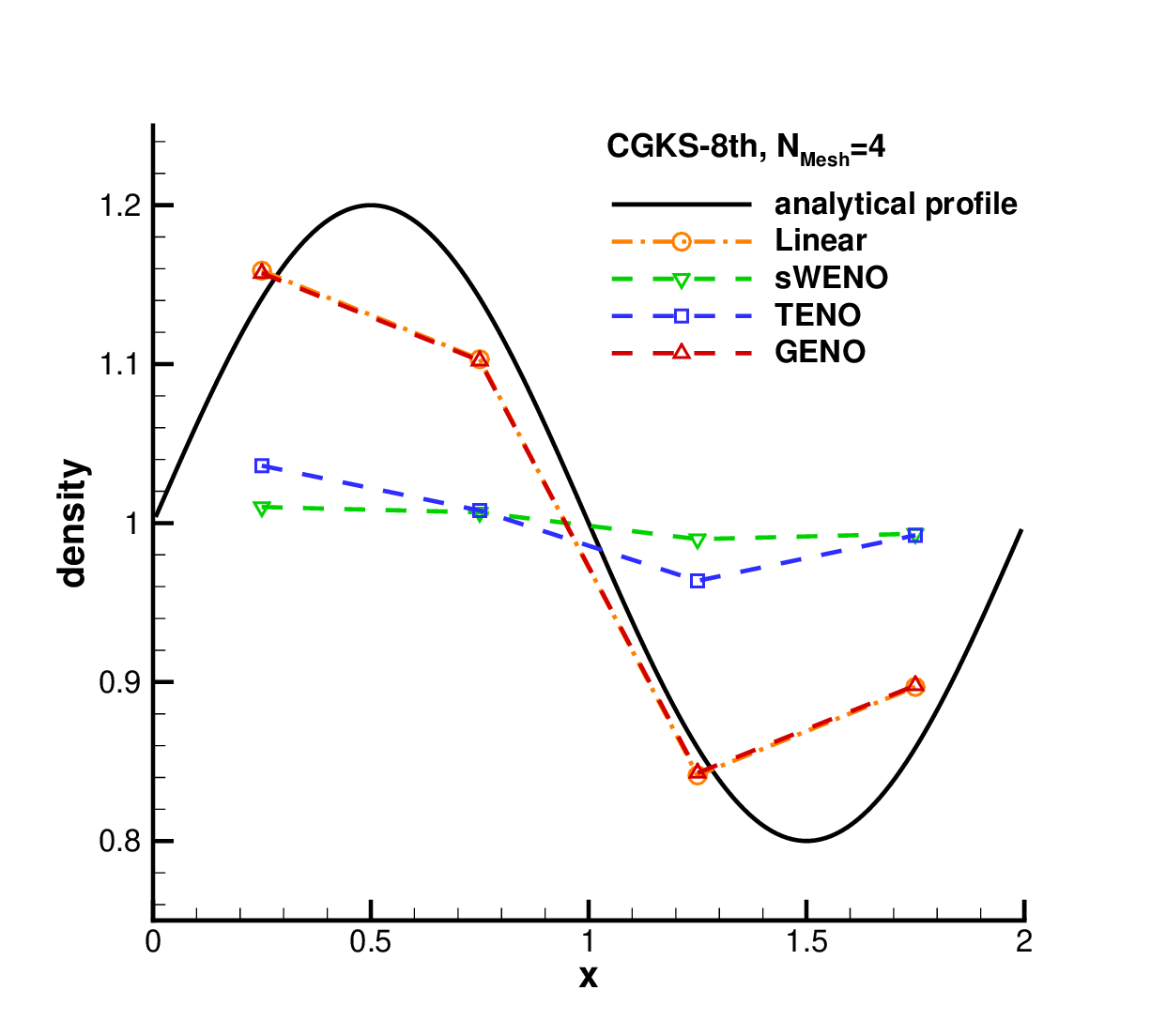}
\caption{\label{1d-accuracy-test-3} 1-D accuracy test case: results of the 8th-order compact GKS with various nonlinear reconstruction methods. Error versus mesh refinement at $t = 2$ (left), and density wave profiles on the coarsest mesh after long-time evolution at $t = 20$ (right). }
\end{figure}

\begin{figure}[!htb]
\centering
\includegraphics[width=0.495\textwidth]{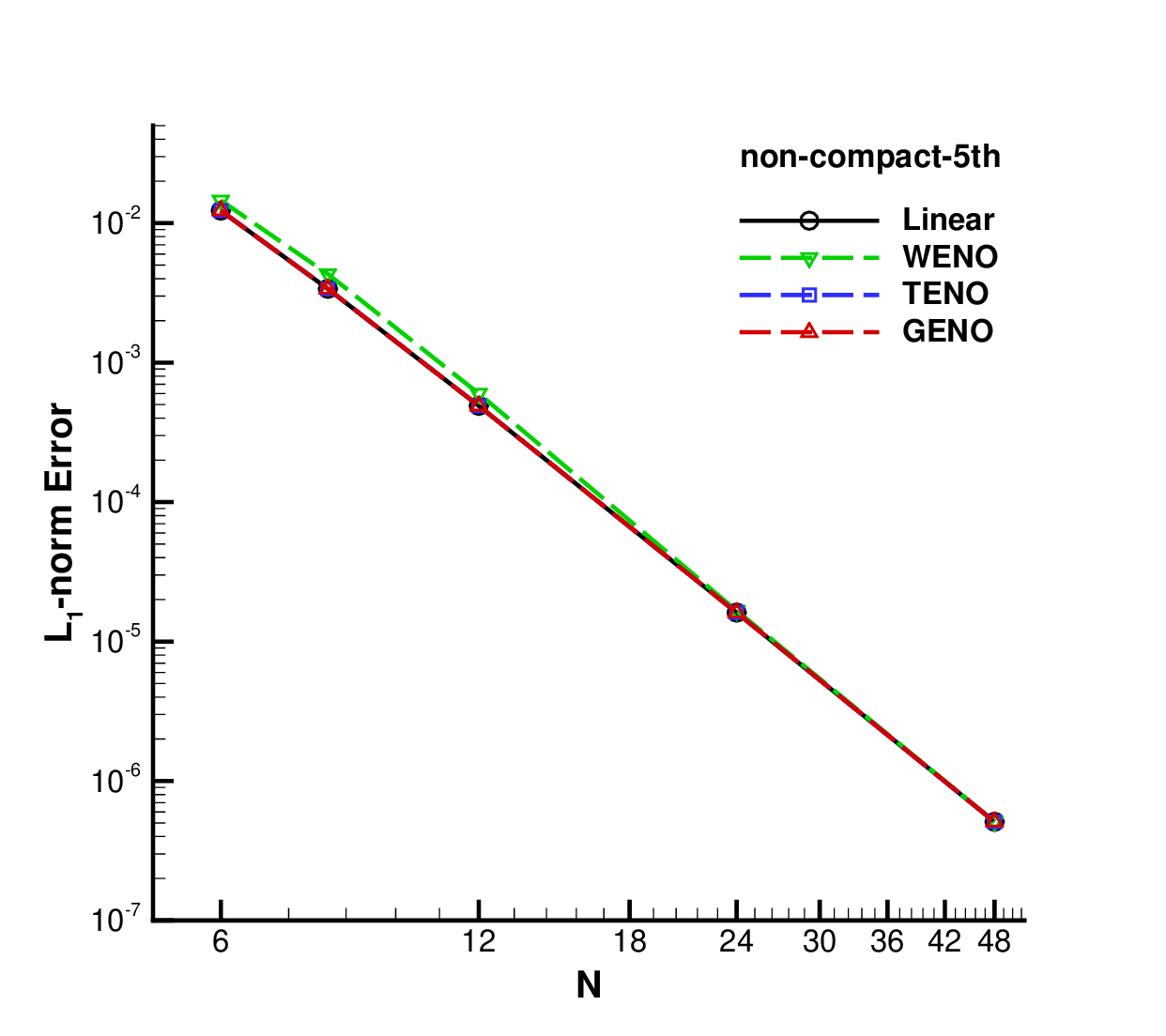}
\includegraphics[width=0.495\textwidth]{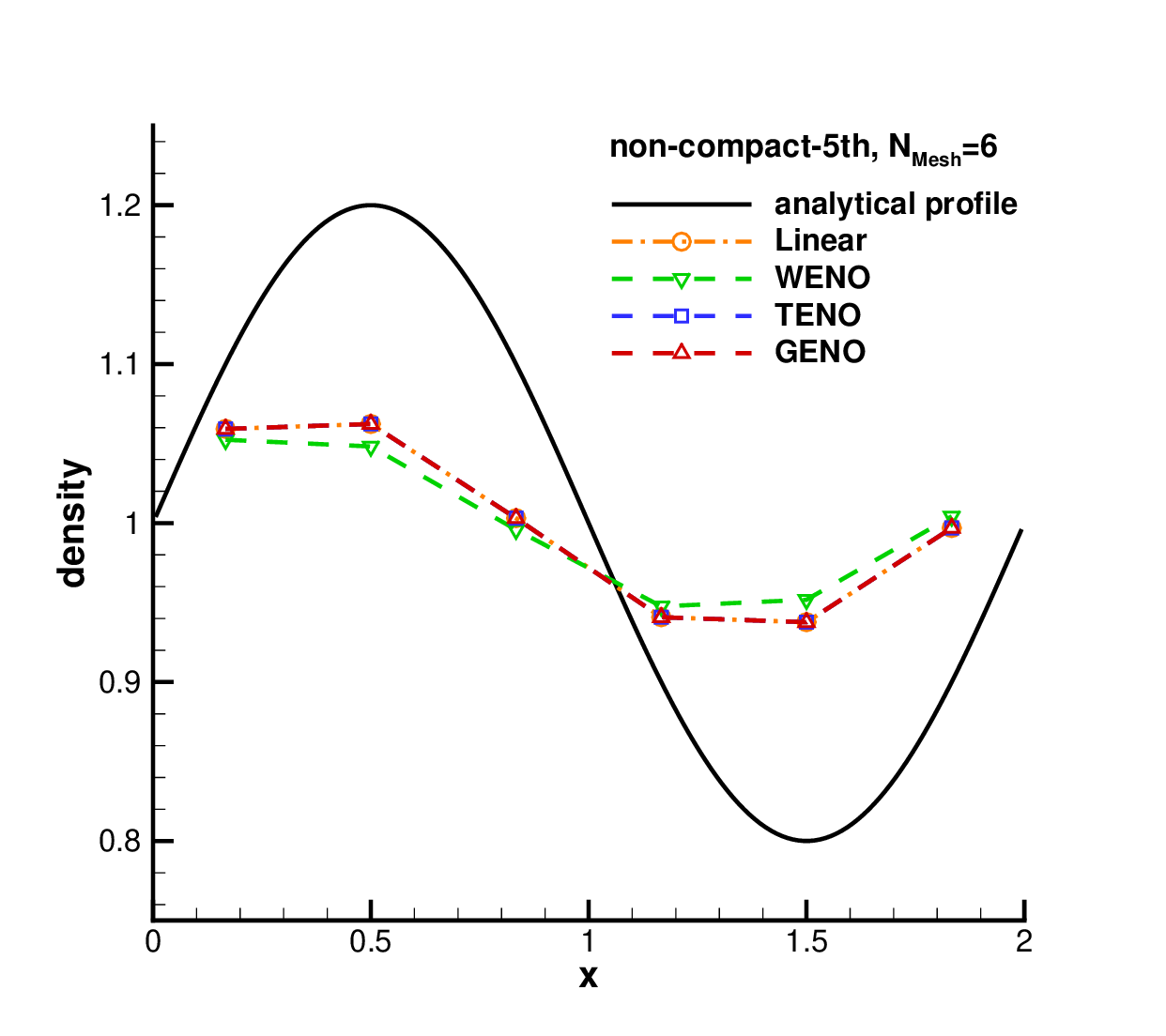}
\caption{\label{1d-accuracy-test-1} 1-D accuracy test case: results of the 5th-order non-compact scheme with various nonlinear reconstruction methods. Error versus mesh refinement at $t = 2$ (left), and density wave profiles on the coarsest mesh after long-time evolution at $t = 20$ (right). }
\end{figure}

\begin{figure}[!htb]
\centering
\includegraphics[width=0.495\textwidth]{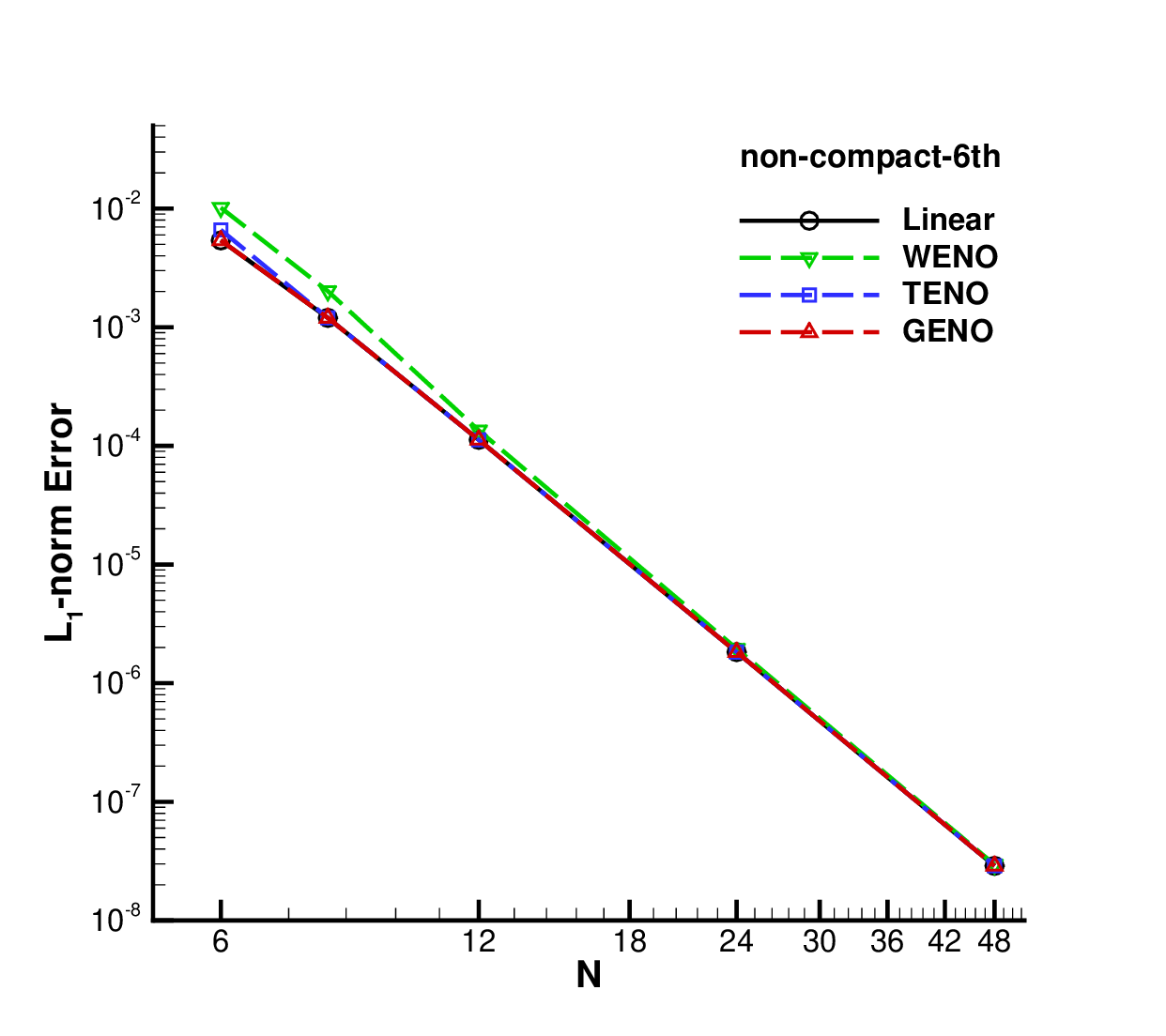}
\includegraphics[width=0.495\textwidth]{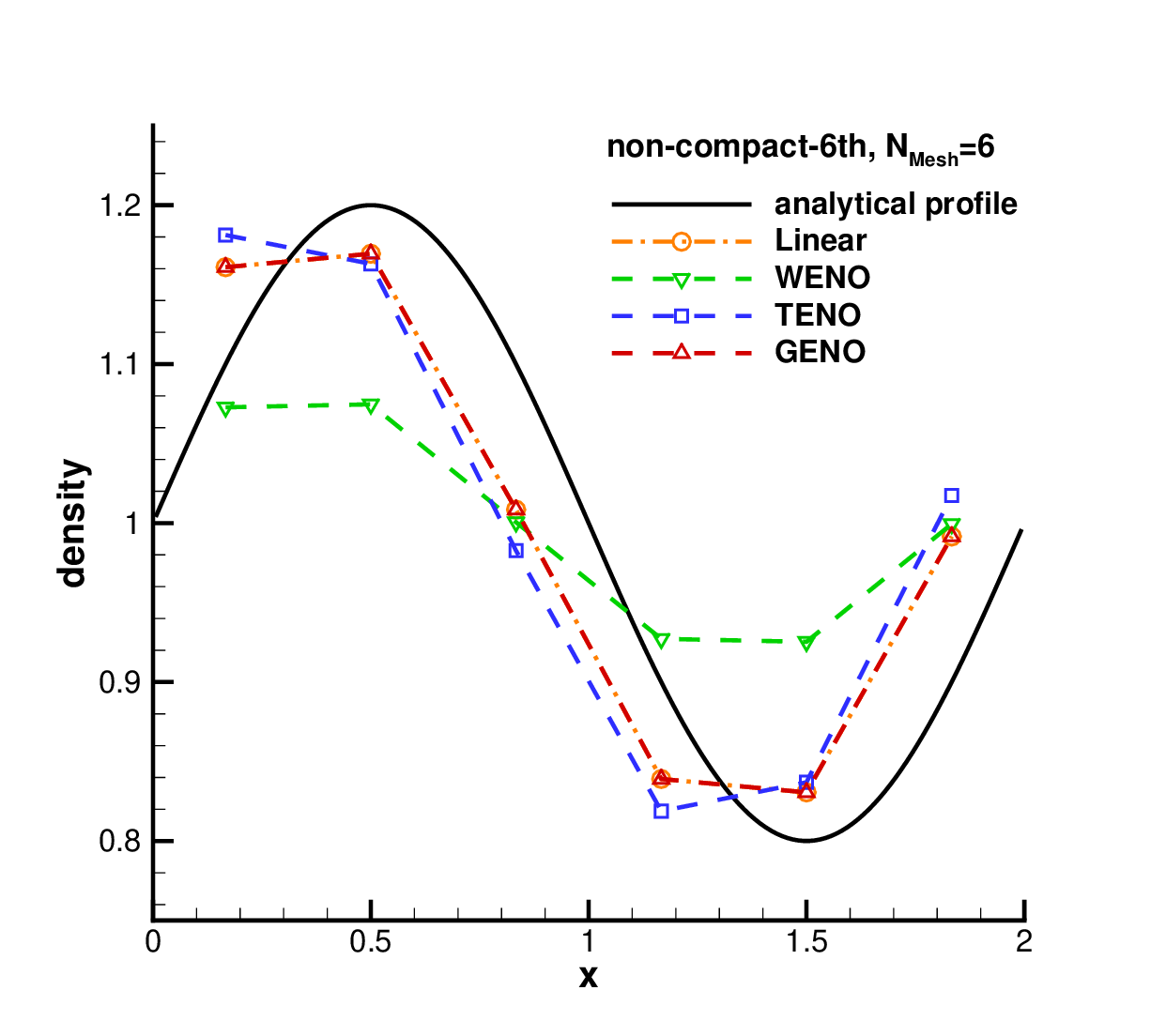}
\caption{\label{1d-accuracy-test-2} 1-D accuracy test case: results of the 6th-order non-compact scheme with various nonlinear reconstruction methods. Error versus mesh refinement at $t = 2$ (left), and density wave profiles on the coarsest mesh after long-time evolution at $t = 20$ (right). }
\end{figure}

Fig. \ref{1d-accuracy-test-3} to Fig. \ref{1d-accuracy-test-2} present results obtained by the 8th-order compact GKS and non-compact 5th- and 6th-order schemes using WENO (sWENO), TENO, and GENO reconstructions. The left of Fig. \ref{1d-accuracy-test-3} to Fig. \ref{1d-accuracy-test-2} depict error convergence with mesh refinement at $t=2$, where the coarsest meshes for the non-compact and compact schemes correspond to $6$ and $4$ cells, respectively. The right plots of Fig. \ref{1d-accuracy-test-3} to Fig. \ref{1d-accuracy-test-2} show results on the respective coarsest meshes after a long computation time ($t = 20$). When the mesh resolution is insufficient to resolve the flow, the accuracy of WENO (sWENO), TENO, and GENO reconstructions differs significantly, where the GENO yields the optimal results. Furthermore, the compact scheme achieves superior results with fewer mesh cells.

\subsection{1-D entropy wave advection}
A density wave convection problem with a complex initial density profile is tested.
The computational domain takes [-800, 1000] with periodic boundary conditions.
The initial condition is
\begin{equation*}
(\rho,U,p)=(1+0.5\mathrm{e}^{\mathrm{-ln2x^2/b^2}},1,1),
\end{equation*}
where the parameter $b$ takes the value $1.5$.
The computation was run up to $t=400$, corresponding to a long-time evolution of $400/3\approx 133.3$ entropy wave propagation periods (where one period is $2b$).
The computation uses $1800$ mesh cells.

\begin{figure}[!htb]
\centering
\includegraphics[width=0.495\textwidth]{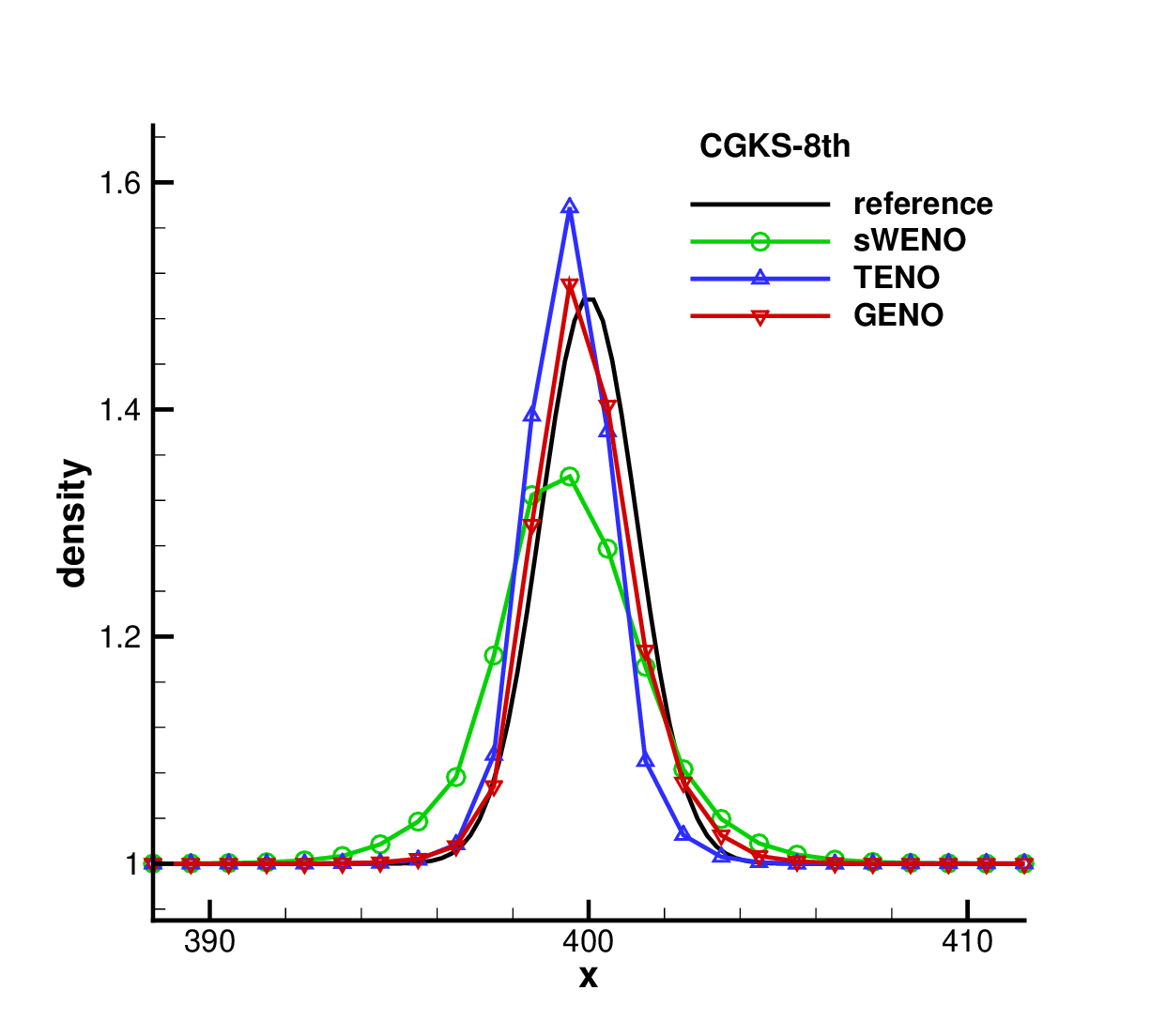}
\includegraphics[width=0.495\textwidth]{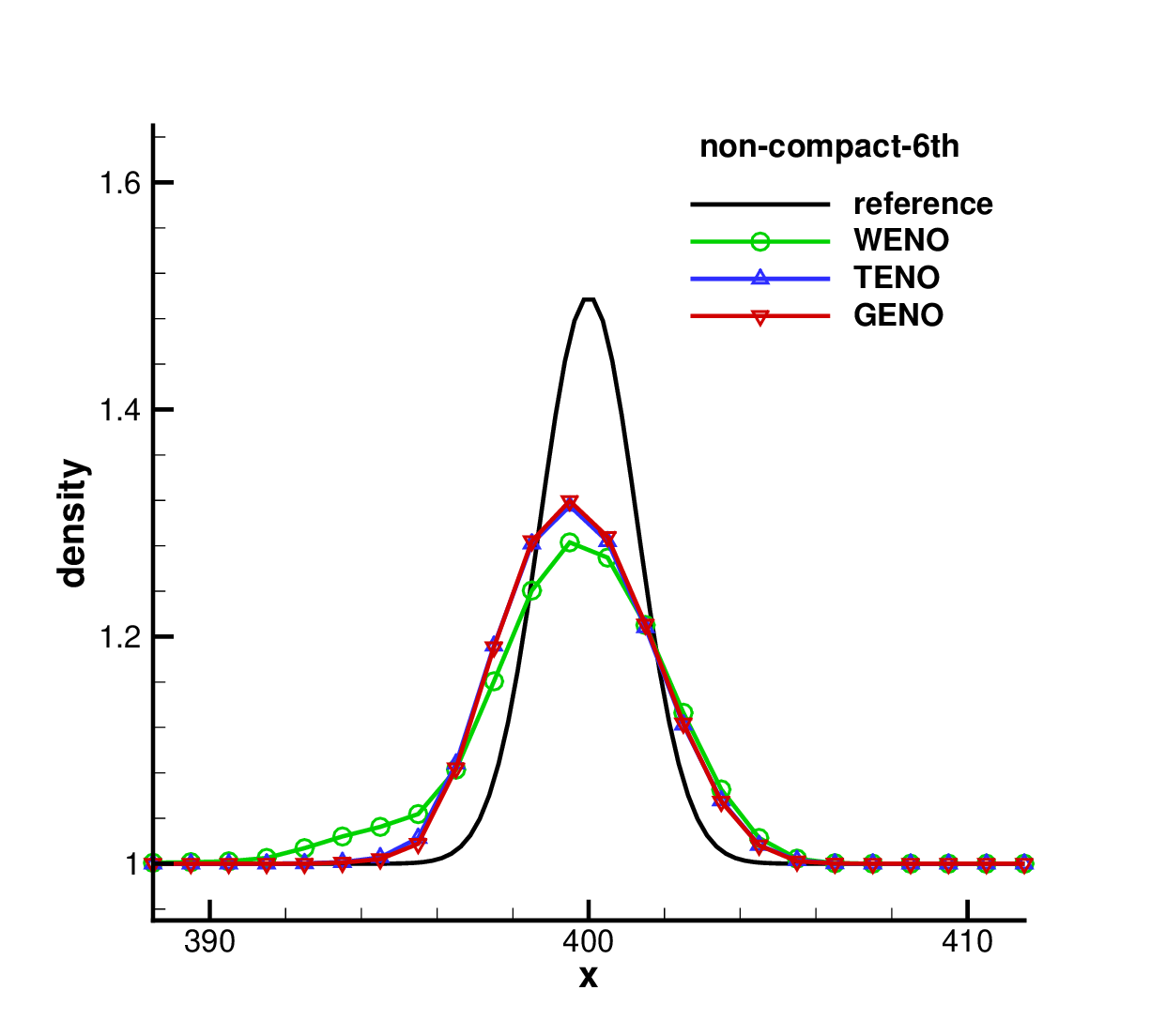}
\caption{\label{1D-entropy-wave-advection} Advection of 1-D entropy wave: the results are obtained at $t=400$ by the 8th-order compact GKS (left) and the 6th-order noncompact scheme (right) using WENO (sWENO), TENO and GENO nonlinear reconstructions after a long period of wave propagation. The mesh cell size is $1$.}
\end{figure}

In Fig. \ref{1D-entropy-wave-advection}, the left and right plots show results for the compact 8th-order GKS and non-compact 6th-order scheme, respectively.
GENO outperforms WENO (sWENO) and TENO in terms of accuracy.
Furthermore, the compact scheme's results are significantly better than those of the non-compact scheme.

\subsection{1-D shock-entropy wave interaction test case}
The Shu-Osher problem is a benchmark case for assessing the performance of high-order schemes in the coexistence of discontinuous and smooth waves.
The initial condition for the Shu-Osher problem is \cite{Case-Shu-Osher}:
\begin{equation*}
(\rho,U,p)=\left\{\begin{array}{ll}
(27/7, 4\sqrt{35}/9, 31/3),  \ \ \ \ &  x \leq 1,\\
(1 + 0.2\sin (5x), 0, 1),  &  1 <x \leq 10.
\end{array} \right.
\end{equation*}
The computational domain takes $[0, 10]$.
Free boundary conditions are imposed at both ends. The simulation is carried out up to $t = 1.8$.

\begin{figure}[!htb]
\centering
\includegraphics[width=0.495\textwidth]{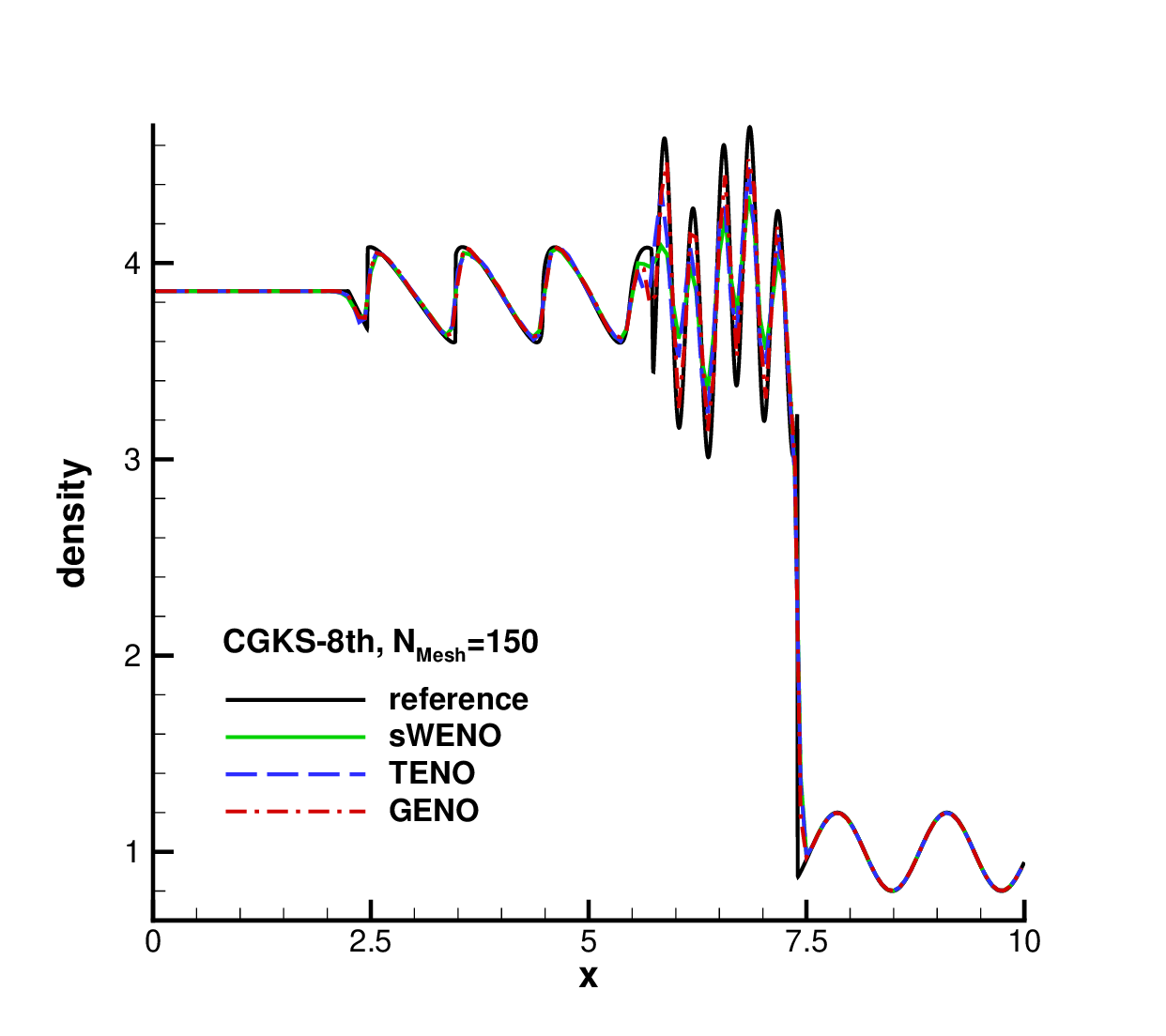}
\includegraphics[width=0.495\textwidth]{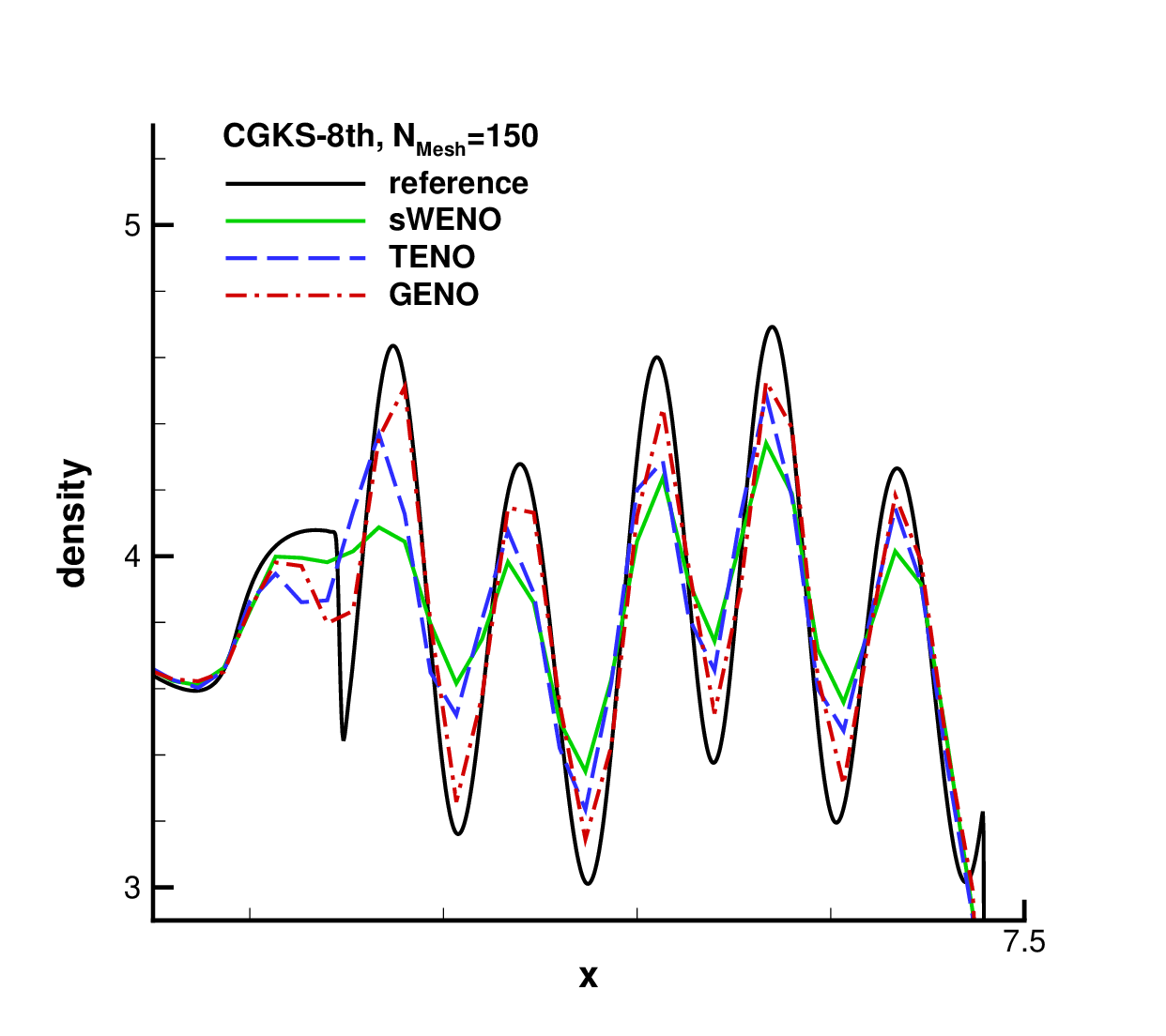}
\caption{\label{1d-shocktube-shuosher-3} Shu-Osher problem: density and the local enlargement with $N=150$ cells by the 8th-order compact GKS with different nonlinear reconstructions. }
\end{figure}

\begin{figure}[!htb]
\centering
\includegraphics[width=0.495\textwidth]{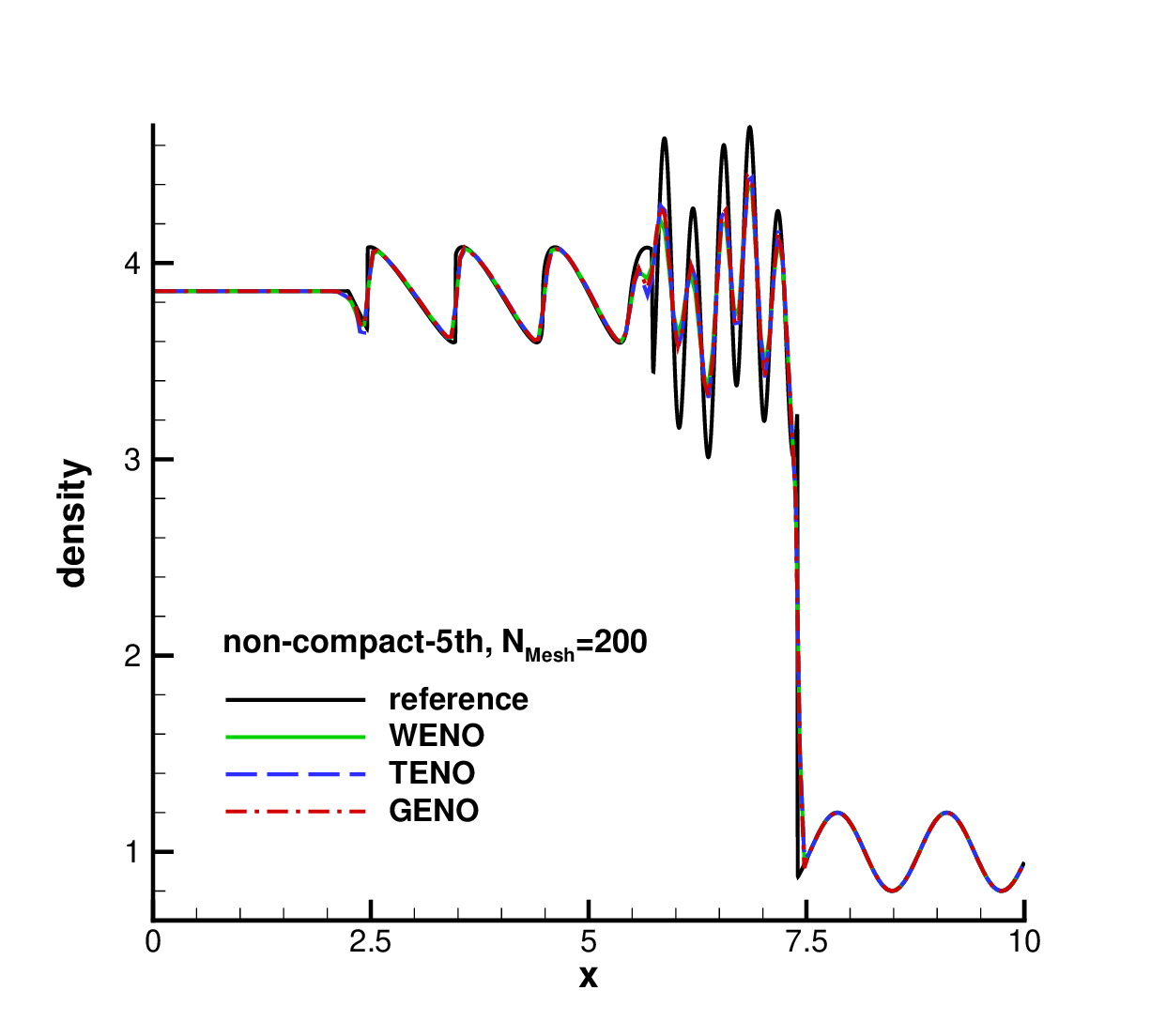}
\includegraphics[width=0.495\textwidth]{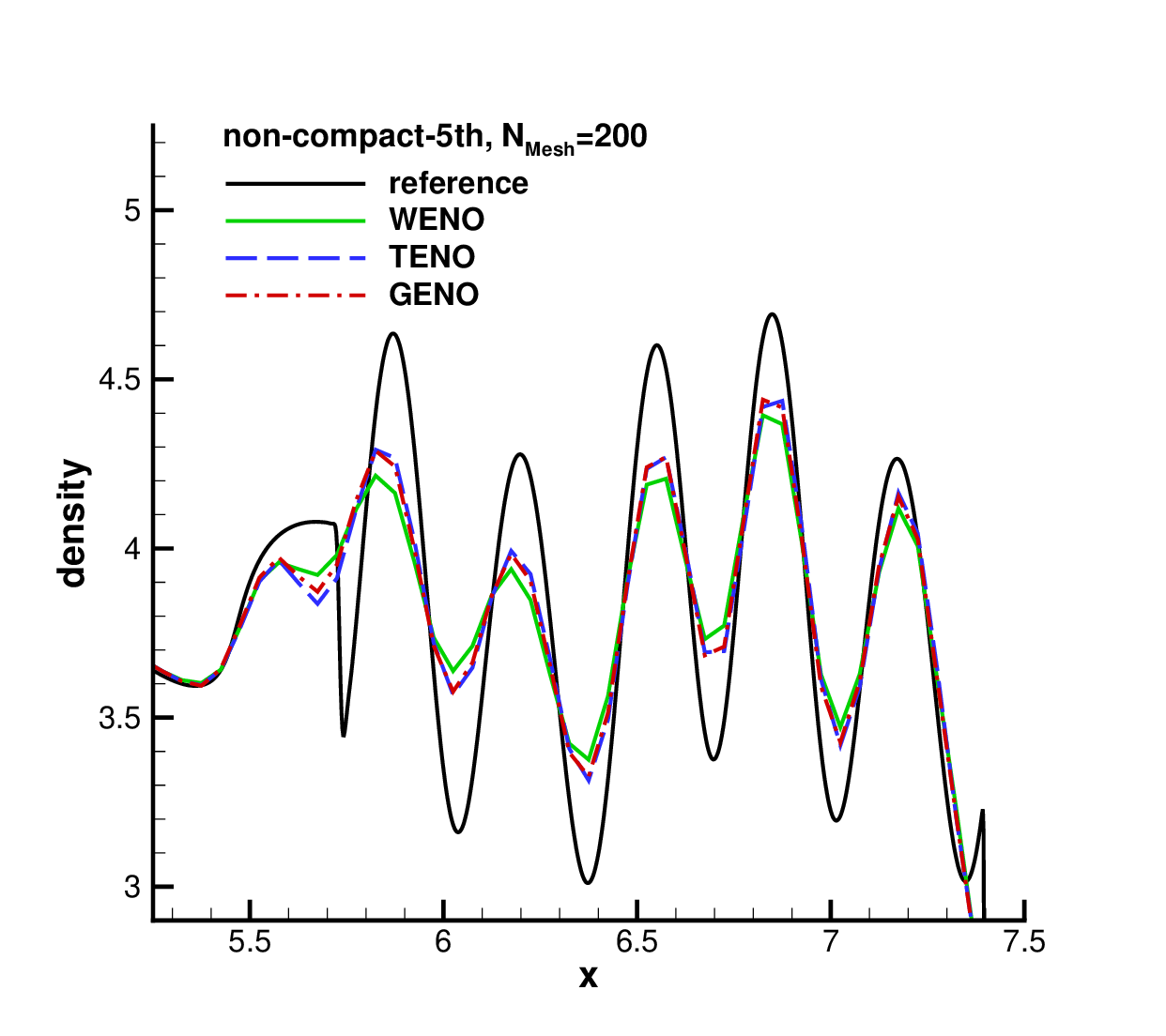}
\caption{\label{1d-shocktube-shuosher-1} Shu-Osher problem: density and the local enlargement with $N=200$ cells by the 5th-order non-compact scheme with different nonlinear reconstructions. }
\end{figure}

\begin{figure}[!htb]
\centering
\includegraphics[width=0.495\textwidth]{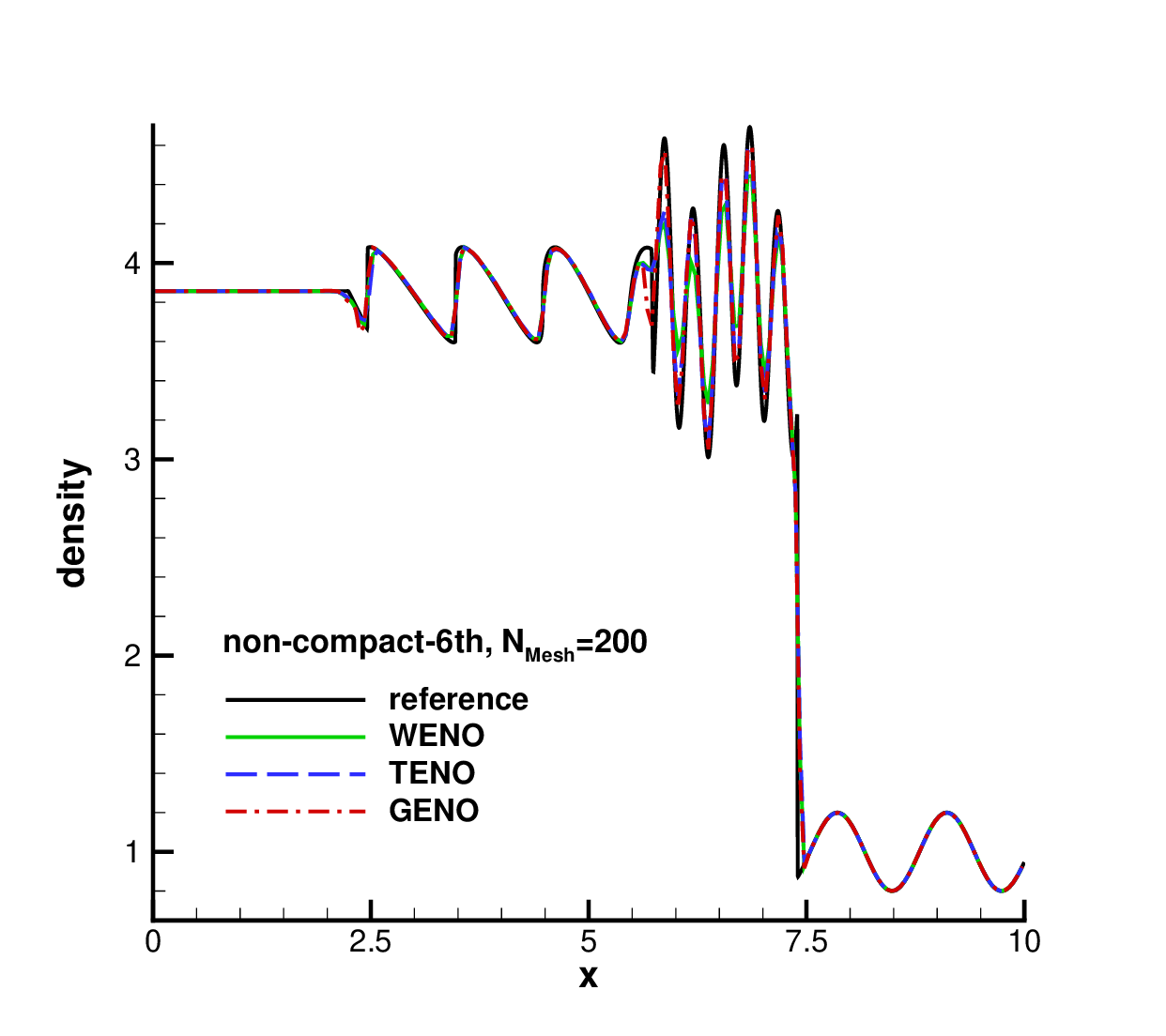}
\includegraphics[width=0.495\textwidth]{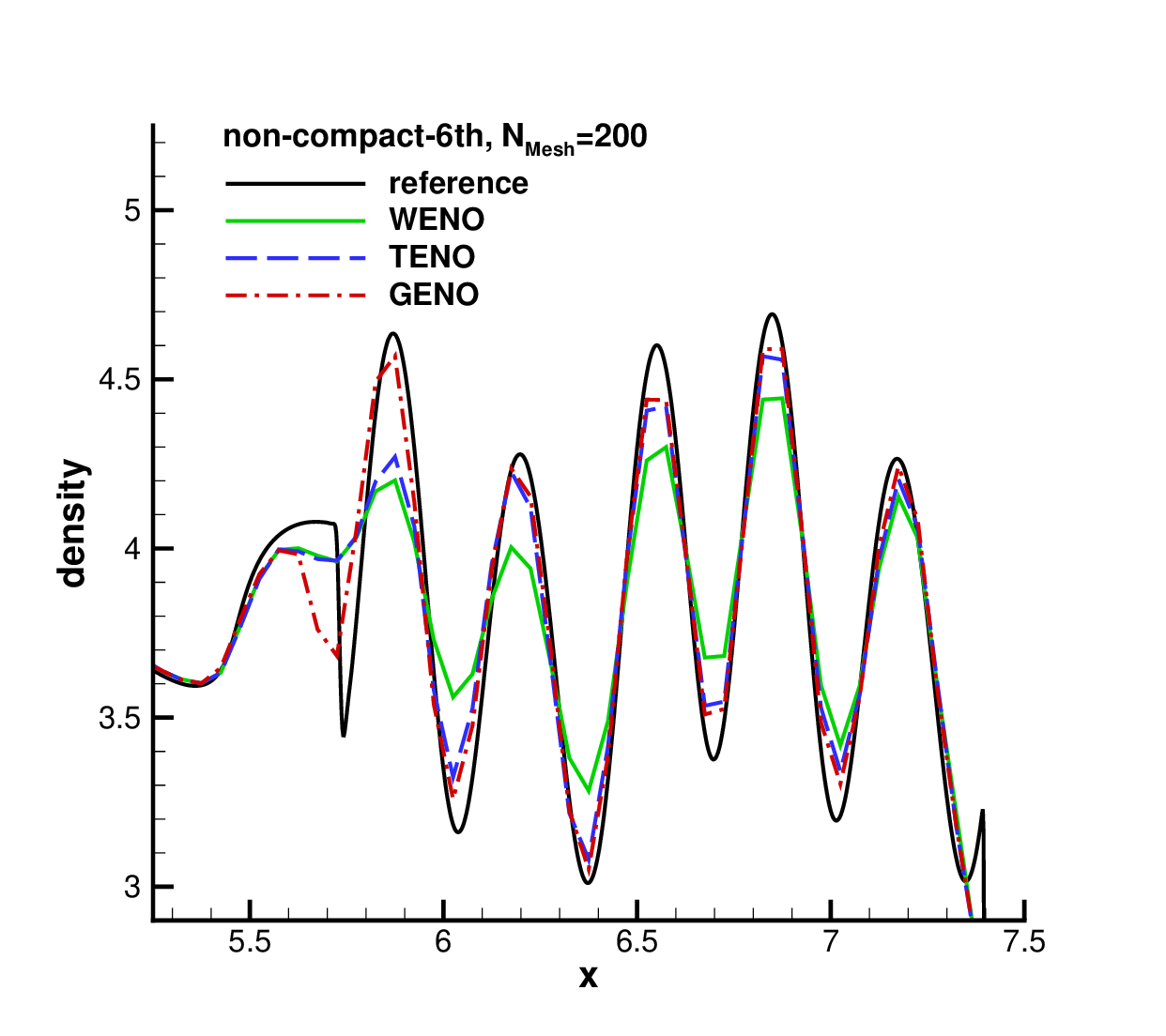}
\caption{\label{1d-shocktube-shuosher-2} Shu-Osher problem: density and the local enlargement with $N=200$ cells by the 6th-order non-compact scheme with different nonlinear reconstructions. }
\end{figure}

Fig. \ref{1d-shocktube-shuosher-3} to Fig. \ref{1d-shocktube-shuosher-2} present results obtained using a compact 8th-order GKS and non-compact 5th- and 6th-order schemes, each employing three nonlinear reconstructions.
The compact and noncompact schemes use $150$ and $200$ cells, respectively.
The left and right plots of Fig. \ref{1d-shocktube-shuosher-3} to Fig. \ref{1d-shocktube-shuosher-2} depict the density and a local enlargement at the output time.
When the mesh resolution is insufficient to resolve the flow,  as the case with the coarse mesh here, the accuracy of WENO (sWENO), TENO, and GENO reconstructions differs significantly, where GENO yields the optimal results.
Furthermore, the compact scheme achieves superior results with fewer mesh cells.

As an extension of the Shu-Osher problem, the Titarev-Toro problem \cite{titarev2004finite} is tested as well, and the initial condition in this case is the following
\begin{equation*}
(\rho,U,p)=\left\{\begin{array}{ll}
(1.515695, 0.523346, 1.805),  \ \ \ \ &  x \leq -4,\\
(1 + 0.1\sin (20\pi x), 0, 1),  &  -4 <x.
\end{array} \right.
\end{equation*}
The computational domain is $[-5,5]$.
The inflow boundary condition is imposed on left end, and the fixed wave profile is given on the right end. The computational time is $t=5.0$.

\begin{figure}[!htb]
\centering
\includegraphics[width=0.495\textwidth]{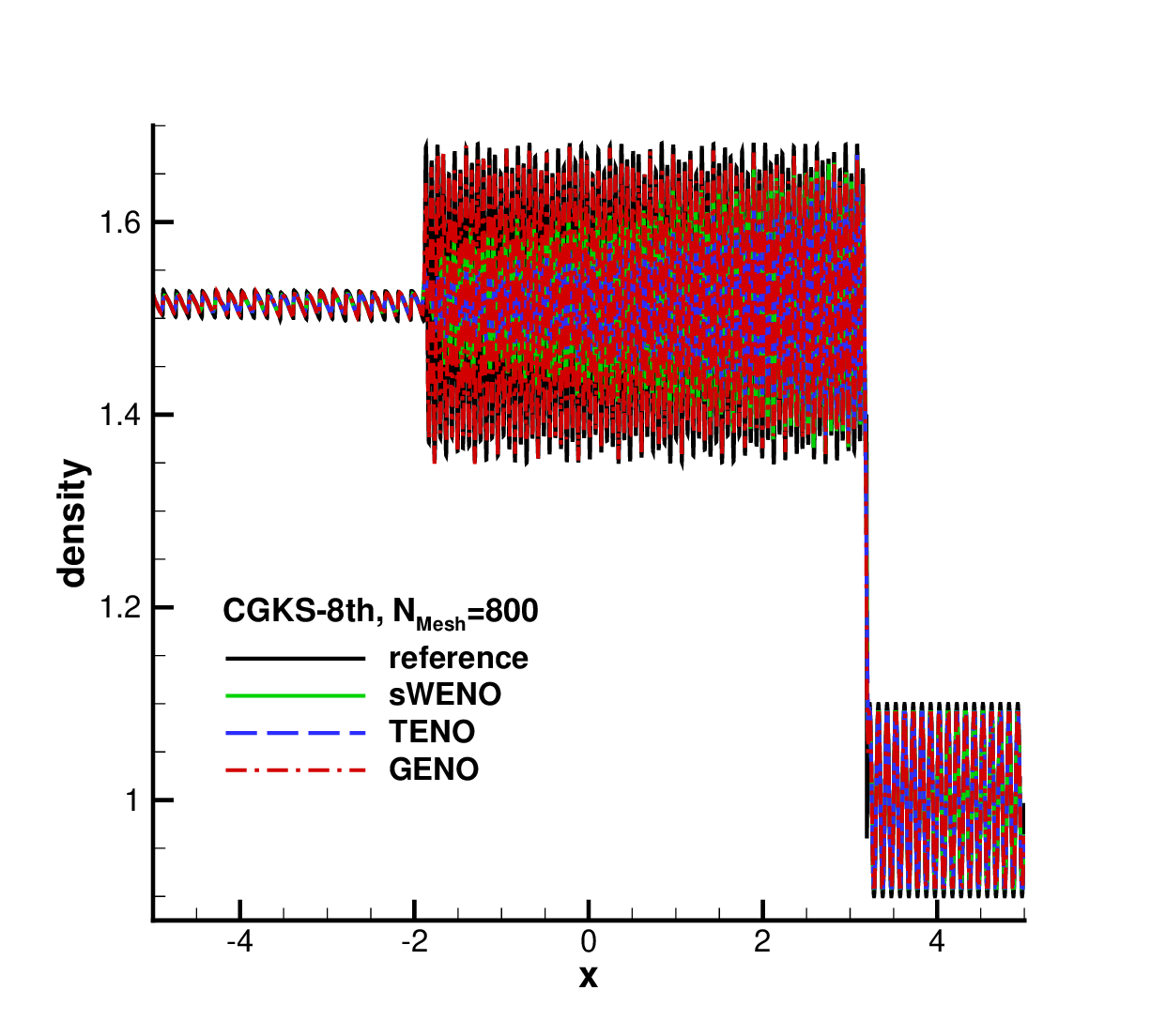}
\includegraphics[width=0.495\textwidth]{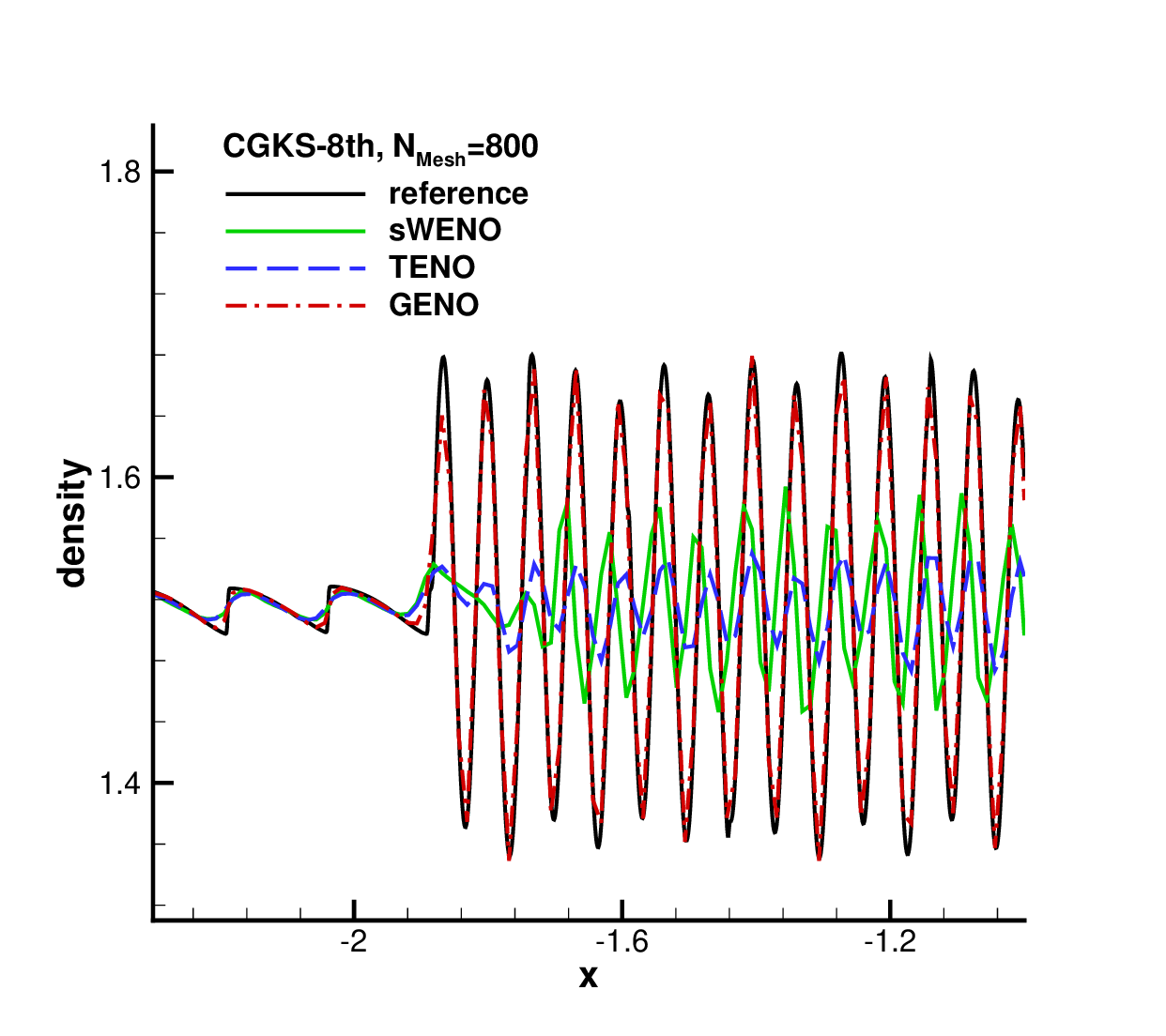}
\caption{\label{1d-shocktube-toro-2} Titarev-Toro problem: density and the local enlargement with $N=800$ cells by the 8th-order compact GKS with different nonlinear reconstructions. }
\end{figure}

\begin{figure}[!htb]
\centering
\includegraphics[width=0.495\textwidth]{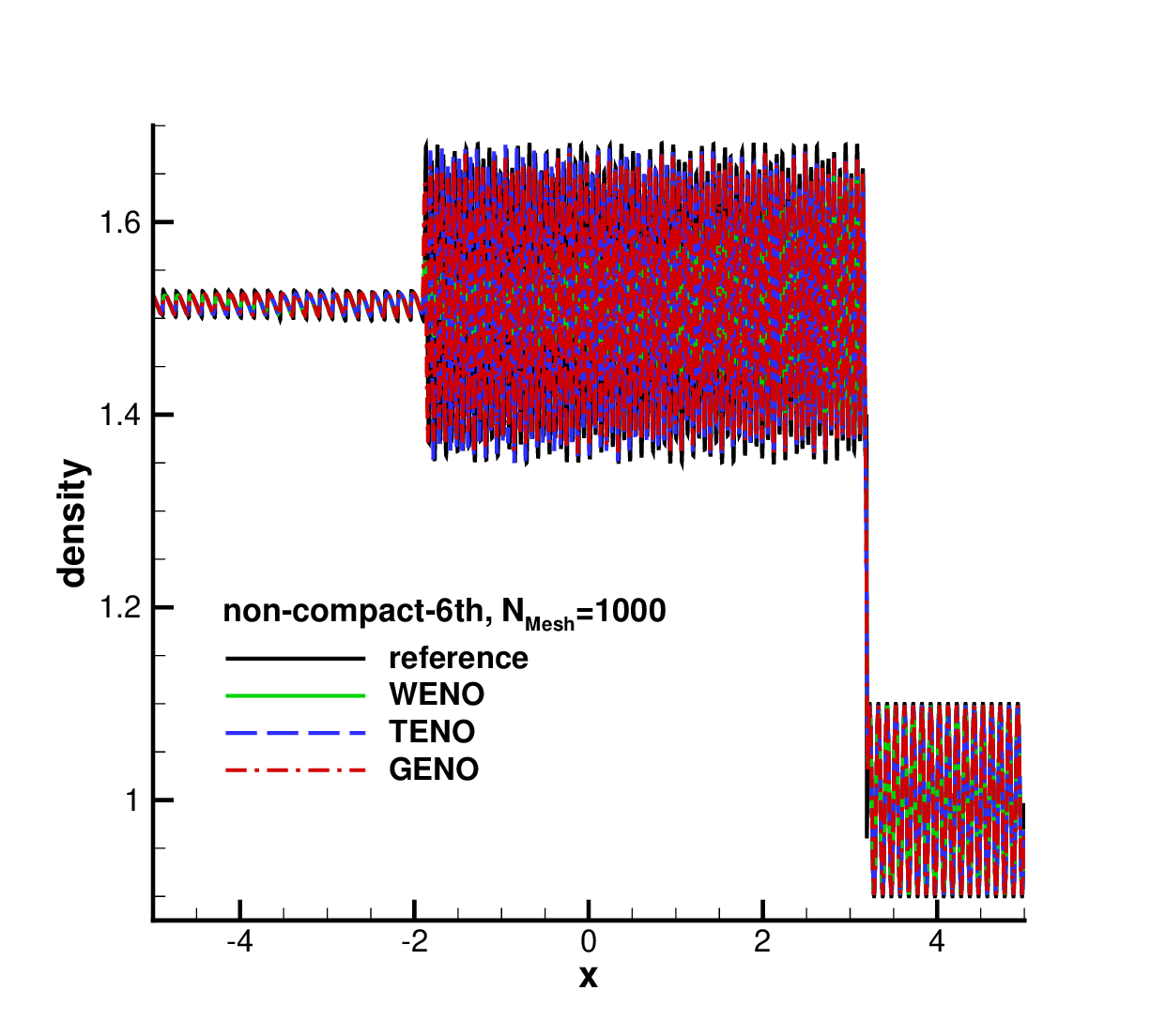}
\includegraphics[width=0.495\textwidth]{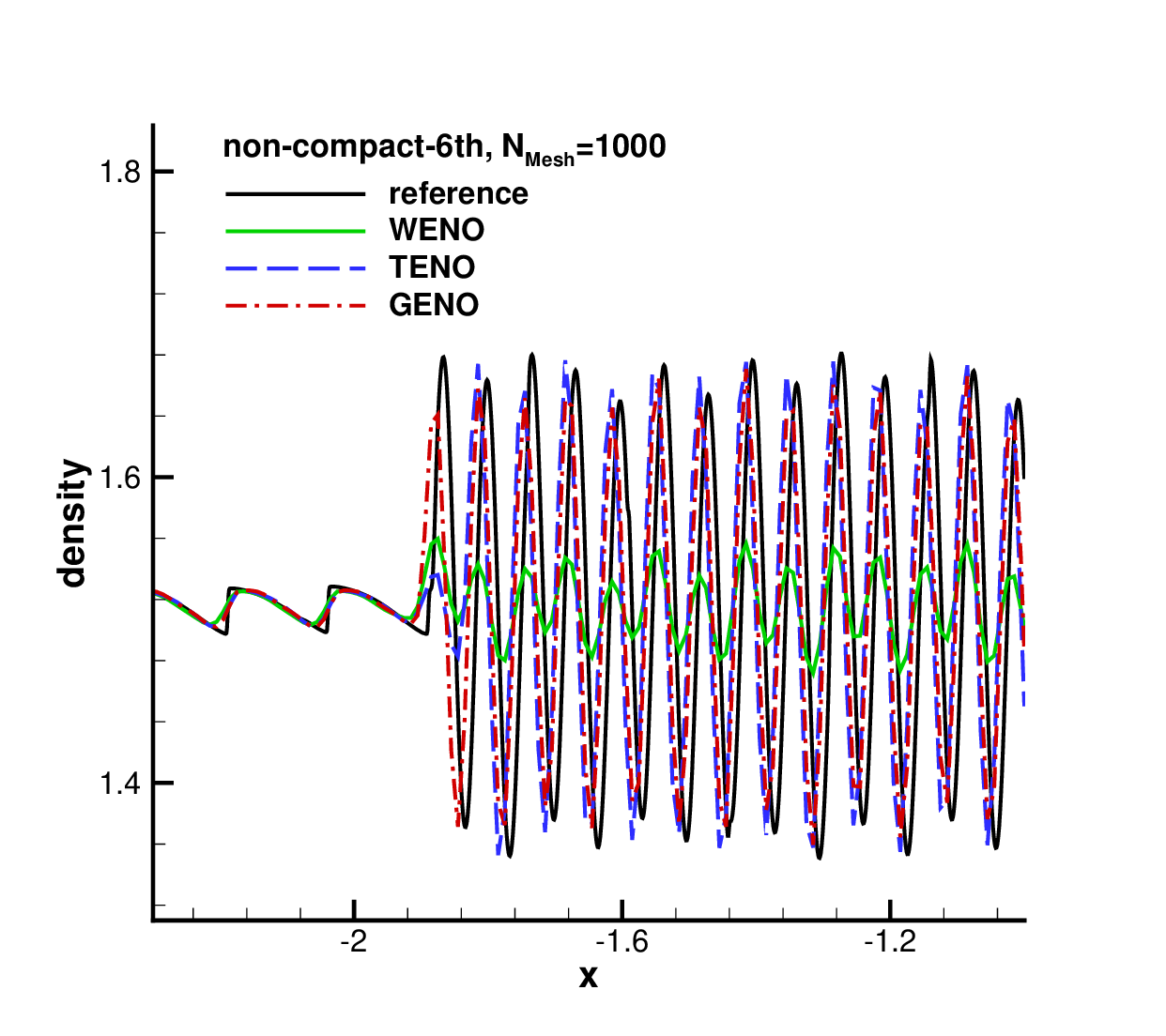}
\caption{\label{1d-shocktube-toro-1} Titarev-Toro problem: density and the local enlargement with $N=1000$ cells by the non-compact 6th-order scheme with different nonlinear reconstructions. }
\end{figure}

Fig. \ref{1d-shocktube-toro-2} and Fig. \ref{1d-shocktube-toro-1} present results obtained using the 8th-order compact GKS and the non-compact 6th-order scheme, each employing three nonlinear reconstructions.
The compact and non-compact schemes use $800$ and $1000$ cells, respectively.
The left and right plots depict the density and a local enlargement.
The accuracy of WENO (sWENO), TENO, and GENO reconstructions differs significantly.
It should be noted that in this specific case, the 8th-order compact GKS with TENO reconstruction requires the setting of the TENO weight parameter $C_T$ to $10^{-3}$
in order to achieve numerical stability. In all other cases, $C_T$ was consistently set to $10^{-6}$. Furthermore, the right plot of Fig. \ref{1d-shocktube-toro-1} shows that TENO exhibits overshoots and undershoots at the peaks and troughs near $x= -1.6$.

\subsection{1-D blast wave test case}
The Woodward-Colella blast wave problem \cite{Case-Woodward} is computed. The test case verifies the robustness of high-order schemes for capturing strong shock wave.
The initial condition is given as follows
\begin{align*}
(\rho,U,p) =\begin{cases}
(1, 0, 1000), & 0\leq x<0.1,\\
(1, 0, 0.01), & 0.1\leq x<0.9,\\
(1, 0, 100),  & 0.9\leq x\leq 1.
\end{cases}
\end{align*}
The computational domain is $[0,1]$, and the reflecting boundary conditions are imposed on both directions. A mesh with $N=400$ cells is used. The computation time for this case is $t=0.038$.

\begin{figure}[!htb]
\centering
\includegraphics[width=0.495\textwidth]{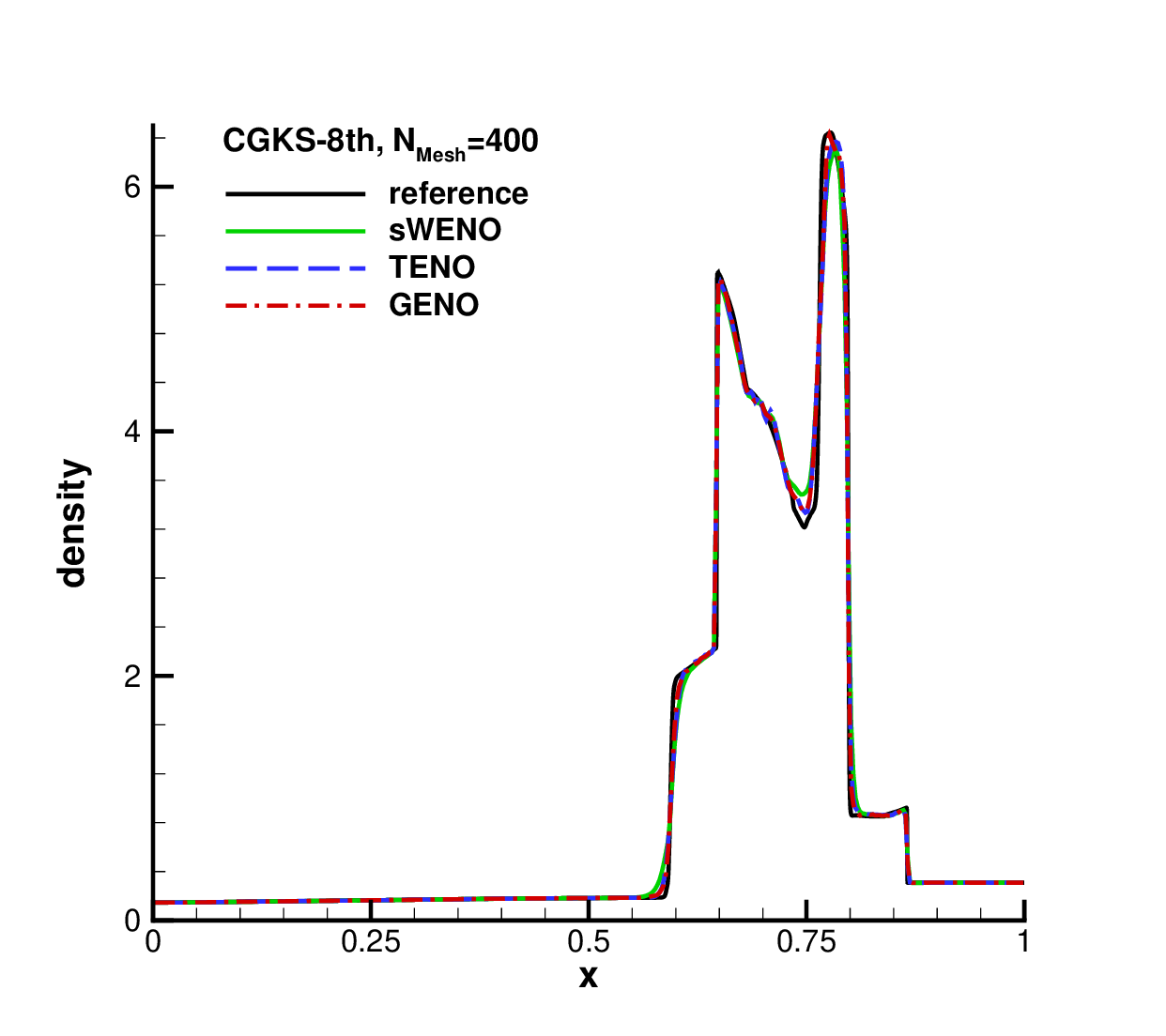}
\includegraphics[width=0.495\textwidth]{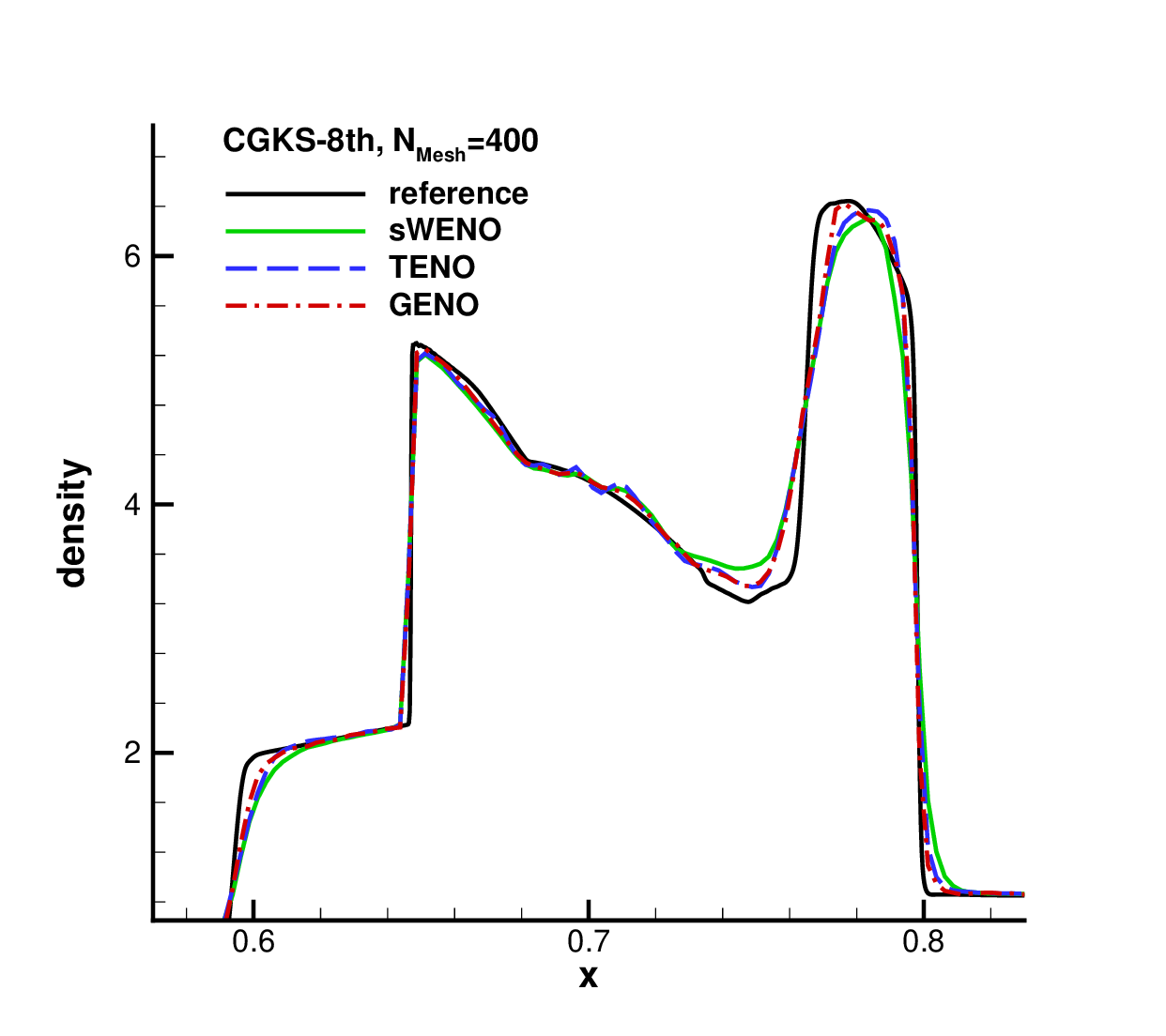}
\caption{\label{1d-shocktube-blast-1} Blast wave problem: density and the local enlargement with $N=400$ cells by the compact GKS with different nonlinear reconstructions. }
\end{figure}

\begin{figure}[!htb]
\centering
\includegraphics[width=0.495\textwidth]{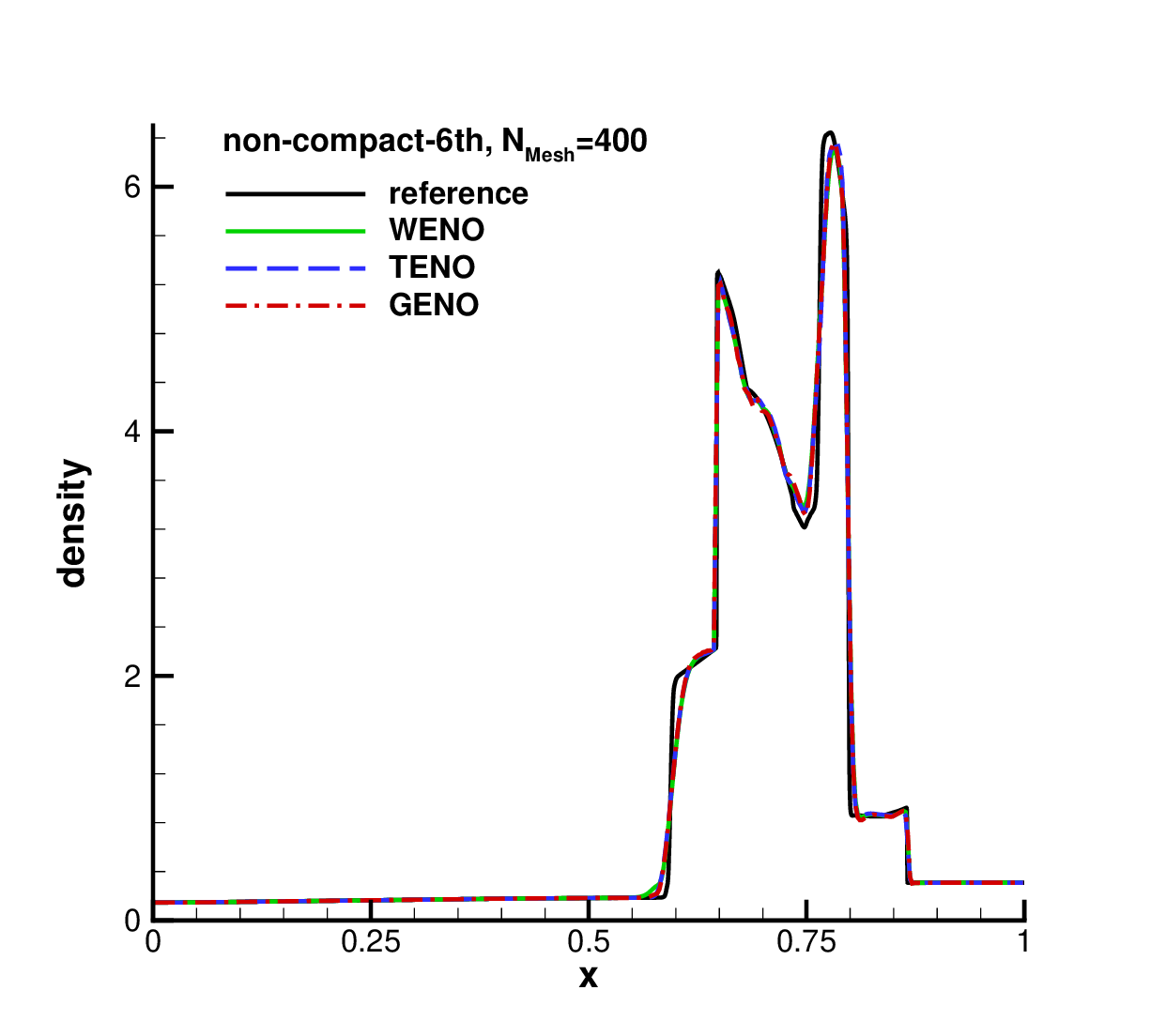}
\includegraphics[width=0.495\textwidth]{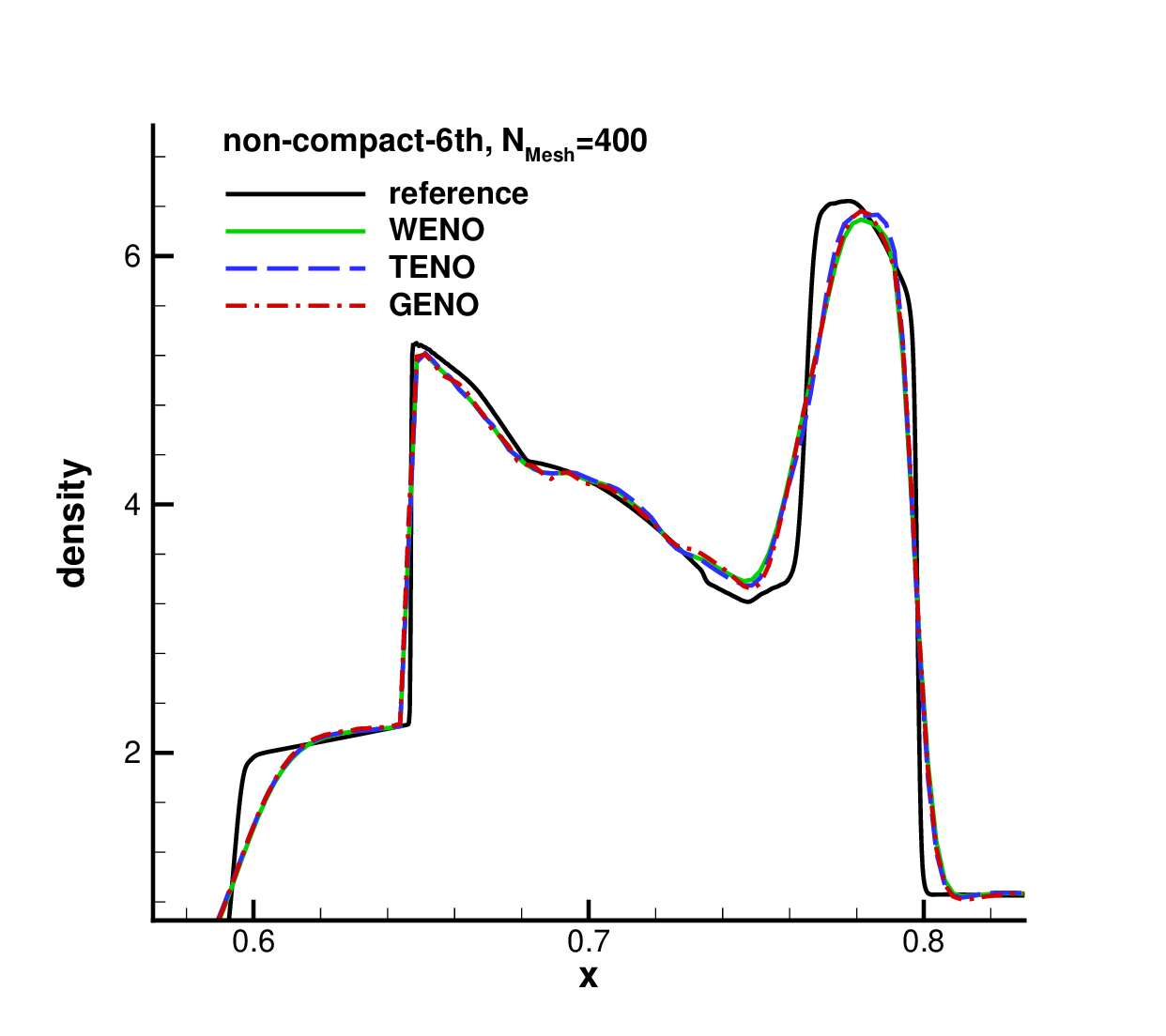}
\caption{\label{1d-shocktube-blast-2} Blast wave problem: density and the local enlargement with $N=400$ cells by the non-compact 6th-order scheme with different nonlinear reconstructions. }
\end{figure}

Fig. \ref{1d-shocktube-blast-1} and Fig. \ref{1d-shocktube-blast-2} present the density and a local enlargement, as obtained by the compact and non-compact high-order schemes, respectively. All three nonlinear reconstructions yield similar results, demonstrating the robustness of GENO in handling strong shock waves.

\subsection{2-D hypersonic flow around a cylinder}
Furthermore, to verify the applicability of GENO method on unstructured mesh, a hypersonic inviscid flow around a cylinder is first considered.
The incoming flow has a Mach number $Ma=8$. The adiabatic reflective boundary condition is imposed on the wall of the cylinder.

\begin{figure}[!htb]
	\centering
	\includegraphics[width=0.32\textwidth]{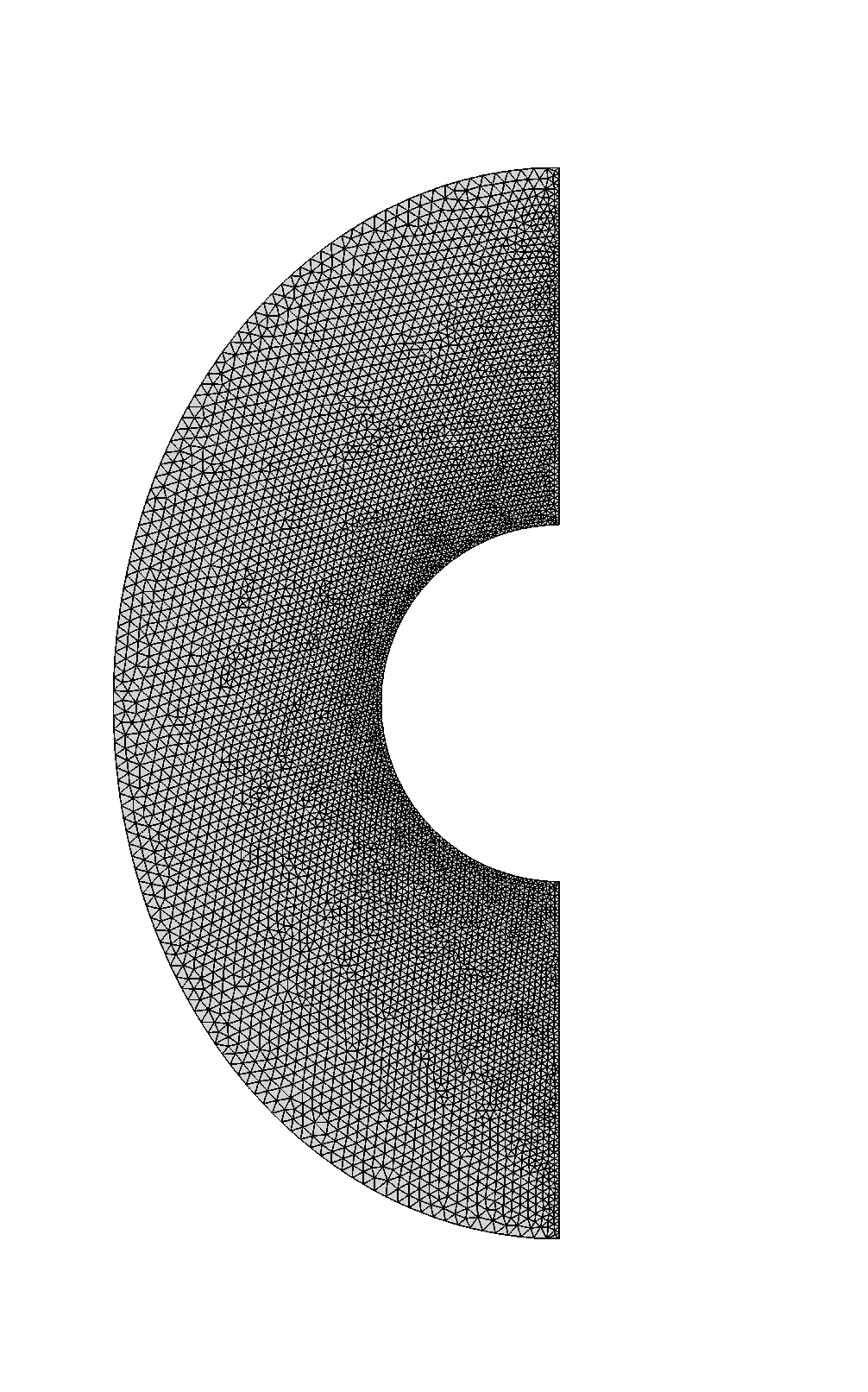}
	\includegraphics[width=0.32\textwidth]{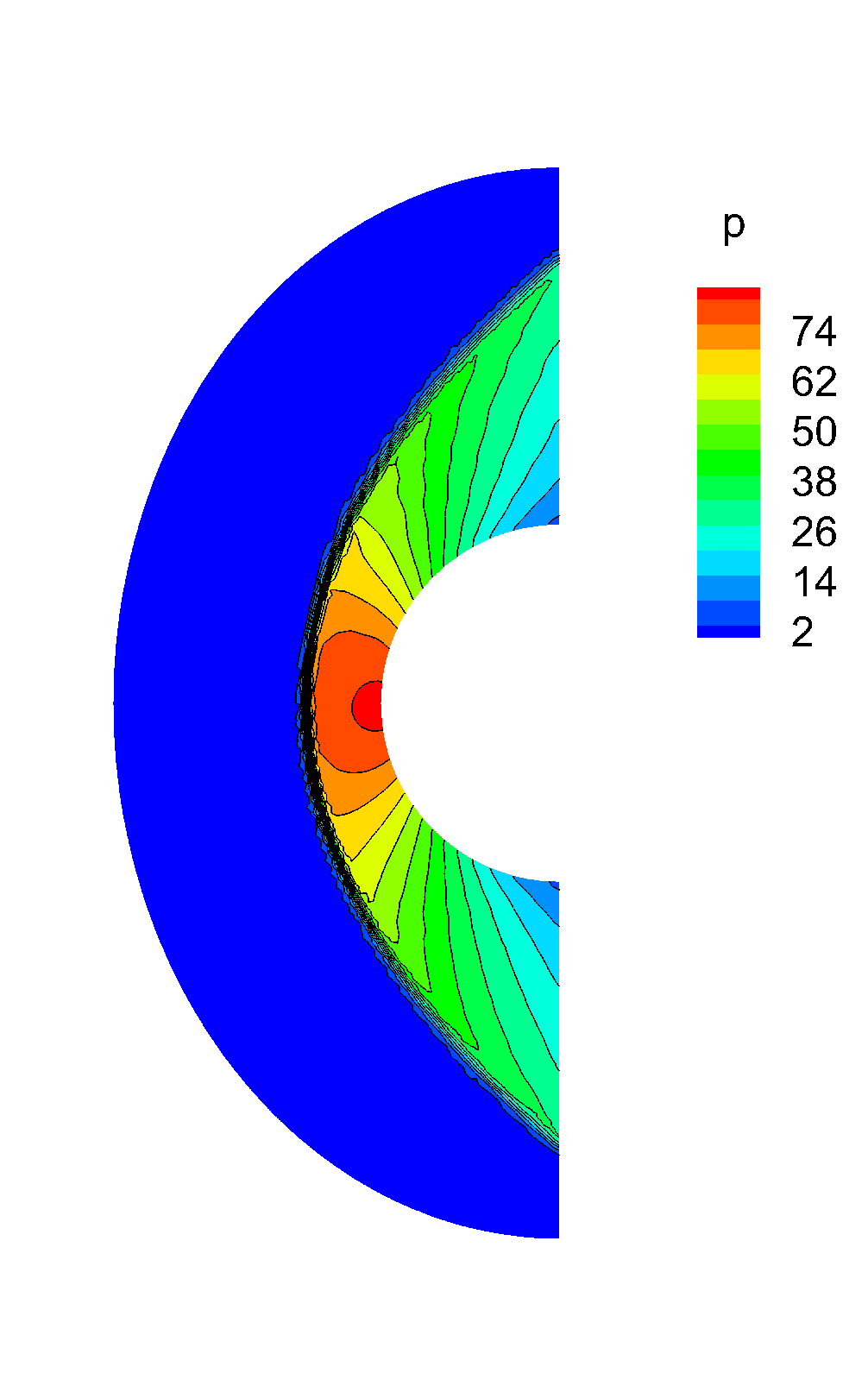}
	\includegraphics[width=0.32\textwidth]{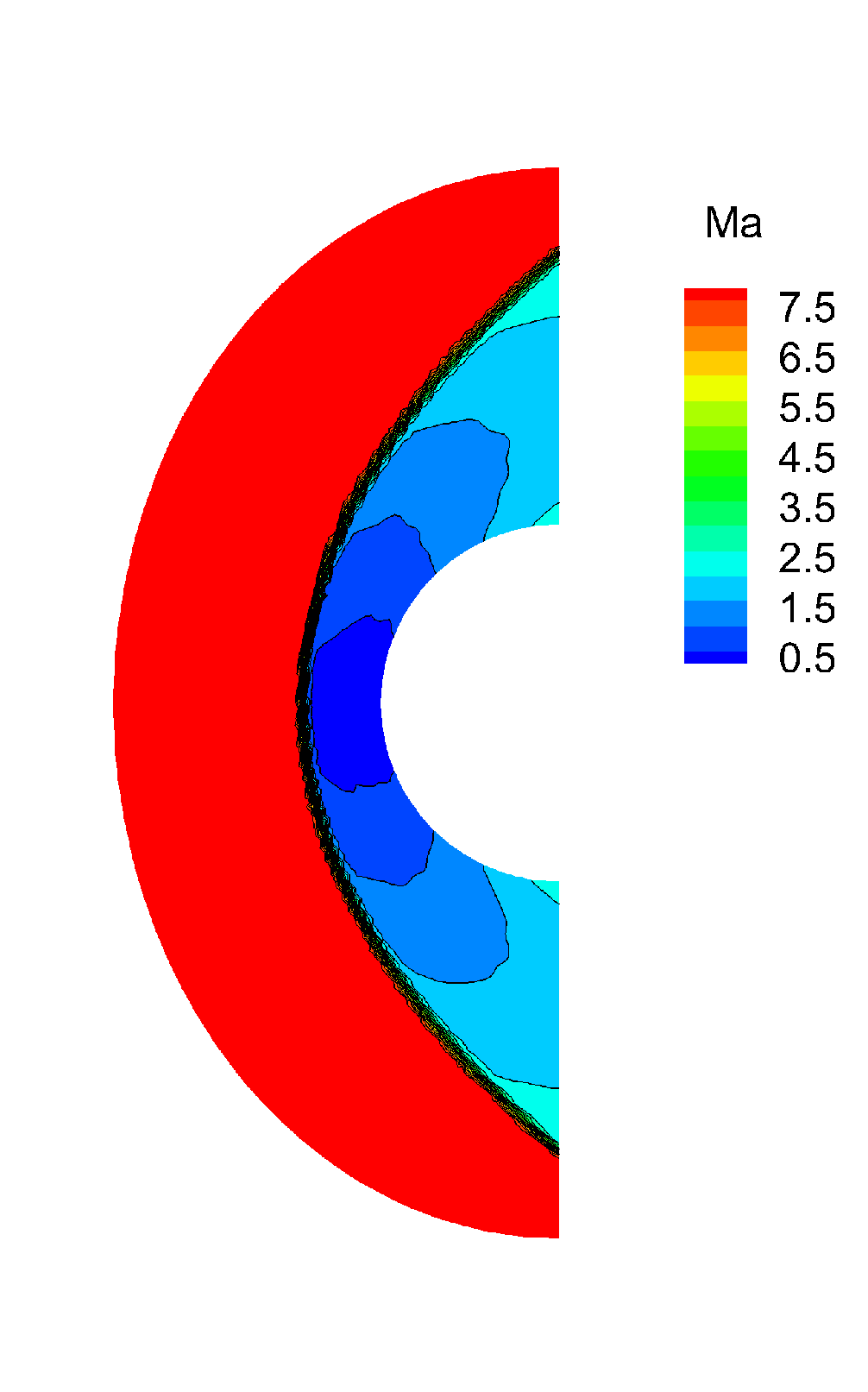}
	\caption{\label{2d-inv-cylinder} Hypersonic inviscid flow passing through a circular cylinder: pressure distributions with Mach number $Ma =8$ on triangular mesh.}
\end{figure}

The mesh is presented in Fig. \ref{2d-inv-cylinder}. The regular triangular mesh is used, and the mesh is refined near the cylinder.
Fig. \ref{2d-inv-cylinder} also presents contours of pressure and Mach number after the computation has reached a steady state, exhibiting good symmetry and no numerical oscillations.

\subsection{2-D Mach 3 step problem}

The step problem, a 2-D benchmark for inviscid flow simulations involving strong shock waves and their interactions \cite{Case-Woodward}. The computational domain is defined as $[0,3]\times[0,1] \setminus [0.6,3]\times[0,0.2]$, representing a wind tunnel of height $1$ and length $3$ with a step located at $x=0.6$, having a height of $0.2$. The initial state within the tunnel is characterized by $\rho=1, U=3, V=0, p=1/1.4$. This same state is enforced as the inflow boundary condition on the left. Slip boundary conditions are applied to the upper and lower walls.
The corner of the step initiates a rarefaction fan. In contrast to the approach taken in \cite{Case-Woodward}, no modifications to the density or velocity magnitude are applied to cells near the corner in the present computation.

\begin{figure}[!htb]
	\centering
    \includegraphics[width=0.80\textwidth]{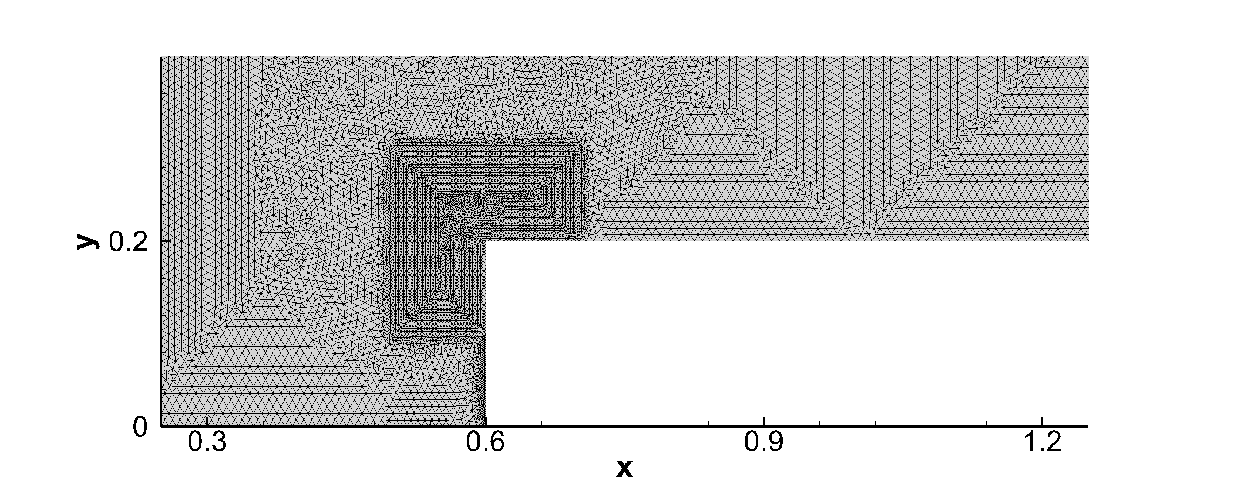}
	\includegraphics[width=0.80\textwidth]{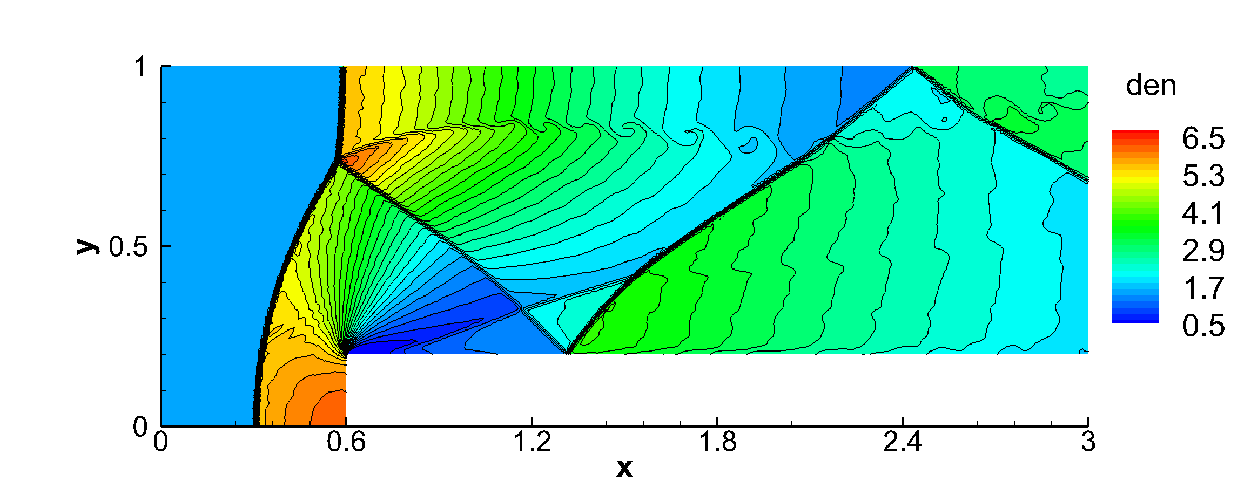}
	\caption{\label{2-forward-step} Mach 3 step problem: Magnified view of the computational mesh (top) and density contours (bottom) obtained by the 4th-order compact GKS with GENO reconstruction. The mesh is a nearly uniform triangular mesh ($h=1/120$) with local refinement to $h=1/300$ near the corner.}
\end{figure}

The computation employs a nearly uniform triangular mesh with a baseline cell size of $h=1/120$, which is locally refined to $h=1/300$ near the corner (Fig. \ref{2-forward-step}, top). The resulting density distribution at $t=4.0$ is displayed. The 4th-order compact GKS, equipped with the GENO reconstruction, produces a high-resolution solution that resolves fine physical instabilities along the slip line. Moreover, the scheme demonstrates robust shock-capturing capabilities and remains free of noticeable oscillations near shock waves.

\subsection{2-D viscous shock tube problem}

Viscous shock tube problem is a complicated flow phenomena associated with shock wave and boundary layer interaction.
This problem requires not only the robustness of the scheme, but also the accuracy of the numerical method. The flow is bounded in a unit square cavity.
The computational domain is set as $[0,1] \times [0,0.5]$ and a symmetrical boundary condition is used on the top boundary.
The non-slip and adiabatic wall conditions are imposed on other boundaries. The initial condition is
\begin{equation*}
(\rho,U,V,p) = \begin{cases}
(120,0,0,120/\gamma),  0\leq x<0.5,\\
(1.2,0,0,1.2/\gamma),  0.5\leq x\leq1.
\end{cases}
\end{equation*}
The viscosity coefficient is $\mu=0.005$ with a corresponding Reynolds number $Re=200$.
The Prandtl number in the current computation is set to be $P_r=1$. The computation time is $t=1.0$.

\begin{figure}[!htb]
\centering
\includegraphics[width=0.495\textwidth]{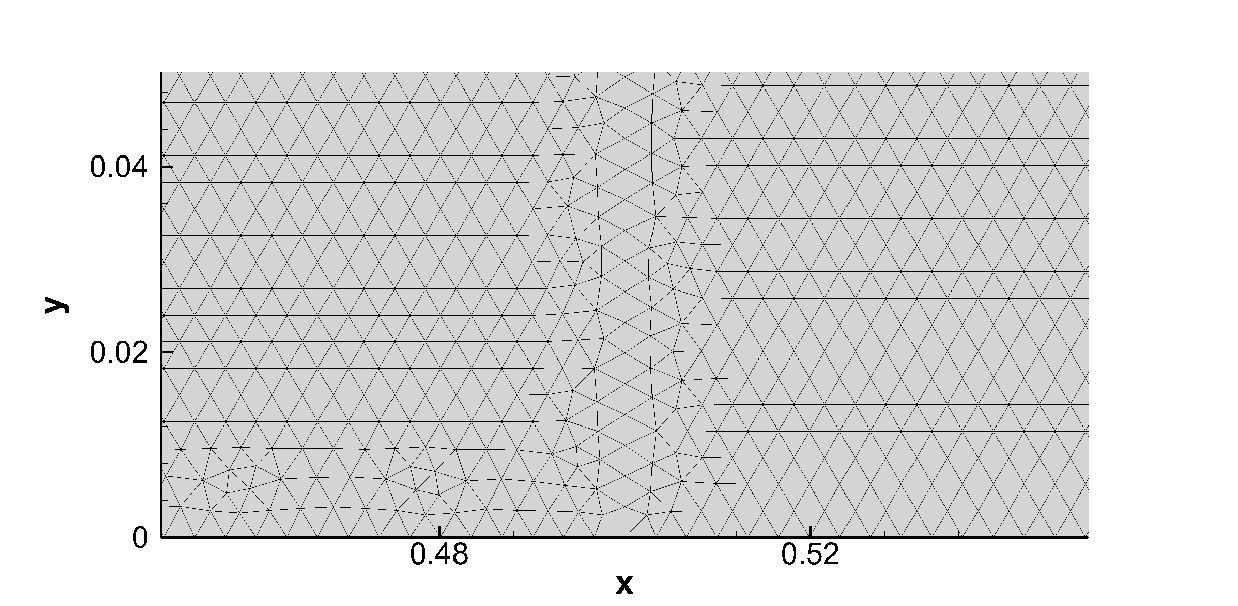}
\includegraphics[width=0.495\textwidth]{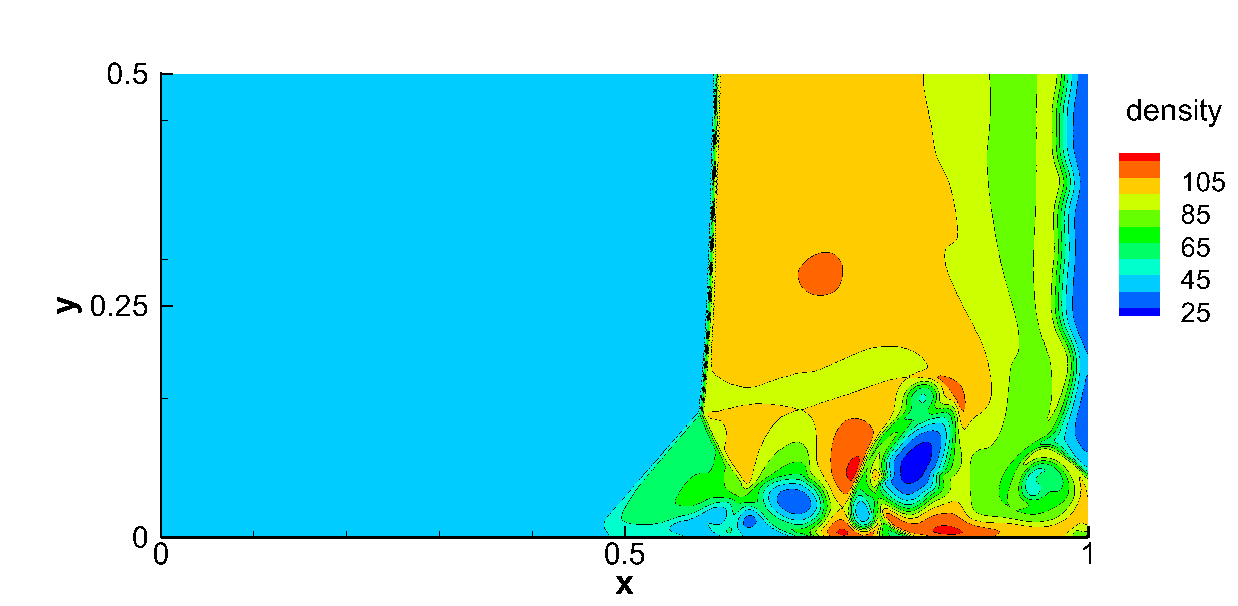}
\caption{\label{2d-vis-shock-tube-1} Viscous shock tube problem: Magnified view of the uniform triangular mesh (cell size $h=1/300$) and density contours at $t=1$, obtained using the 4th-order compact GKS with GENO reconstruction. }
\end{figure}

\begin{figure}[!htb]
\centering
\includegraphics[width=0.495\textwidth]{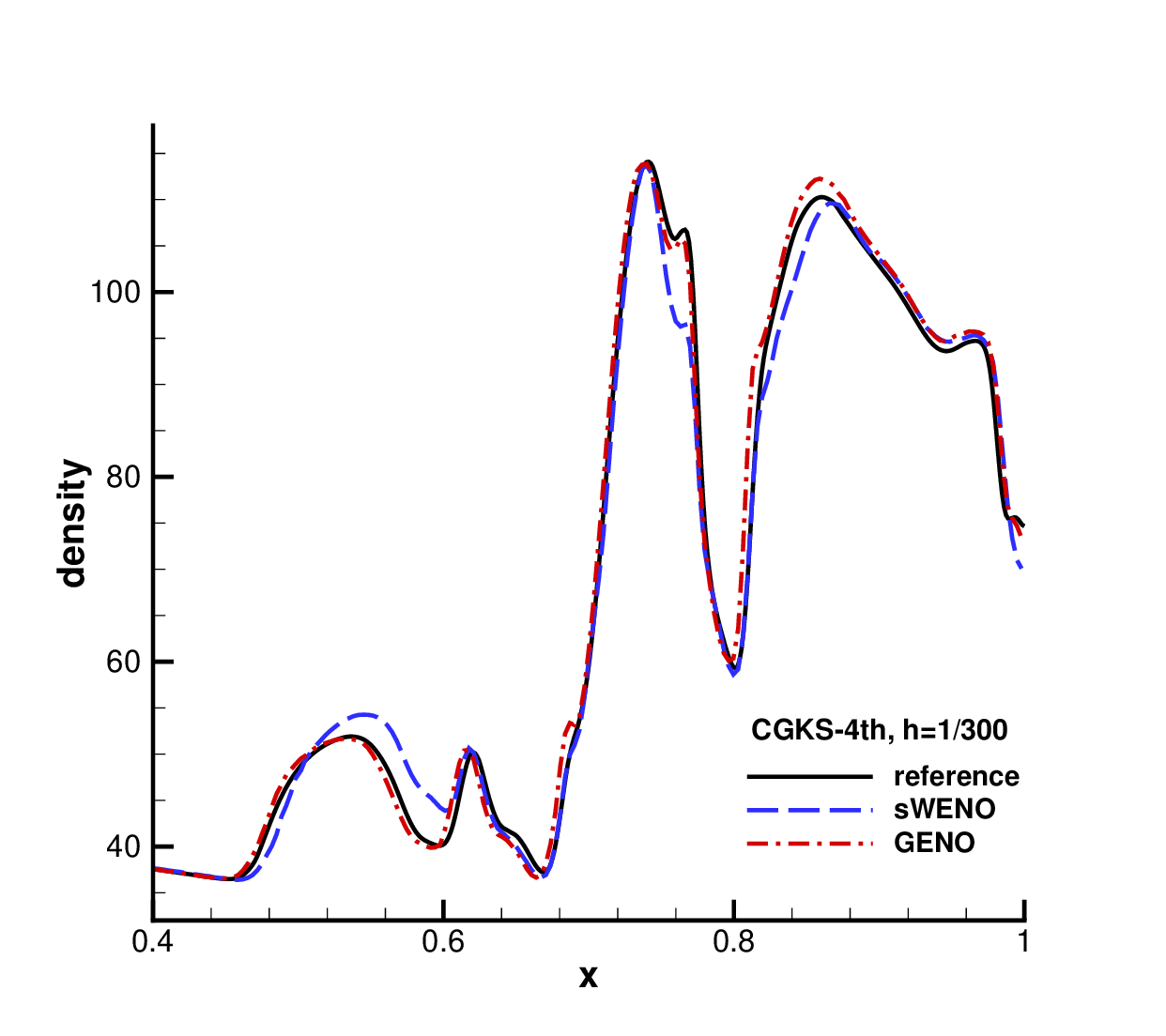}
\includegraphics[width=0.495\textwidth]{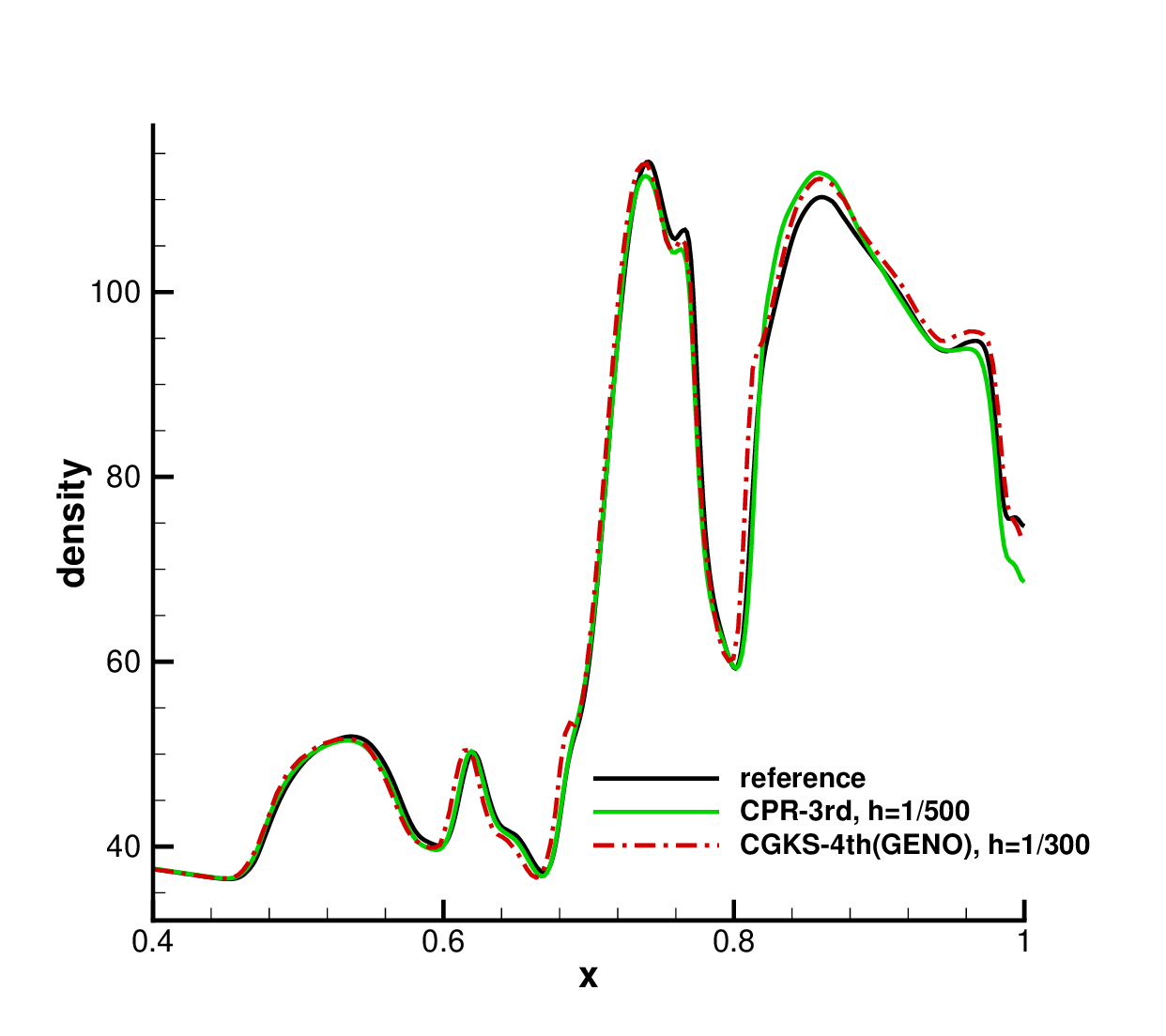}
\caption{\label{2d-vis-shock-tube-2} Viscous shock tube problem: Comparative performance of the 4th-order compact GKS ($h=1/300$) and the 3rd-order CPR ($h=1/500$), with the GKS distribution function used to compute the numerical flux for the CPR. }
\end{figure}

Fig. \ref{2d-vis-shock-tube-1} presents an amplified view of the uniform triangular mesh ($h=1/300$) employed in the computations, along with density contours at $t=1$ obtained using the 4th-order compact GKS with GENO reconstruction.
Quantitative comparisons of the results are provided in Fig. \ref{2d-vis-shock-tube-2}.
The left plot shows results obtained with the compact GKS using both sWENO and GENO reconstructions on the mesh with cell size $h=1/300$. The results corresponding to GENO reconstruction exhibit improved accuracy near the locations with density extremes.
Furthermore, the right plot of Fig. \ref{2d-vis-shock-tube-2} illustrates results from the 4th-order compact GKS with $h=1/300$ and the 3rd-order CPR with $h=1/500$.
Despite the CPR method utilizing more effective degrees of freedom per cell and being computed on a finer mesh (with a cell edge length 0.6 times that of the CGKS method), the compact GKS with GENO reconstruction yields results of nearly identical accuracy.

\begin{figure}[!htb]
\centering
\includegraphics[width=0.495\textwidth]{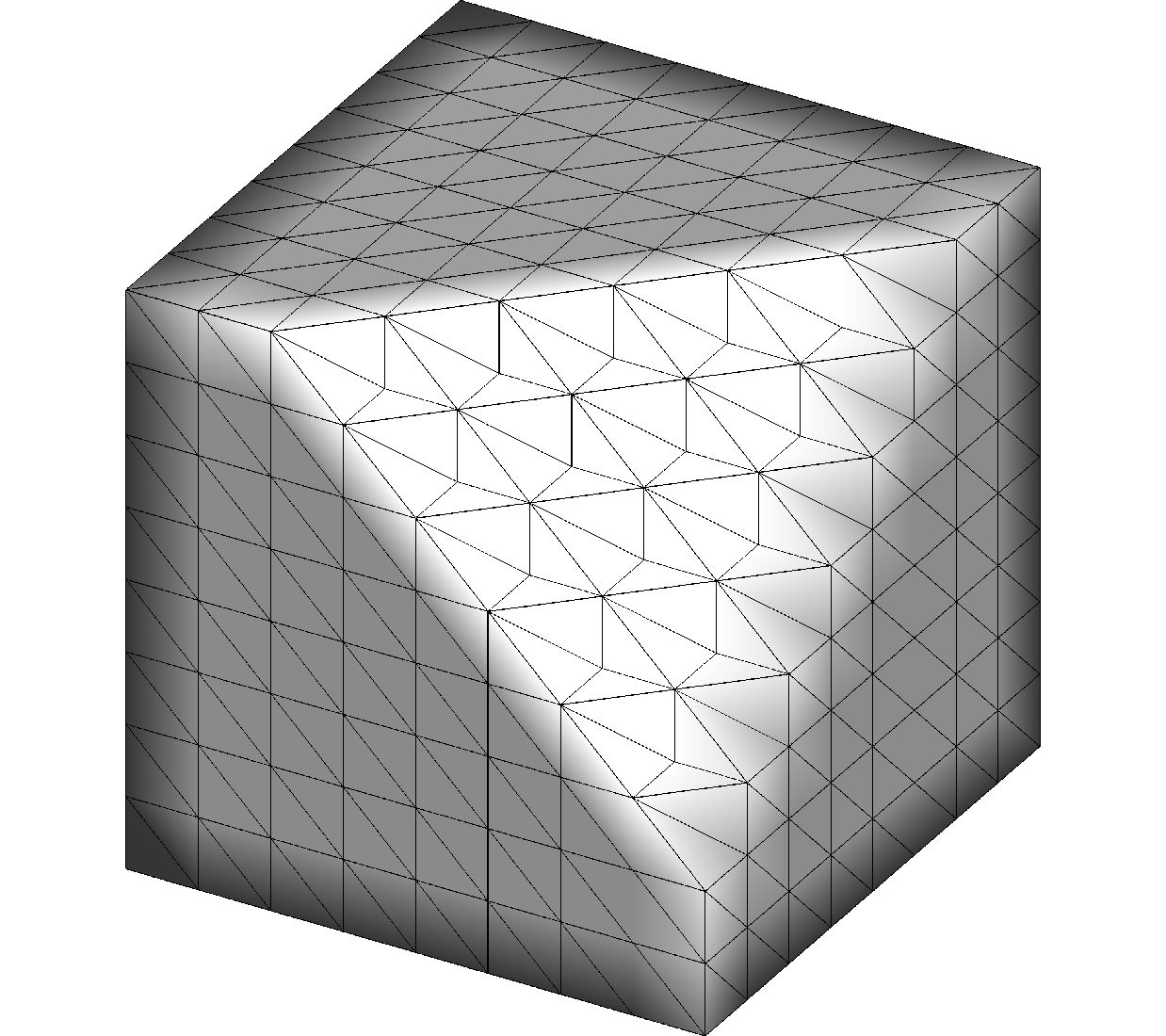}
\includegraphics[width=0.495\textwidth]{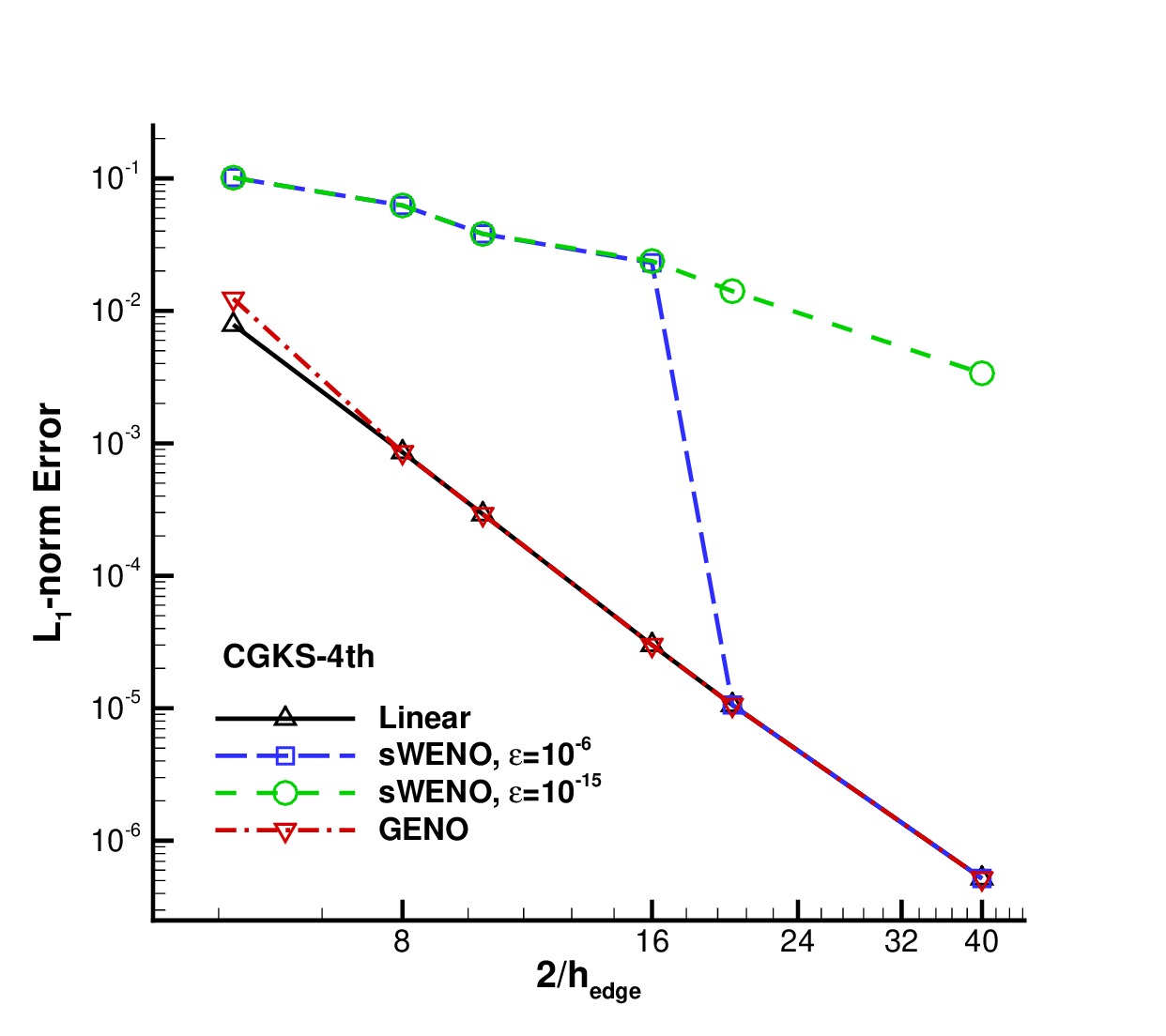}
\caption{\label{3d-Accuracy-dof} 3-D accuracy test: A partial view of the coarse tetrahedral mesh (left) with a cell size of $h=0.25$. The plot on the right presents the error versus mesh refinement at $t=2$. The 4th-order compact GKS with various nonlinear reconstruction methods is employed.}
\end{figure}

\subsection{3-D accuracy test case}
The GENO reconstruction is applied to the 3D high-order compact GKS. First, the accuracy test case is presented.
A 3D density wave propagation problem is used to investigate the accuracy of GENO reconstruction, particularly in the case of coarse mesh.
The initial conditions are given by
\begin{align*}
\begin{split}
\rho(x,y,z)=1+0.2\sin(\pi (x+y+z)),~p(x,y,z)=1, \\
U(x,y,z)=V(x,y,z)=W(x,y,z)=1.
\end{split}
\end{align*}
The computational domain is $[0,2]^3$. Periodic boundary conditions are applied on all boundaries.
A tetrahedral mesh, generated by dividing each hexahedral cell of a Cartesian mesh into six tetrahedra, is used in the computation.

Fig. \ref{3d-Accuracy-dof} shows the convergence of the error with mesh refinement for the 4th-order compact GKS.
GENO reconstruction achieves accuracy consistent with a linear scheme on a relatively coarse mesh with cell size $2/8$.
In comparison, the accuracy of sWENO reconstruction depends on parameters in the nonlinear weights and only reaches the accuracy of a linear reconstruction when the mesh is refined to a size of $2/20$.
In contrast to the one-dimensional case, low-order reconstruction in three dimensions is based on only a linear polynomial, which poses a challenge for WENO-type reconstructions.

\begin{figure}[!htb]
\centering
\includegraphics[width=0.400\textwidth]{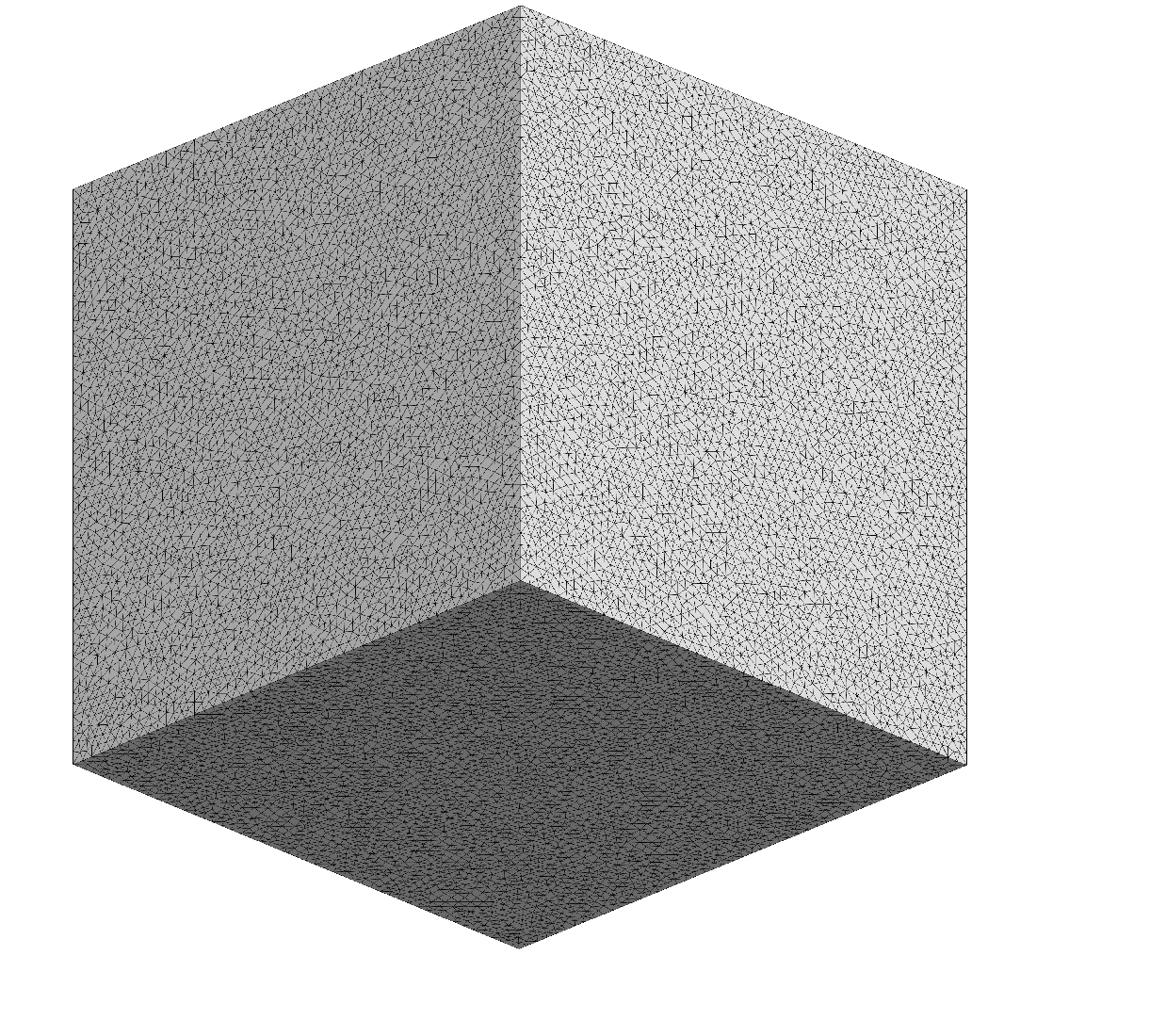}
\includegraphics[width=0.400\textwidth]{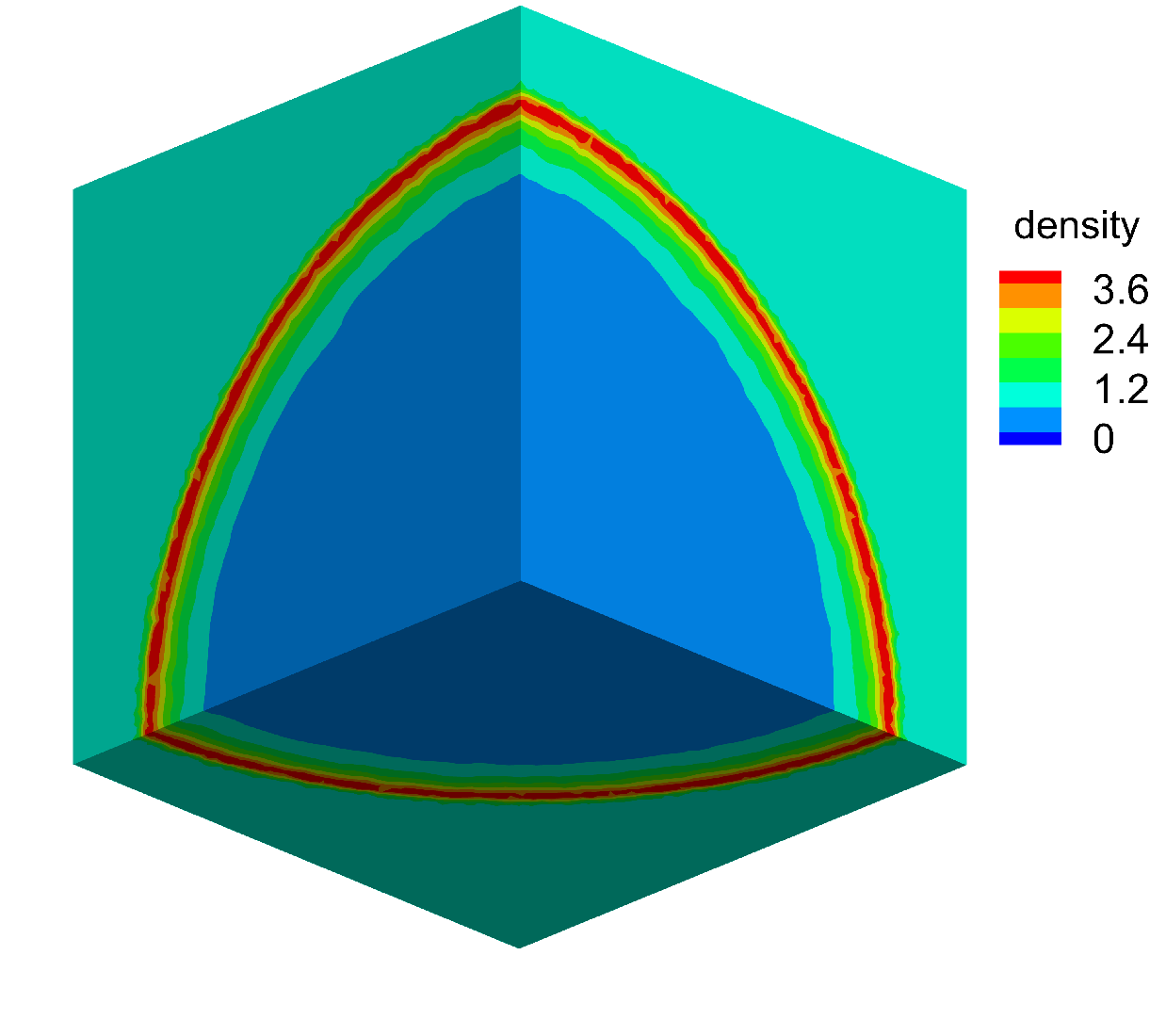}\\
\includegraphics[width=0.495\textwidth]{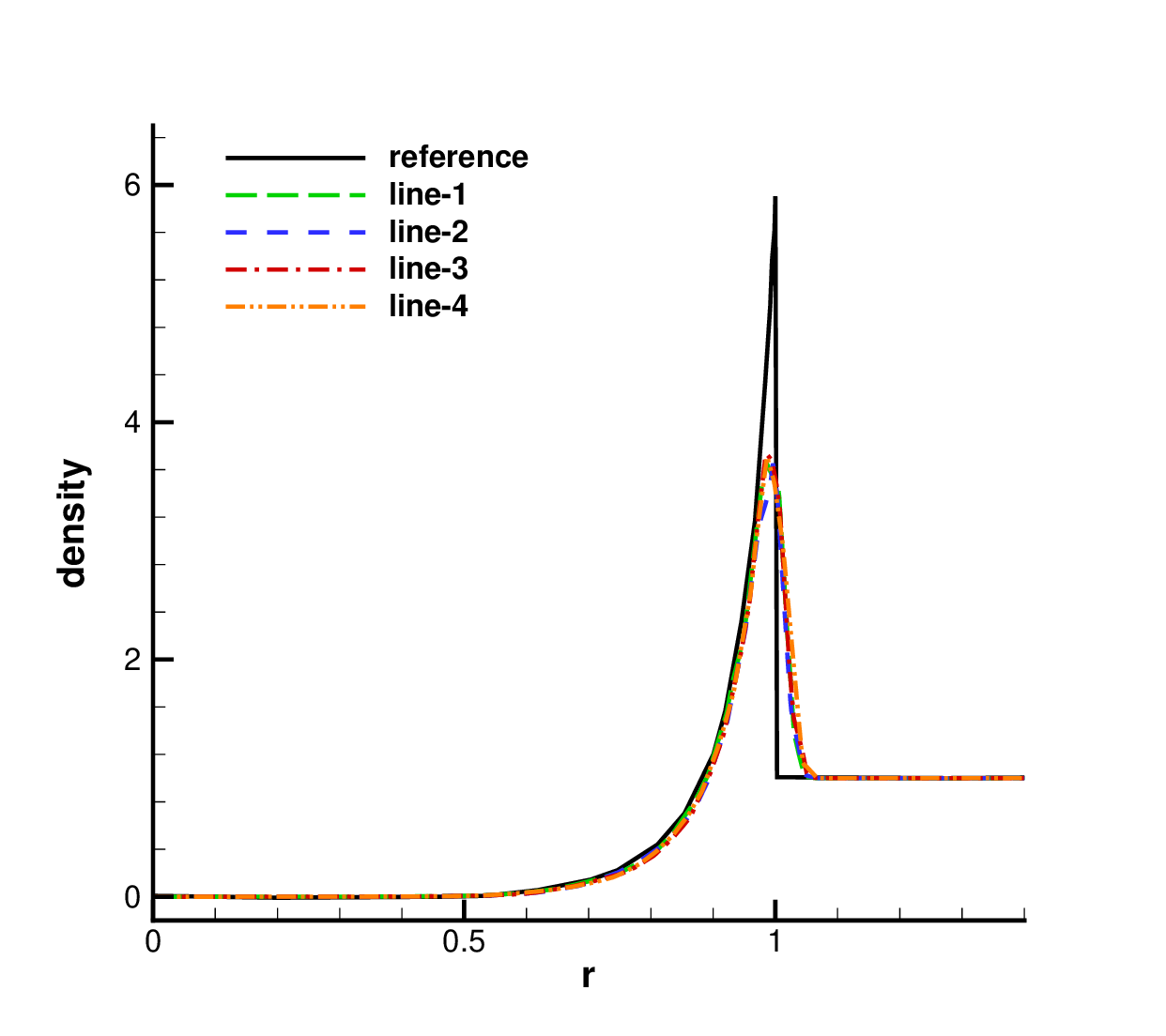}
\includegraphics[width=0.495\textwidth]{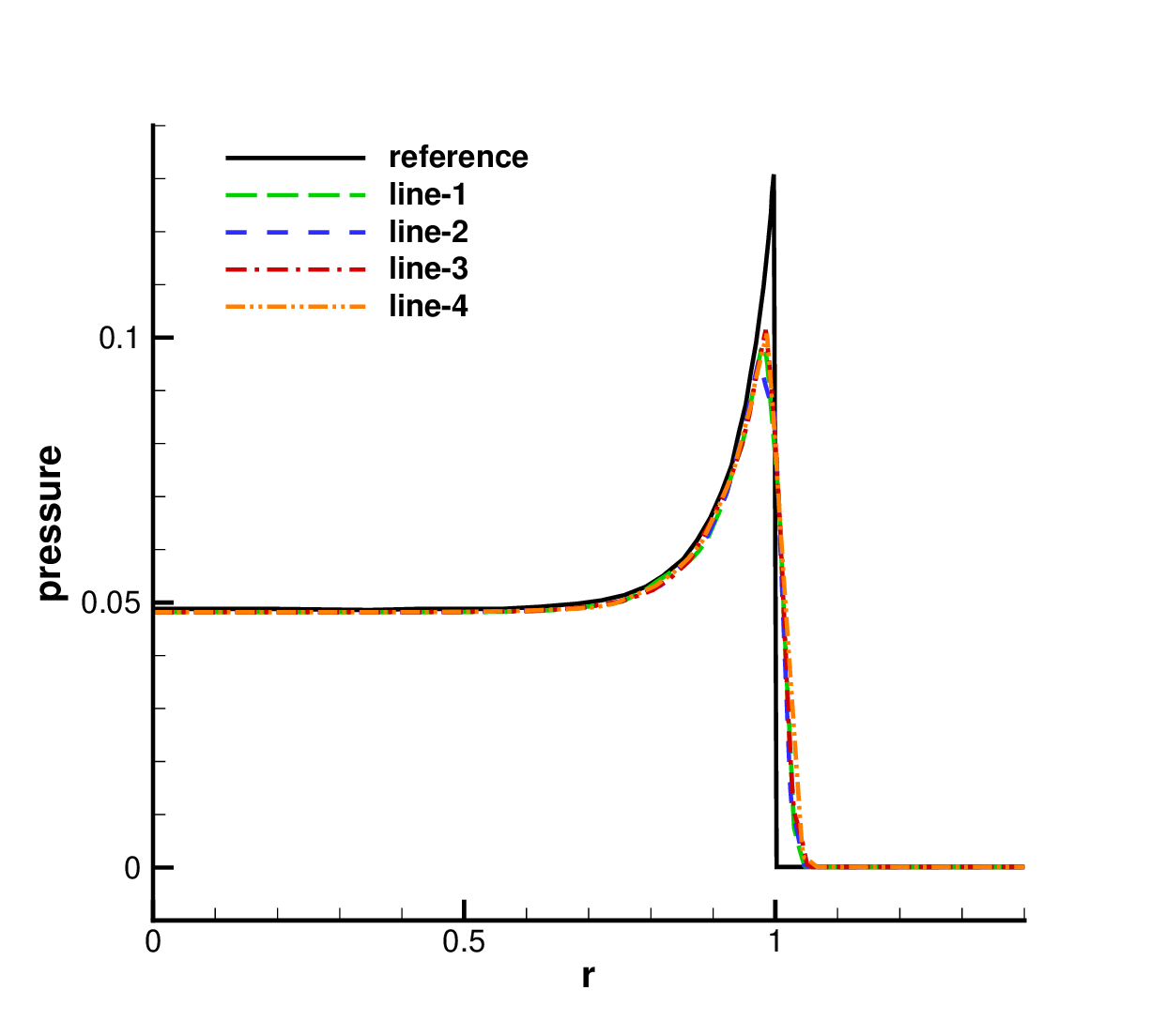}
\caption{\label{1-3d-explosion} Sedov problem: 3-D density distribution (left) and density distribution along different lines (right) at $t=1$. The average side length of cells of the tetrahedral mesh is $h=0.025$.}
\end{figure}

\subsection{3-D explosion and implosion problems}
The 3-D Sedov problem is an explosion case to model blast wave from deposited energy at a singular point \cite{sedov}.
The computational domain is $[0,1.2]^3$. The constant initial conditions with $\rho=1$, $p=1\times10^{-4}$ and $U=V=W=0$ are imposed in the whole domain except the cells containing the origin. The pressure of the cells containing the origin is set as $p=(\gamma-1)\epsilon /V_0$, where $\epsilon=0.106384$ and $V_0$ is the total volume of those cells.
The inviscid wall condition is adopted along the boundaries $x=0$, $y=0$ and $z=0$. The free boundary condition is imposed to the other boundaries.
The cell size is $h=0.025$ defined by the average side length of tetrahedral cells. The computational output time is $t=1.0$.

The 3D density distributions along different lines are shown in Fig. \ref{1-3d-explosion}, where the reference solution is from \cite{sedov}.
Lines 1, 2, and 3 are defined by connecting the origin to the points $(1.2, 0.6, 0.6)$, $(0.6, 1.2, 0.6)$, and $(0.6, 0.6, 1.2)$, respectively. Line 4 connects the origin to $(1.2, 1.2, 1.2)$.
The strong robustness of the 4th-order compact GKS is demonstrated by its use of a large CFL number of $CFL=0.6$ without additional limiting technique.
In addition, the high resolution of the compact GKS for strong shock waves is verified by the high post-shock density peak, and the numerical shock wave remains sharp and spans only two mesh cells.
Compared to the results of the 5th-order finite volume scheme in \cite{charest2015high}, the present 4th-order compact GKS demonstrates better accuracy, with comparable mesh cell sizes used in both computations.

\begin{figure}[!htb]
\centering
\includegraphics[width=0.400\textwidth]{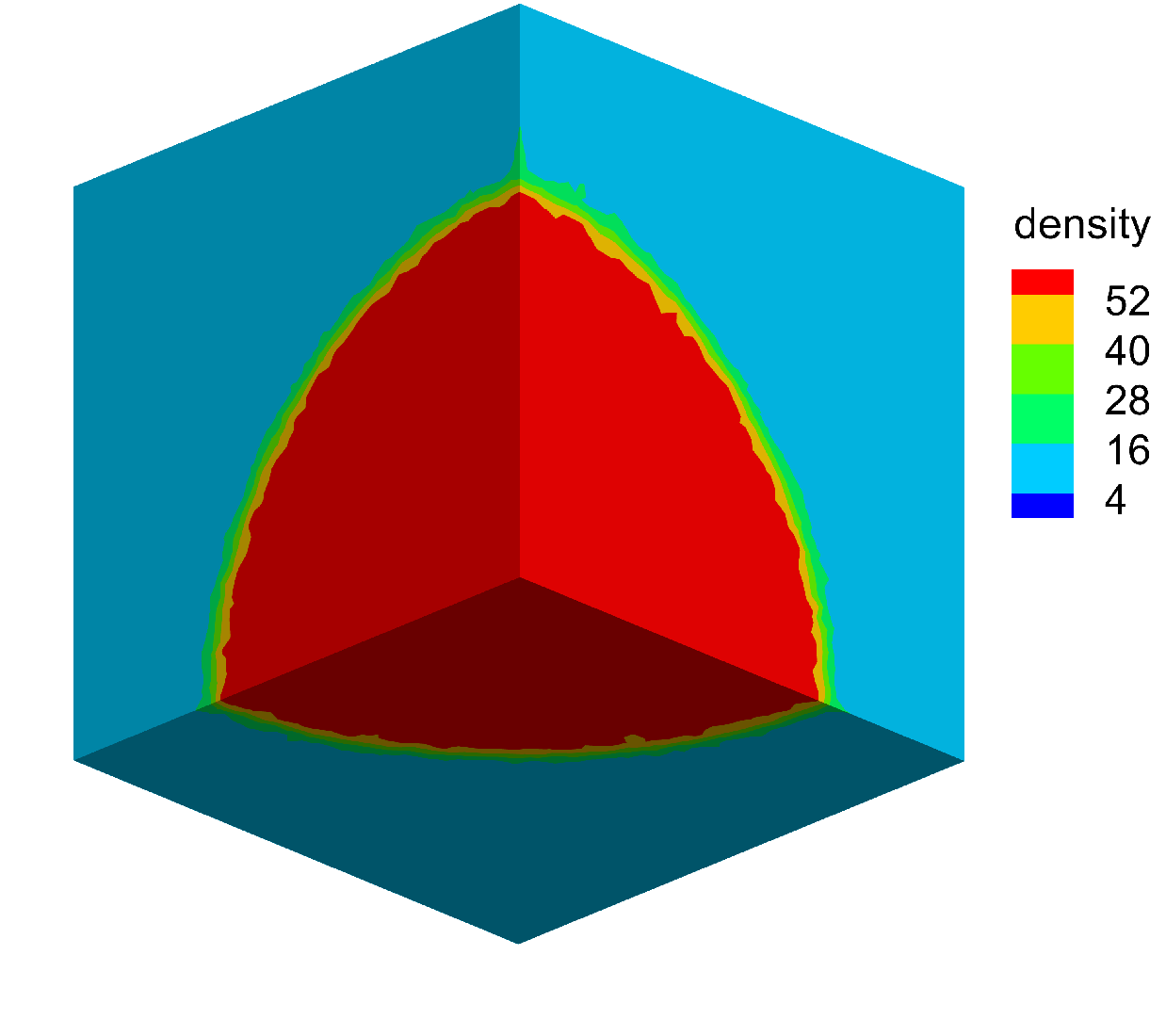}
\includegraphics[width=0.495\textwidth]{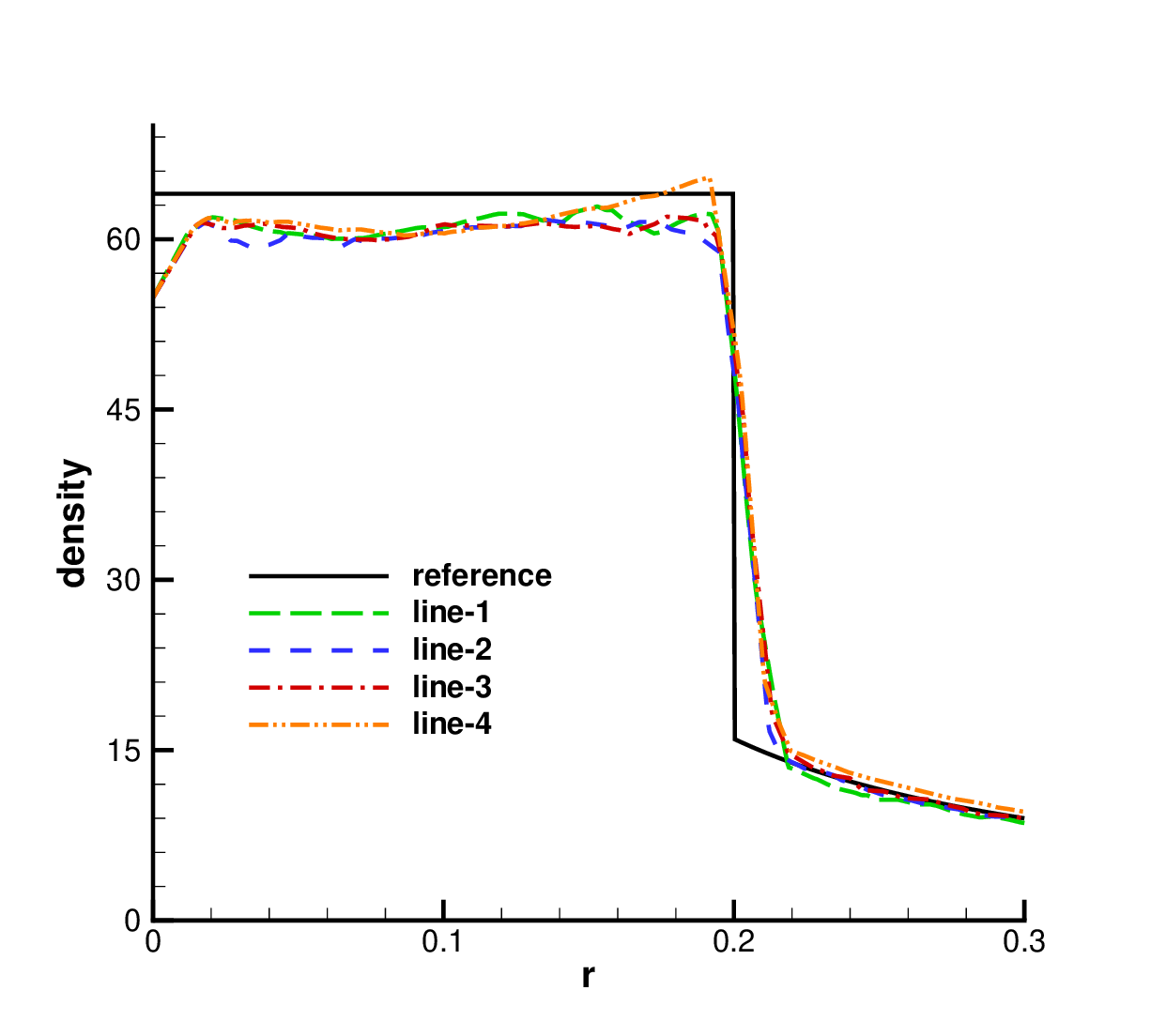}
\caption{\label{1-3d-implosion} Noh problem: 3-D density distribution (left) and density distribution along different lines (right) at $t=0.6$. The average side length of cells of the tetrahedral mesh is $h=0.01$.}
\end{figure}

The 3-D Noh problem is an implosion test to model the gas compression with constant radial velocity towards a spherical center, where a moving strong shock wave is formed \cite{Noh}.
The computational domain is $[0,0.3]^3$. The initial density and pressure are $\rho=1$ and $p=1\times10^{-3}$, and the velocity is $(U,V,W)=(-x,-y,-z)/\sqrt{x^2+y^2+z^2}$. The ratio of the specific heat is $\gamma=5/3$.
The inviscid wall condition is adopted along the boundaries $x=0$, $y=0$ and $z=0$. The supersonic inflow boundary condition is imposed to the other boundaries with the same pressure and velocity as the initial conditions and the analytical solution of density,
\begin{equation*}
\rho = \begin{cases}
64,  r<t/3,\\
(1+t/r)^2,  r>t/3.
\end{cases}
\end{equation*}
The mesh cell size is $1/100$ defined by the average side length of tetrahedral cells.
The final computational output time is $t=0.6$.

The 3-D density distributions along different lines are shown in Fig. \ref{1-3d-implosion}.
Lines 1, 2, and 3 are defined by connecting the origin to the points $(0.3, 0.15, 0.15)$, $(0.15, 0.3, 0.15)$, and $(0.15, 0.15, 0.3)$, respectively. Line 4 connects the origin to $(0.3, 0.3, 0.3)$.
Again, in this test case, a large CFL number of $CFL=0.6$ without additional limiting technique is used.
The present 4ht-order compact GKS with GENO reconstruction yields sufficiently accurate results, outperforming several well-developed high-order schemes with orders up to $8$ reported in \cite{johnsen2010assessment}, even though the computations in \cite{johnsen2010assessment} were performed on a mesh with a cell size of $0.002$.

\subsection{3-D Taylor-Green vortex flow}
The Taylor-Green vortex flow is a canonical benchmark for evaluating the performance of high-order schemes on viscous flows.
Evolving from a smooth initial condition, this flow progressively develops a broadband of fine-scale structures, making it an ideal test for a scheme's ability to capture the complex flow structures in turbulence.
The initial condition is set as
\begin{equation*}
\begin{split}
&U= U_0\mathrm{sin}x \mathrm{cos}y \mathrm{cos}z,\\
&V=-U_0\mathrm{cos}x \mathrm{sin}y \mathrm{cos}z,\\
&W=0,\\
&\rho=\rho_0+\frac{\gamma Ma_{ref}^2}{16}(\mathrm{cos}2x+\mathrm{cos}2y)(\mathrm{cos}2z+2),
\end{split}
\end{equation*}
where $U_0=1$, and the initial temperature is $T_0=1$.
The computational domain is a periodic cube defined by $\Omega_0=[-\pi,\pi]^3$.
The reference Mach number is set to $Ma_{ref}=0.1$, and the Reynolds number is $Re=1600$.
These are defined by $Ma_{ref}=U_0/\sqrt{\gamma p_0/\rho_0}$ and $Re=\rho_0 U_0 L_0/\mu$, respectively,
where the characteristic length is $L_0=2\pi$ and $\mu$ is the dynamic viscosity.

The averaged kinetic energy $E_k$ and its dissipation rate $\varepsilon(E_k)$ are used to quantitatively assess the accuracy of the numerical results.
The averaged kinetic energy is defined as
\begin{align*}
E_k=\frac{1}{\rho_0|\Omega|}\int_{\Omega} \frac{1}{2}\rho \mathbf{U}\cdot \mathbf{U}\mathrm{d}V.
\end{align*}
The dissipation rate of kinetic energy is given by
\begin{align*}
\varepsilon(E_k)=-\frac{\mathrm{d} E_k}{\mathrm{d} t}.
\end{align*}

\begin{figure}[!htb]
\centering
\includegraphics[width=0.475\textwidth]{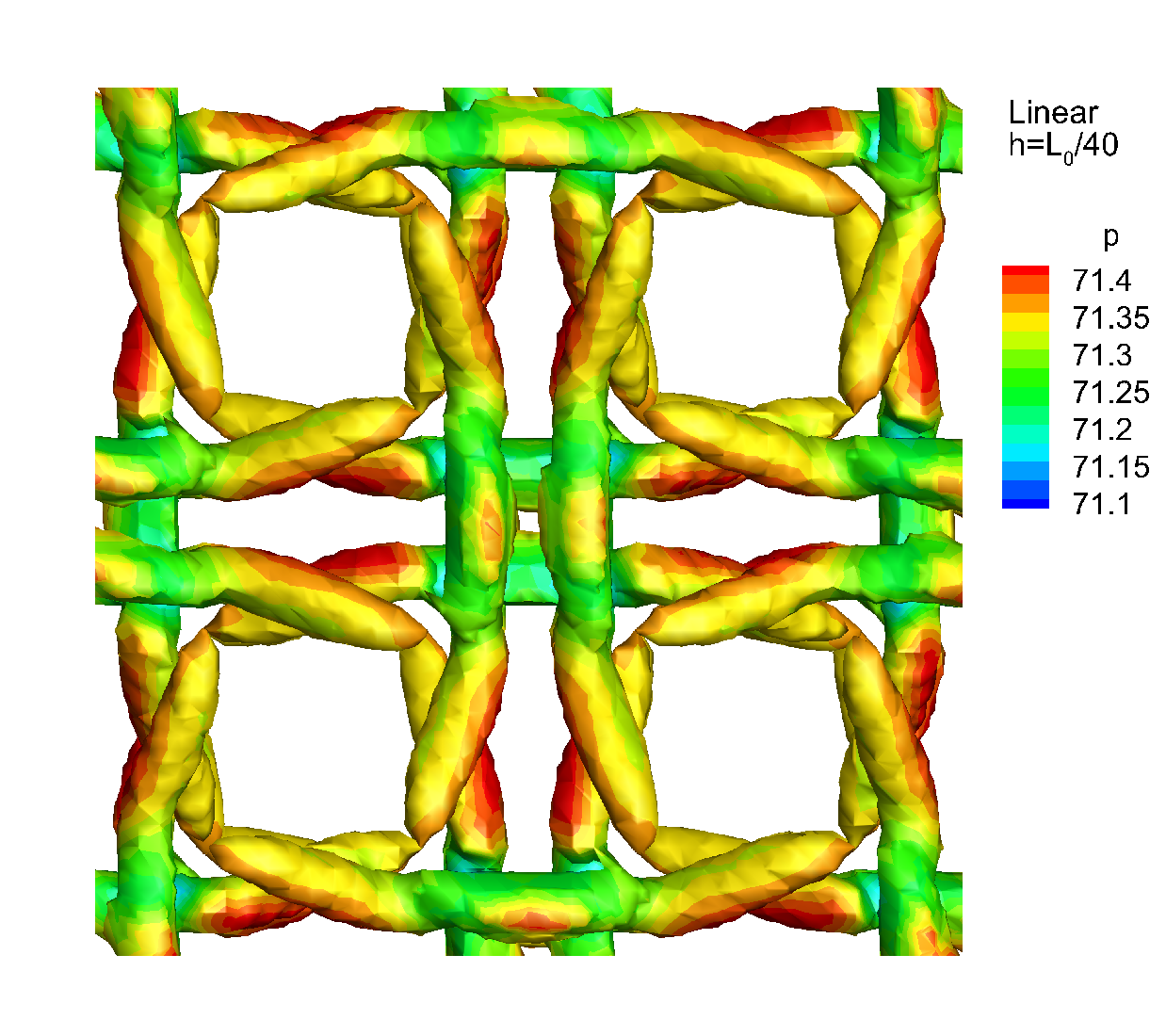}
\includegraphics[width=0.475\textwidth]{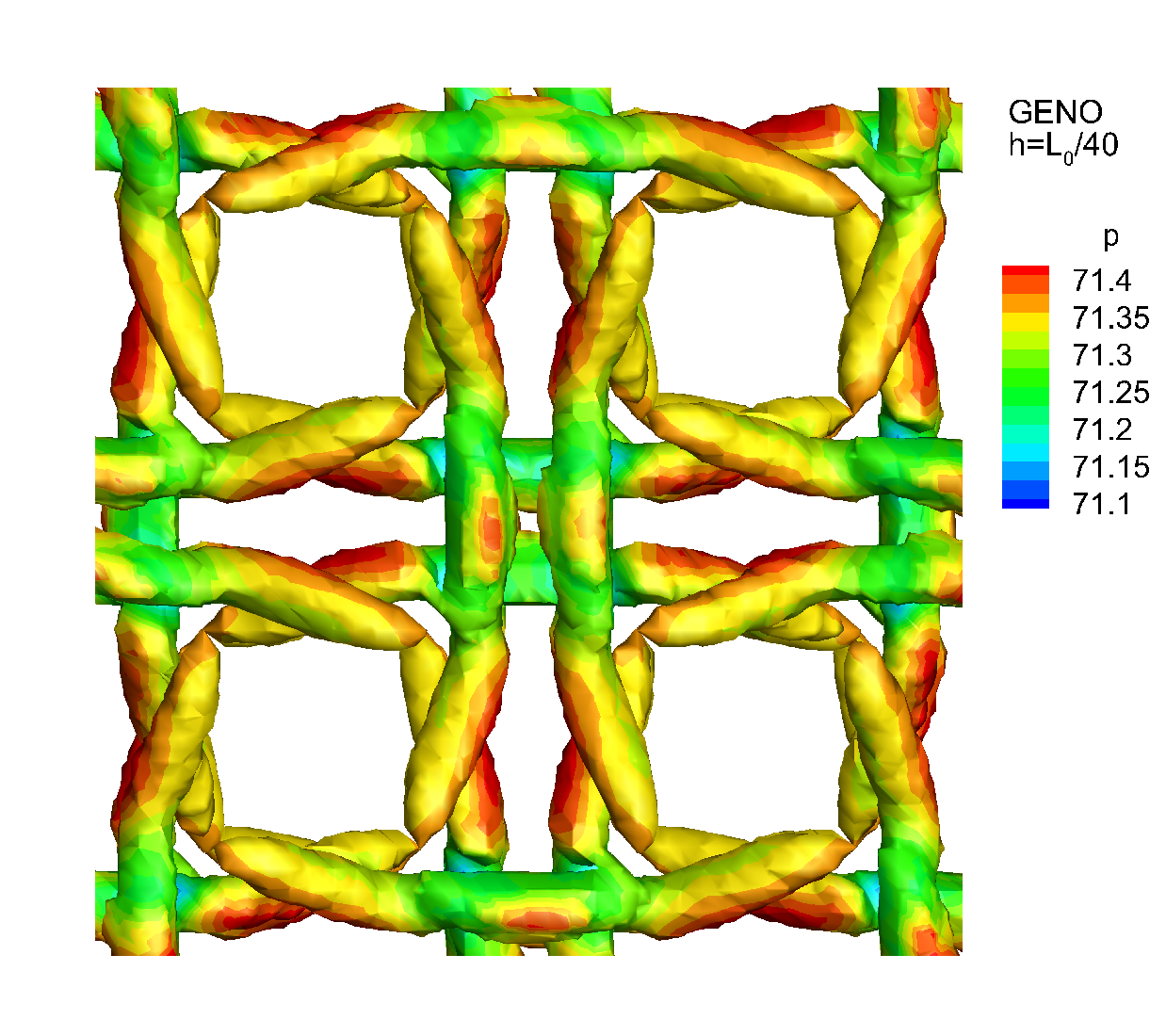}
\caption{\label{3d-tgv-1} Taylor-Green vortex flow at $Re=1600$: Q-criterion isosurfaces ($Q=1$) colored by pressure at $t =5$. The left and right panels correspond to the 4th-order compact GKS results using linear and GENO reconstructions, respectively.}
\end{figure}

\begin{figure}[!htb]
\centering
\includegraphics[width=0.475\textwidth]{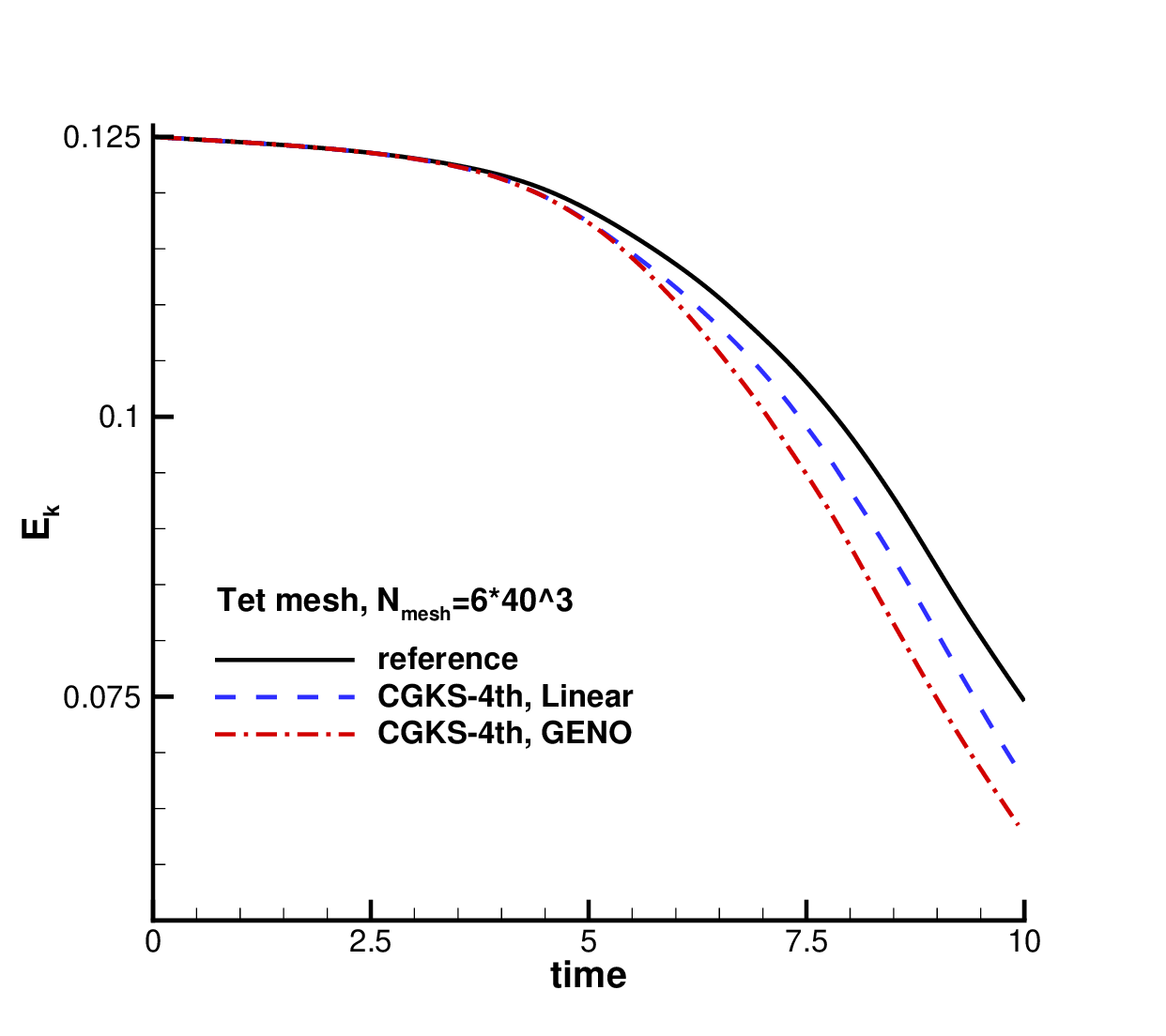}
\includegraphics[width=0.475\textwidth]{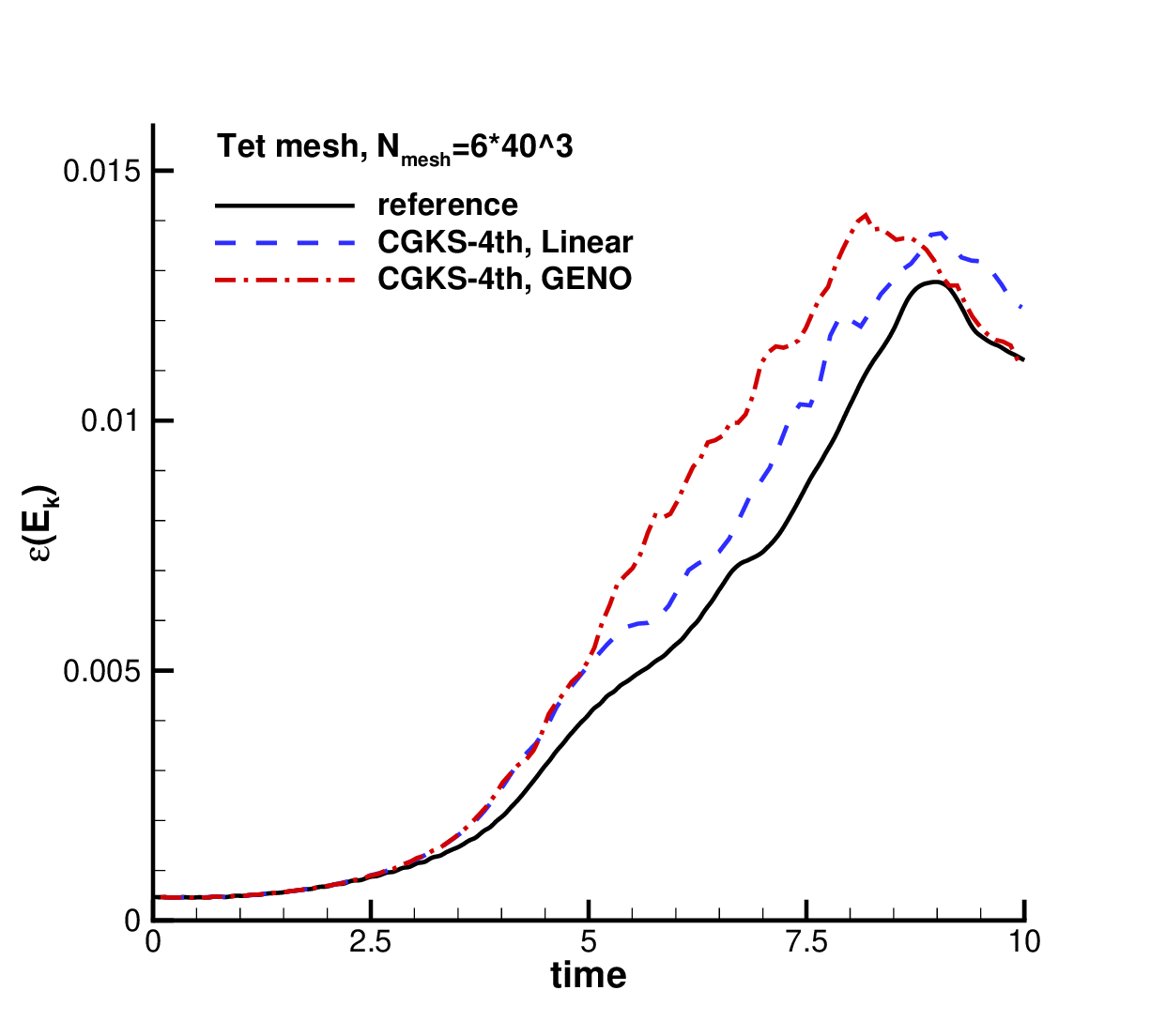}
\caption{\label{3d-tgv-2} Taylor-Green vortex flow at $Re=1600$: time history of kinetic energy and kinetic energy dissipation rate. The 4th-order compact GKS is used.}
\end{figure}

To assess the performance of GENO reconstruction when the mesh resolution is insufficient to fully resolve flow structures in vortex dynamics evolution, a coarse mesh is employed in the computation.
The mesh, consisting of $6 \times 40^3$ cells, was generated using the same procedure as the meshes for the 3-D accuracy tests shown in Fig. \ref{3d-Accuracy-dof}.
Fig. \ref{3d-tgv-1} illustrates the vortex structures identified by the Q-criterion isosurfaces at $t = 5$.
The left and right figures display the results obtained using the 4th-order compact GKS with linear reconstruction and GENO reconstruction, respectively.
As can be seen, the linear and GENO reconstructions yield nearly identical results.
Fig. \ref{3d-tgv-2} shows the time history of $E_k$ and $\varepsilon(E_k)$. The results obtained with GENO reconstruction deviate only slightly from those obtained with linear reconstruction, and the level of deviation is acceptable on the present coarse mesh.

\section{Conclusion}

This study introduces GENO, a novel formulation for constructing high-order nonlinear reconstructions that enables idealized adaptive transitions between high-order linear and robust lower-order reconstructions when solving problems ranging from smooth to discontinuous solutions.
GENO extends the ENO paradigm by generalizing its smoothness-based adaptive stencil selection and reconstruction approach. The method adaptively employs high-order linear reconstructions in smooth regions and ENO-property-possessing lower-order reconstructions in discontinuous regions. A linear-reconstruction-preserving path function ensures smooth transitions between these reconstruction types in intermediate regions.
This approach yields superior accuracy compared to WENO-type schemes for under-resolved flow features while maintaining robust, non-oscillatory behavior when handling discontinuous solutions with strong shocks. The direct coupling between high-order and low-order reconstructions in GENO significantly simplifies the construction of high-order nonlinear schemes, particularly beneficial for very high orders and unstructured mesh implementations.
Comparative analysis demonstrates that GENO outperforms WENO-type methods in managing smooth-to-discontinuous transitions, achieving reduced deviation from linearity under challenging smoothness conditions. Comprehensive numerical validation across one-, two-, and three-dimensional unstructured meshes confirms that GENO, when implemented within a high-order compact GKS framework, delivers both enhanced accuracy and robust shock-capturing capabilities. Additionally, its successful application to non-compact schemes validates GENO's consistency while highlighting the inherent advantages of compact formulations.

\section*{Acknowledgments}

The current research is supported by National Key R\&D Program of China (Grant Nos. 2022YFA1004500), National Science Foundation of China (12172316, 92371107), and Hong Kong research grant council (16301222, 16208324).

\bibliographystyle{ieeetr}
\bibliography{compact-GENO}

\end{document}